\documentclass[a4paper, 12pt]{article}
\usepackage[russian]{babel}
\usepackage{amsmath}
\usepackage{amsfonts}
\newcommand{\mx}{\operatorname{\mathfrak M}}
\newcommand{\cmx}{\operatorname{\mathfrak C}}
\newcommand{\mn}{\operatorname{\mathfrak m}}
\newcommand{\cmn}{\operatorname{\mathfrak c}}
 
\newcommand{\e}{\operatorname{\mathfrak e}}
\newcommand{\pq}{\operatorname{\mathfrak q}}
\newcommand{\q}{\operatorname{q^*}}
\newcommand{\p}{\operatorname{p^*}}
\newcommand{\s}{\operatorname{\mathcal s}}
\newcommand{\card}{\operatorname{card}}
\newcommand{\diam}{\operatorname{diam}}
\newcommand{\supvrai}{\operatornamewithlimits{sup\,vrai}}
\newcommand{\N}{\mathbb N}
\newcommand{\Z}{\mathbb Z}
\newcommand{\R}{\mathbb R}
\newcommand{\Nu}{\mathcal N}

\newcommand{\J}{\mathcal J}
\newcommand{\mes}{\operatorname{mes}}

\begin{document}

\author{ С. Н. Кудрявцев }
\title{Теорема типа Литтлвуда-Пэли для ортопроекторов
на взаимно ортогональные подпространства
кусочно-полиномиальных функций
и следствие из не\"е}
\date{}\maketitle
\begin{abstract}
В статье доказано утверждение, представляющее собой
аналог теоремы Литтлвуда-Пэли для ортопроекторов на взаимно
ортогональные подпространства кусочно-полиномиальных функций на кубе $ I^d. $
Исходя из него установлена оценка сверху нормы функций в
$ L_p(I^d) $ через соответствующие нормы проекций на подпространства
кусочно-полиномиальных функций нескольких переменных.
С помощью этих соотношений получены верхние оценки колмогоровских
поперечников классов Бесова (непериодических) функций,
удовлетворяющих смешанным условиям Гельдера.
\end{abstract}

Ключевые слова: ортопроектор, взаимно ортогональные подпространства,
кусочно-полиномиальные функции, теорема Литтлвуда-Пэли, поперечник

\bigskip

\centerline{Введение}

Как известно (см., например, [1], [2] и др. работы), важное значение для
вывода порядковых оценок
точности приближения в $L_p$ способами,
основанными на применении кратных тригонометрических рядов, классов периодических функций
нескольких переменных с условиями на смешанные производные (разности) имеет
теорема Литтлвуда-Пэли для кратных рядов Фурье и следствия из нее.
По поводу теоремы Литтлвуда-Пэли для кратных рядов Фурье
см., например, [3, п. 1.5.2], [4] и приведенную там литературу.
Для получения соответствующих оценок точности приближения классов
(непериодических) функций, Заданных на кубе $ I^d, $ подчиненных условиям на смешанные разности,
полезно иметь аналоги этих результатов для средств приближения таких классов
непериодических функций.
С этой целью в работе доказан аналог теоремы Литтлвуда-Пэли для
семейств ортопроекторов $ \{\mathcal E_\kappa^{d,l}, \kappa \in \Z_+^d \},
d \in \N, l \in \Z_+^d $ (см. п. 2.2.)
на взаимно ортогональные подпространства кусочно-полиномиальных функций
на кубе $ I^d, $ построенных в [5] (см. также ниже п. 1.4.),
а именно, показано, что при $ d \in \N, l \in \Z_+^d, 1 < p < \infty $
для любой функции $ f \in L_p(I^d) $ выполняются неравенства
\begin{equation*} \tag{1}
c_{1} \| f \|_{L_p(I^d)} \le
( \int_{I^d} (\sum_{\kappa \in \Z_+^d} ((\mathcal E_\kappa^{d,l} f)(x))^2 )^{p/2}
dx)^{1/p} \le c_{2} \| f \|_{L_p(I^d)}
\end{equation*}
с некоторыми константами $ c_1, c_2 > 0, $ зависящими от $ d,l,p.$

Исходя из (1), установлена оценка сверху нормы функции $ f $ в
$ L_p(I^d) $ через соответствующие нормы проекций $ \mathcal E_\kappa^{d,l} f,
\kappa \in \Z_+^d, $ на подпространства кусочно-полиномиальных функций,
которая имеет вид
\begin{equation*} \tag{2}
\| f \|_{L_p(I^d)} \le c_{3} (\sum_{\kappa \in \Z_+^d }
\| \mathcal E_\kappa^{d,l} f \|_{L_p(I^d)}^{\p} )^{1/\p},
\end{equation*}
(ср. с аналогом для операторов взятия частных сумм ряда Фурье в [2]), где $ \p = \min(2,p). $
В качестве иллюстрации применения с помощью (2) получена оценка сверху
колмогоровского $n$-поперечника в $ L_q(I^d) $ множества
$\mathcal S_{p,\theta}^\alpha \mathcal B(I^d) $ -- единичного шара
относительно полунормы, определяющей пространство Бесова функций на кубе $ I^d, $
удовлетворяющих смешанным условиям Гельдера (см. п. 3.1.).
Эта оценка при некоторых условиях на
$ d \in \N, \alpha \in \R_+^d, 1 \le \theta \le \infty, 1 \le p < \infty,
1 \le q < \infty $ такова:
\begin{multline*} \tag{3}
d_n(\mathcal S_{p,\theta}^\alpha \mathcal B(I^d), L_q(I^d))\\ \le c_{4}
n^{-\mn(\alpha -(p^{-1} -\pq^{-1})_+ \e)}
(\log n)^{(\mn(\alpha -(p^{-1} -\pq^{-1})_+ \e) +(1/\pq -1/\theta)_+)
(\cmn -1)},
\end{multline*}
где $ \mn(x) = \min_{j=1,\ldots,d} x_j, x \in \R^d, \cmn = \card
\{j =1, \ldots, d: \alpha_j = \mn(\alpha) \}, \e = (1, \ldots,1),
\pq = \min(2, \max(p, q)). $ Порядковая оценка (3) при $ q \ge 2 $
точна (см. (3.3.19')). Для сравнения (3) со случаем периодических
функций см. [2].

Работа состоит из введения и трех параграфов. В \S 1 приведены
предварительные сведения об ортопроекторах на подпространства кусочно-полиномиальных
функций и кратных рядах, рассматриваемых в работе. В п. 2.2. \S 2 на
основании сведений из \S 1  и п. 2.1.
устанавливается справедливость (1) и выводится (2). В \S 3 с применением
(2) доказывается (3) (см. п. 3.3.). Перейдем к точным формулировкам
и доказательствам.
\bigskip

\centerline{\S 1. Операторы проектирования на подпространства }
\centerline{кусочно-полиномиальных функций}
\bigskip

 1.1. В этом пункте вводятся обозначения, используемые
в настоящей работе.

Для $ d \in \N $ через $ \Z_+^d $ обозначим множество
$$
\Z_+^d =\{\lambda =(\lambda_1,  \ldots, \lambda_d) \in \Z^d:
\lambda_j \ge0, j=1, \ldots, d\}.
$$
 Обозначим также при  $ d \in \N $ для $ l \in \Z_+^d $ через
$ \Z_+^d(l) $ множество
 $$
 \Z_+^d(l) =\{ \lambda  \in \Z_+^d: \lambda_j \le l_j, j=1, \ldots,  d\}.
 $$

Для $  d \in  \N, l \in \Z_+^d $ через $ \mathcal P^{d,l} $ будем
обозначать пространство вещественных  полиномов, состоящее из всех
функций $ f: \R^d \mapsto \R $ вида
$$
f(x) =\sum_{\lambda  \in \Z_+^d(l)}a_{\lambda}\cdot x^{\lambda},
x\in \R^d,
$$
где  $   a_{\lambda}   \in   \R,  x^{\lambda} =x_1^{\lambda_1}
\ldots x_d^{\lambda_d}, \lambda \in \Z_+^d(l). $

При $ d \in \N, l \in \Z_+^d  $   для области $ D\subset\R^d $
через $\mathcal   P^{d,l}(D) $  обозначим пространство функций $
f, $ определенных в $  D, $   для каждой из  которых существует
полином $ g\in\mathcal P^{d,l} $
 такой, что  сужение $ g\mid_D = f.$

При $ n \in \N, 1 \le p \le \infty $ через $ l_p^n $ обозначим
пространство $ \R^n $ с фиксированной в нем нормой
$$
\|x\|_{l_p^n} = \begin{cases} (\sum_{j=1}^n
|x_j|^p)^{1/p}, \text{ при } p < \infty; \\
\max_{j=1, \ldots,n} |x_j|, \text{ при } p = \infty,
\end{cases} x \in \R^n.
$$

 Для измеримого по Лебегу множества $ D \subset \R^d$
при $ 1\le p\le\infty $ через  $  L_p(D),$  как обычно,
обозначается пространство  всех  вещественных измеримых на $ D $
функций $f,$ для которых определена норма
$$
\|f\|_{L_p(D)}= \begin{cases} (\int_D |f(x)|^p
dx)^{1/p}, 1 \le p<\infty;\\
\supvrai_{x \in D}|f(x)|, p=\infty. \end{cases}
$$

Введем еще следующие обозначения.

Для $ x,y \in \R^d $ положим $ xy =x \cdot y = (x_1 y_1, \ldots,
x_d y_d), $ а для $ x \in \R^d $ и $ A \subset \R^d $ определим
$$
x A = x \cdot A = \{xy: y \in A\}.
$$

Для $ x \in \R^d: x_j \ne 0, $ при $ j=1,\ldots,d,$ положим $
x^{-1} = (x_1^{-1},\ldots,x_d^{-1}). $

При $ d \in \N $ для $ x,y \in \R^d $ будем обозначать
$$
(x,y) = \sum_{j=1}^d x_j y_j.
$$

Обозначим через $ \R_+^d $ множество $ x \in \R^d, $ для которых $
x_j >0 $ при $ j=1,\ldots,d,$ и для $ a \in \R_+^d, x \in \R^d $
положим $ a^x = a_1^{x_1} \ldots a_d^{x_d}.$

При $ d \in \N $ определим множества
\begin{eqnarray*}
I^d &=& \{x \in \R^d: 0 < x_j < 1,j=1,\ldots,d\},\\
\overline I^d &=& \{x \in \R^d: 0 \le x_j \le 1,j=1,\ldots,d\},\\
B^d &=& \{x \in \R^d: -1 \le x_j \le 1,j=1,\ldots,d\}.
\end{eqnarray*}

Через $ \e $ будем обозначать вектор в $ \R^d, $ задаваемый
равенством $ \e =(1,\ldots,1). $

Далее, напомним, что для области $ D \subset \R^d $ и вектора $ h
\in \R^d $ через $ D_h $ обозначается множество
$$
D_h = \{x \in D: x +th \in D \forall t \in \overline I\}.
$$

Для функции $ f, $ заданной в области $ D \subset \R^d, $ и
вектора $ h \in \R^d $ определим в $ D_h $ ее разность $ \Delta_h
f $ с шагом $ h, $ полагая
$$
(\Delta_h f)(x) = f(x+h) -f(x), x \in D_h,
$$
а для $ l \in \N: l \ge 2, $ в $ D_{lh} $ определим $l$-ую
разность $ \Delta_h^l f $ функции $ f $ с шагом $ h $ равенством
$$
(\Delta_h^l f)(x) = (\Delta_h (\Delta_h^{l-1} f))(x), x \in
D_{lh},
$$
положим также $ \Delta_h^0 f = f. $

При $ d \in \N $ для $ j=1,\ldots,d$ через $ e_j $ будем
обозначать вектор $ e_j = (0,\ldots,0,1_j,0,\ldots,0).$

При $ d \in \N $ для $ x \in \R^d $ через $\s(x) $ обозначим
множество $ \s(x) = \{j =1,\ldots,d: x_j \ne 0\}, $ а для
множества $ J \subset \{1,\ldots,d\} $ через $ \chi_J $ обозначим
вектор из $ \R^d $ с компонентами
$$
(\chi_J)_j = \begin{cases} 1, &\text{для} j \in J; \\ 0,
&\text{для} j \in (\{1,\ldots,d\} \setminus J). \end{cases}
$$

При $ d \in \N $ для $ x \in \R^d $ и $ J = \{j_1,\ldots,j_k\}
\subset \N: 1 \le j_1 < j_2 < \ldots < j_k \le d, $ через $ x^J $
обозначим вектор $ x^J = (x_{j_1},\ldots,x_{j_k}) \in \R^k, $ а
для множества $ A \subset \R^d $ положим $ A^J = \{x^J: x \in A\}. $

В заключение этого пункта введем еще несколько обозначений.

Для банахова пространства $ X $ (над $ \R$) обозначим $ B(X) = \{x
\in X: \|x\|_X \le 1\}. $

Для банаховых пространств $ X,Y $ через $ \mathcal B(X,Y) $
обозначим банахово пространство, состоящее из непрерывных линейных
операторов $ T: X \mapsto Y, $ с нормой
$$
\|T\|_{\mathcal B(X,Y)} = \sup_{x \in B(X)} \|Tx\|_Y.
$$
Отметим, что если $ X=Y,$ то $ \mathcal B(X,Y) $ является
банаховой алгеброй. Отметим еще, что при $ Y =  \R $ пространство
$ \mathcal B(X, \R) $ обозначается также $ X^*. $
\bigskip

1.2. В этом пункте содержатся сведения о кратных
рядах, которыми будем пользоваться в дальнейшем.

Для $ d \in \N, y \in \R^d $ положим
$$
\mn(y) = \min_{j=1,\ldots,d} y_j
$$
и для банахова пространства $ X, $ вектора $ x \in X $ и семейства
$ \{x_\kappa \in X, \kappa \in \Z_+^d\} $ будем писать $ x =
\lim_{ \mn(\kappa) \to \infty} x_\kappa, $ если для любого $
\epsilon >0 $ существует $ n_0 \in \N $ такое, что для любого $
\kappa \in \Z_+^d, $ для которого $ \mn(\kappa) > n_0, $
справедливо неравенство $ \|x -x_\kappa\|_X  < \epsilon. $

Пусть $ X $ -- банахово пространство (над $ \R $), $ d \in \N $ и
$ \{ x_\kappa \in X: \kappa \in \Z_+^d\} $ -- семейство векторов.
Тогда под суммой ряда $ \sum_{\kappa \in \Z_+^d} x_\kappa $ будем
понимать вектор $ x \in X, $ для которого выполняется равенство $
x = \lim_{\mn(k) \to \infty} \sum_{\kappa \in \Z_+^d(k)} x_\kappa. $

При $ d \in \N $ через $ \Upsilon^d $ обозначим множество
$$
\Upsilon^d = \{ \epsilon \in \Z^d: \epsilon_j \in \{0,1\},
j=1,\ldots,d\}.
$$

Имеет место

   Лемма 1.2.1

Пусть $ X $ -- банахово пространство, а вектор $ x \in X $ и
семейство $ \{x_\kappa \in X: \kappa \in \Z_+^d\} $ таковы, что $
x = \lim_{ \mn(\kappa) \to \infty} x_\kappa, $ Тогда для семейства $
\{ \mathcal X_\kappa \in X, \kappa \in \Z_+^d \}, $ определяемого
равенством
$$
\mathcal X_\kappa = \sum_{\epsilon \in \Upsilon^d: \s(\epsilon)
\subset \s(\kappa)} (-\e)^\epsilon x_{\kappa -\epsilon}, \kappa
\in \Z_+^d,
$$
справедливо равенство
$$
x = \sum_{\kappa \in \Z_+^d} \mathcal X_\kappa.
$$

Лемма является следствием того, что при $ k \in \Z_+^d $
выполняется равенство
\begin{equation*} \tag{1.2.1}
\sum_{\kappa \in \Z_+^d(k)}
\mathcal X_\kappa = x_k \text{(см. [5])}.
\end{equation*}

Замечание.

Легко заметить, что для любого семейства чисел $ \{x_\kappa \in
\R: x_\kappa \ge 0, \kappa \in \Z_+^d\} , $ если ряд $
\sum_{\kappa \in \Z_+^d} x_\kappa $ сходится, т.е. существует
предел $ \lim\limits_{\mn(k) \to \infty} \sum_{\kappa \in
\Z_+^d(k)} x_\kappa, $ то для любой последовательности подмножеств
$ \{Z_n \subset \Z_+^d, n \in \Z_+\}, $
 таких, что $ \card Z_n < \infty,
Z_n \subset Z_{n+1},  n \in \Z_+, $
и $ \cup_{ n \in \Z_+} Z_n = \Z_+^d, $
справедливо равенство
$ \sum_{\kappa \in \Z_+^d} x_\kappa =
\lim_{ n \to \infty} \sum_{\kappa \in Z_n} x_\kappa. $
Отсюда несложно понять, что если для семейства векторов
$ \{x_\kappa \in X, \kappa \in \Z_+^d\}  $ банахова пространства $ X $
ряд $ \sum_{\kappa \in \Z_+^d} \| x_\kappa \|_X $ сходится,
то для любой последовательности подмножеств
$ \{Z_n \subset \Z_+^d, n \in \Z_+\}, $
 таких, что $ \card Z_n < \infty,
Z_n \subset Z_{n+1},  n \in \Z_+, $
и $ \cup_{ n \in \Z_+} Z_n = \Z_+^d, $
в $ X $ соблюдается равенство
$ \sum_{\kappa \in \Z_+^d} x_\kappa =
\lim_{ n \to \infty} \sum_{\kappa \in Z_n} x_\kappa. $

При $ d \in \N $ для $ x \in \R^d $ обозначим
\begin{eqnarray*}
\mx(x)  &=& \max_{j=1,\ldots,d} x_j, \\
\cmx(x) &=& \card \{j \in \{1,\ldots,d\}: x_j = \mx(x)\}, \\
\cmn(x) &=& \card \{j \in \{1,\ldots,d\}: x_j = \mn(x)\}.
\end{eqnarray*}

Следующие три леммы используются в \S 3 при оценке точности приближений
на основе рассмотренных там кратных рядов.

Лемма 1.2.2

Пусть $ d \in \N, \beta \in \R_+^d, \alpha \in \R^d $ и $
\mx(\beta^{-1} \alpha) >0. $ Тогда существуют константы $
c_1(d,\alpha,\beta)  >0 $ и $ c_2(d,\alpha,\beta) >0 $ такие, что
для $ r \in \N $ соблюдается неравенство
\begin{equation*} \tag{1.2.2}
c_1 2^{ \mx(\beta^{-1} \alpha) r} r^{ \cmx(\beta^{-1} \alpha) -1}
\le \sum_{ \kappa \in \Z_+^d: (\kappa, \beta) \le r} 2^{(\kappa,
\alpha)} \le c_2 2^{ \mx(\beta^{-1} \alpha) r} r^{ \cmx(\beta^{-1}
\alpha) -1}.
\end{equation*}

Лемма 1.2.3

Пусть $ d \in \N, \alpha, \beta \in \R_+^d. $ Тогда существуют
константы $ c_3(d,\alpha,\beta) >0 $ и $ c_4(d,\alpha,\beta) >0 $
такие, что при $ r \in \N $ справедливо неравенство
\begin{equation*} \tag{1.2.3}
c_3 2^{-\mn(\beta^{-1} \alpha) r} r^{\cmn(\beta^{-1} \alpha) -1}
\le \sum_{\kappa \in \Z_+^d: (\kappa, \beta) > r} 2^{-(\kappa,
\alpha)} \le c_4 2^{-\mn(\beta^{-1} \alpha) r} r^{\cmn(\beta^{-1}
\alpha) -1}.
\end{equation*}

Доказательство лемм 1.2.2 и 1.2.3 приведено в [6].

Лемма 1.2.4

Пусть $ d \in \N, \beta \in \R_+^d, \epsilon >0, \tau \ge 0. $
Пусть еще $ J \subset \{1,\ldots,d\} $ -- непустое подмножество, $
J^\prime = \{1,\ldots,d\} \setminus J $ и $ \cmn = \card J. $
Тогда существует константа $ c_5(d,\beta,\epsilon,\tau,J) >0 $
такая, что для любого $ s \in \N $ выполняется неравенство
\begin{equation*} \tag{1.2.4}
\sum_{ \kappa \in \Z_+^d: s-1 < (\kappa, \beta) \le s}
2^{-\epsilon (\kappa^{J^\prime}, \beta^{J^\prime})} (1
+(\kappa^{J^\prime}, \beta^{J^\prime}))^\tau \le c_5 s^{\cmn -1}.
\end{equation*}

Доказательство.

В самом деле, используя неравенство
\begin{multline*}
\card \{ \kappa \in \Z_+^d: r -1 \le (\kappa, \alpha) \le r +1 \} \le
c_0(d, \alpha) (r +1)^{d -1}, \alpha \in \R_+^d, r \in \Z_+, d \in \N,
\end{multline*}
имеем
\begin{multline*}
\sum_{ \kappa \in \Z_+^d: s-1 < (\kappa, \beta) \le s} 2^{-\epsilon
(\kappa^{J^\prime}, \beta^{J^\prime})} (1 +(\kappa^{J^\prime},
\beta^{J^\prime}))^\tau \\
\le \sum_{i=1}^s \sum_{\kappa^J \in (\Z_+^d)^J, \kappa^{J^\prime}
\in (\Z_+^d)^{J^\prime}: i-1 \le (\kappa^J, \beta^J) \le i, s-i-1
\le (\kappa^{J^\prime}, \beta^{J^\prime}) \le s-i+1} 2^{-\epsilon
(\kappa^{J^\prime}, \beta^{J^\prime})} (1 +(\kappa^{J^\prime},
\beta^{J^\prime}))^\tau \\
= \sum_{i=1}^s \sum_{\kappa^J \in (\Z_+^d)^J: i-1 \le (\kappa^J,
\beta^J) \le i} \sum_{\kappa^{J^\prime} \in (\Z_+^d)^{J^\prime}:
s-i-1 \le (\kappa^{J^\prime}, \beta^{J^\prime}) \le s-i+1}
2^{-\epsilon (\kappa^{J^\prime}, \beta^{J^\prime})} (1
+(\kappa^{J^\prime}, \beta^{J^\prime}))^\tau \\
\le \sum_{i=1}^s \sum_{\kappa^J \in (\Z_+^d)^J: i-1 \le (\kappa^J,
\beta^J) \le i} \sum_{\kappa^{J^\prime} \in (\Z_+^d)^{J^\prime}:
s-i-1 \le (\kappa^{J^\prime}, \beta^{J^\prime}) \le s-i+1}
2^{-\epsilon(s-i-1)} (s-i+2)^\tau \\
\le \sum_{i=1}^s 2^{-\epsilon(s-i-1)} (s-i+2)^\tau \card
\{\kappa^J \in (\Z_+^d)^J: i-1 \le (\kappa^J, \beta^J) \le i+1
\}\times\\
\card \{\kappa^{J^\prime} \in (\Z_+^d)^{J^\prime}: s-i-1 \le
(\kappa^{J^\prime}, \beta^{J^\prime}) \le s-i+1\} \\
\le \sum_{i=1}^s 2^{-\epsilon(s-i-1)} (s-i+2)^\tau (i +1)^{\cmn
-1} (s -i +1)^{d -\cmn -1}\\
\le (s +1)^{\cmn -1} \sum_{k=0}^{s-1}
2^{-\epsilon(k -1)} (k +2)^\tau (k +1)^{d -\cmn -1} \\
\le 2^{\cmn
-1} s^{\cmn -1} \sum_{k=0}^\infty 2^{-\epsilon(k -1)} (k +2)^\tau
(k +1)^{d -\cmn -1} = c_6 s^{\cmn -1}. \square
\end{multline*}
\bigskip

1.3. В этом пункте приведем некоторые вспомогательные утверждения,
которые используются в следующем пункте и далее.

Как показано в  [7], [8], справедлива

   Лемма 1.3.1

    Пусть $ d \in \N, 1 \le p < \infty. $ Тогда

1) при  $ j=1,\ldots,d$ для любого непрерывного линейного
оператора $ T: L_p(I) \mapsto L_p(I) $ существует единственный
непрерывный линейный оператор $ \mathcal T^j: L_p(I^d) \mapsto
L_p(I^d), $ для которого для любой функции $ f \in L_p(I^d) $
почти для всех $ (x_1,\ldots,x_{j-1},x_{j+1}, \ldots,x_d) \in
I^{d-1} $ в $ L_p(I) $ выполняется равенство
\begin{equation*} \tag{1.3.1}
(\mathcal
T^j\!f)(x_1,\!\ldots,\!x_{j-1},\cdot,x_{j+1},\!\ldots,\!x_d)\!=\!
(T(f(x_1,\!\ldots,\!x_{j-1},\cdot,x_{j+1},\!\ldots,\!x_d)))(\cdot),
\end{equation*}

2) при этом, для каждого $ j=1,\ldots,d $ отображение $ V_j^{L_p},
$ которое каждому оператору $ T \in \mathcal B(L_p(I), L_p(I)) $
ставит в соответствие оператор $ V_j^{L_p}(T) = \mathcal T^j \in
\mathcal B(L_p(I^d), L_p(I^d)), $ удовлетворяющий (1.3.1),
является непрерывным гомоморфизмом банаховой алгебры $ \mathcal
B(L_p(I), L_p(I)) $ в банахову алгебру $ \mathcal B(L_p(I^d),
L_p(I^d)), $

3) причем, для любых операторов $ S,T \in \mathcal B(L_p(I),
L_p(I)) $ при любых $ i,j =1,\ldots,d: i \ne j, $ соблюдается
равенство
\begin{equation*} \tag{1.3.2}
(V_i^{L_p}(S) V_j^{L_p}(T))f = (V_j^{L_p}(T) V_i^{L_p}(S))f, f \in
L_p(I^d).
\end{equation*}

Замечание.

Если при $ d \in \N, 1 \le p \le q < \infty, $ оператор $ T \in
\mathcal B(L_p(I), L_p(I)) \cap \mathcal B(L_q(I), L_q(I)), $ то
при $ j=1,\ldots,d $ для $ f \in L_q(I^d)$ справедливо равенство
\newline
$ (V_j^{L_p} T)f = (V_j^{L_q} T)f. $ Поэтому символы $ L_p, L_q $
в качестве индексов у $ V_j $ можно опускать.
\bigskip

1.4. В этом пункте будут построены семейства операторов проектирования
на подпространства кусочно-полиномиальных функций, и описаны их
свойства, которые используются в п. 2.2. при доказательстве основных результатов
работы и в п. 3.3. при применении этих результатов.

Отметим сразу, что справедливость всех утверждений, приведенных в этом
пункте без доказательства, установлена в [5].

Для $ d \in \N, x,y \in \R^d $ будем писать $ x \le y (x < y), $
если для каждого $ j=1,\ldots,d $ выполняется неравенство $ x_j
\le y_j (x_j < y_j). $

Для $ d \in \N, m,n \in \Z^d: m \le n, $ обозначим
$$
\Nu_{m,n}^d = \{ \nu \in \Z^d: m \le \nu \le n \} = \prod_{j=1}^d
\Nu_{m_j,n_j}^1.
$$

При $ d \in \N $ для $ t \in \R^d $ через $ 2^t $ будем обозначать
вектор $ 2^t = (2^{t_1}, \ldots, 2^{t_d}). $

Для $ d \in \N, \kappa \in \Z^d, \nu \in \Z^d $ обозначим через
$ \chi_{\kappa, \nu}^d $ характеристическую функцию множества $
Q_{\kappa, \nu}^d = 2^{-\kappa} \nu +2^{-\kappa} I^d. $ Понятно,
что при $ d \in \N, \kappa \in \Z^d, \nu \in \Z^d $ имеют место
равенства
$$
Q_{\kappa,\nu}^d = \prod_{j=1}^d
Q_{\kappa_j,\nu_j}^1, \\
\chi_{\kappa,\nu}^d(x) = \prod_{j=1}^d
\chi_{\kappa_j,\nu_j}^1(x_j), x \in \R^d.
$$

Введем в рассмотрение пространства кусочно-полиномиальных функций.

Для $ d \in \N, l \in \Z_+^d $ и $ \kappa \in \Z_+^d$ через $
\mathcal P_\kappa^{d,l}$ обозначим линейное подпространство в $
L_\infty(I^d),$ состоящее из функций $ f \in L_\infty(I^d), $ для
каждой из которых существует набор  полиномов $ \{f_\nu \in
\mathcal P^{d,l}, \nu \in \Nu_{0, 2^\kappa -\e}^d\} $ такой, что
\begin{equation*} \tag{1.4.1}
f = \sum_{\nu \in \Nu_{0, 2^\kappa -\e}^d} f_\nu
\chi_{\kappa,\nu}^d.
\end{equation*}

При определении операторов проектирования на подпространства
$ \mathcal P_\kappa^{d,l} $ используются
операторы из следующего предложения.

Предложение 1.4.1

   Пусть $ d \in \N, l \in \Z_+^d. $ Тогда
для любого $ \delta \in \R_+^d $ и $ x^0 \in \R^d $ для $ Q = x^0
+\delta I^d $ существует единственный линейный оператор $
P_{\delta, x^0}^{d,l}: L_1(Q) \mapsto \mathcal P^{d,l}, $
обладающий следующими свойствами:

1) для $ f \in \mathcal P^{d,l} $ имеет место равенство
\begin{equation*}
P_{\delta, x^0}^{d,l}(f \mid_Q) = f,
  \end{equation*}

2)
\begin{equation*}
\ker P_{\delta,x^0}^{d,l} = \biggl\{\,f \in L_1(Q) : \int
\limits_{Q} f(x) g(x) \,dx =0\ \forall g \in \mathcal
P^{d,l}\,\biggr\},
\end{equation*}

причем, существуют константы $ c_1(d,l) >0 $ и $ c_2(d,l) >0 $
такие, что

   3) при $ 1 \le p \le \infty $ для $ f \in
L_p(Q) $ справедливо неравенство
  \begin{equation*}
\|P_{\delta, x^0}^{d,l}f \|_{L_p(Q)} \le c_1 \|f\|_{L_p(Q)},
  \end{equation*}

4) при $ 1 \le p < \infty $ для $ f \in L_p(Q) $ выполняется
неравенство
   \begin{equation*}   \|f -P_{\delta, x^0}^{d,l}f \|_{L_p(Q)} \le c_2 \sum_{j=1}^d
\delta_j^{-1/p} (\int_{\delta_j B^1} \int_{Q_{(l_j +1) \xi e_j}}
|\Delta_{\xi e_j}^{l_j +1} f(x)|^p dx d\xi)^{1/p}.
\end{equation*}

Для $ d \in \N, l, \kappa \in \Z_+^d, \nu \in \Nu_{0, 2^\kappa
-\e}^d $ определим непрерывный линейный оператор $ S_{\kappa,
\nu}^{d,l}: L_1(I^d) \mapsto \mathcal P^{d,l}(I^d) \cap L_\infty
(I^d), $ полагая для $ f \in L_1(I^d) $ значение
$$
S_{\kappa, \nu}^{d,l} f = P_{\delta, x^0}^{d,l} (f \mid_{(x^0
+\delta I^d)})
$$
при $ \delta = 2^{-\kappa}, x^0 = 2^{-\kappa} \nu $ (см. предложение
1.4.1).

Определим при $ d \in \N, l, \kappa \in \Z_+^d $ линейный
непрерывный оператор $ E_\kappa^{d,l}: L_1(I^d) \mapsto \mathcal
P_\kappa^{d,l} \cap L_\infty (I^d) $ равенством
$$
E_\kappa^{d,l} f = \sum_{\nu \in \Nu_{0, 2^\kappa -\e}^d}
(S_{\kappa, \nu}^{d,l} f) \chi_{\kappa, \nu}^d, f \in L_1(I^d).
$$

Следующая лемма будет полезна как в этом пункте, так и в п. 2.2.

Леммаа 1.4.2

Пусть $ d \in \N $ и $ \kappa, \kappa^\prime, \nu, \nu^\prime \in \Z^d $
таковы, что $ \kappa^\prime \le \kappa, $ а $ Q_{\kappa, \nu}^d \cap
Q_{\kappa^\prime, \nu^\prime}^d \ne \emptyset. $
Тогда имеет место включение
\begin{equation*} \tag{1.4.2}
Q_{\kappa, \nu}^d \subset
Q_{\kappa^\prime, \nu^\prime}^d.
\end{equation*}

Доказательство.
В условиях леммы, выбирая $ x \in Q_{\kappa, \nu}^d \cap
Q_{\kappa^\prime, \nu^\prime}^d, $ имеем
$$
2^{-\kappa} \nu < x < 2^{-\kappa^\prime} \nu^\prime +2^{-\kappa^\prime}
$$
и
$$
2^{-\kappa^\prime} \nu^\prime < x < 2^{-\kappa} \nu +2^{-\kappa},
$$
откуда
$$
 \nu < 2^{\kappa -\kappa^\prime} \nu^\prime +2^{\kappa -\kappa^\prime}
$$
и
$$
2^{\kappa -\kappa^\prime} \nu^\prime <  \nu +\e,
$$
а, следовательно,
$$
 \nu \le 2^{\kappa -\kappa^\prime} \nu^\prime +2^{\kappa -\kappa^\prime} -\e
$$
и
$$
2^{\kappa -\kappa^\prime} \nu^\prime \le  \nu
$$
или
$$
2^{\kappa -\kappa^\prime} \nu^\prime \le  \nu \le
2^{\kappa -\kappa^\prime} \nu^\prime +2^{\kappa -\kappa^\prime} -\e.
$$
Поэтому для $ x \in Q_{\kappa,\nu}^d $ выполняются
соотношения
$$
2^{-\kappa} \nu < x < 2^{-\kappa} \nu +2^{-\kappa}
$$
а, значит,
$$
2^{-\kappa} 2^{\kappa -\kappa^\prime} \nu^\prime \le  2^{-\kappa} \nu
< x < 2^{-\kappa} \nu +2^{-\kappa} \le
2^{-\kappa} (2^{\kappa -\kappa^\prime} \nu^\prime +2^{\kappa -\kappa^\prime}
-\e) +2^{-\kappa}
$$
или
$$
2^{-\kappa^\prime} \nu^\prime < x <
2^{-\kappa^\prime} \nu^\prime
+2^{-\kappa^\prime},
$$
т.е. $ x \in Q_{\kappa^\prime, \nu^\prime}^d. \square $

Из (1.4.2) следует, что при $ d \in \N, l \in \Z_+^d $ для $  \kappa,
\kappa^\prime  \in \Z_+^d: \kappa^\prime \le \kappa, $ справедливо
включение
\begin{equation*} \tag{1.4.3}
\mathcal P_{\kappa^\prime}^{d,l} \subset \mathcal P_\kappa^{d,l}.
\end{equation*}

Заметим, что ввиду (1.4.3) при $ d \in \N, l, \kappa \in \Z_+^d $
для $ \epsilon \in \Upsilon^d: \s(\epsilon) \subset \s(\kappa), $
справедливо включение $ \mathcal P_{\kappa -\epsilon}^{d,l}
\subset  \mathcal P_\kappa^{d,l}. $

Принимая во внимание это обстоятельство, при $ d \in \N, l, \kappa
\in \Z_+^d $ определим линейный непрерывный оператор $ \mathcal
E_\kappa^{d,l}: L_1(I^d) \mapsto \mathcal P_\kappa^{d,l} \cap
L_\infty(I^d), $ полагая
$$
\mathcal E_\kappa^{d,l} = \sum_{\epsilon \in \Upsilon^d:
\s(\epsilon) \subset \s(\kappa)} (-\e)^\epsilon E_{\kappa
-\epsilon}^{d,l}.
$$

Лемма 1.4.3

Пусть $ d \in \N, l, \kappa \in \Z_+^d. $ Тогда имеет место
равенство
\begin{equation*} \tag{1.4.4}
\mathcal E_\kappa^{d,l} =
\prod_{j=1}^d V_j(\mathcal E_{\kappa_j}^{1,l_j}).
\end{equation*}

Лемма 1.4.4

Пусть  $ d \in \N, l \in \Z_+^d.$ Тогда
справедливы следующие утверждения:

1) при $ \kappa \in \Z_+^d $ для $ f \in \mathcal P_\kappa^{d,l} $
соблюдается равенство
\begin{equation*} \tag{1.4.5}
E_\kappa^{d,l} f = f;
\end{equation*}

2) при $ \kappa \in \Z_+^d $ ядро
\begin{equation*}
\ker E_\kappa^{d,l} = \biggl\{ \, f \in L_1(I^d):
\int\limits_{I^d} f(x) g(x) \,dx =0 \ \forall g \in \mathcal
P_\kappa^{d,l} \,\biggr\}.
\end{equation*}

Лемма  1.4.5

При $ d \in \N, l \in \Z_+^d $ для $ \kappa, \kappa^\prime \in
\Z_+^d $ соблюдаются равенства
\begin{equation*} \tag{1.4.6}
\mathcal E_\kappa^{d,l} \mathcal E_{\kappa^\prime}^{d,l} =
\begin{cases} \mathcal E_\kappa^{d,l},
 &\text{при $ \kappa = \kappa^\prime $}; \\
       0, &\text{при $ \kappa \ne \kappa^\prime $}.
\end{cases}
\end{equation*}

При $ d \in \N, l, \kappa \in \Z_+^d $ положим
$$
\mathfrak P_\kappa^{d,l} = \mathcal E_\kappa^{d,l} (\mathcal
P_\kappa^{d,l} ) \subset \mathcal P_\kappa^{d,l}.
$$

Лемма 1.4.6

Пусть $ d \in \N, l, \kappa \in \Z_+^d. $ Тогда справедливы
следующие соотношения:

1)
\begin{equation*} \tag{1.4.7}
\Im \mathcal E_\kappa^{d,l} = \mathfrak P_\kappa^{d,l};
\end{equation*}

2) для любых $ f, g \in L_1(I^d) $ выполняется равенство
\begin{equation*} \tag{1.4.8}
\int_{I^d} (\mathcal E_\kappa^{d,l} f) \cdot g dx = \int_{I^d} f
\cdot (\mathcal E_\kappa^{d,l} g) dx;
\end{equation*}

3)
\begin{equation*}
\ker \mathcal E_\kappa^{d,l} = \{ f \in L_1(I^d): \int_{I^d} f g
dx =0 \forall g \in \mathfrak P_\kappa^{d,l};
\end{equation*}

4) для $ \kappa^\prime \in \Z_+^d: \kappa^\prime \ne \kappa, $ и $
f, g \in L_1(I^d) $ верно равенство
\begin{equation*} \tag{1.4.9}
\int_{I^d} (\mathcal E_\kappa^{d,l} f) (\mathcal
E_{\kappa^\prime}^{d,l} g) dx =0.
\end{equation*}

Отметим некоторые особенности подпространств $ \mathfrak
P_\kappa^{d,l}, d \in \N, l,\kappa \in \Z_+^d. $

При $ d \in \N, l \in \Z_+^d $ для $ \kappa \in \Z_+^d, J =
\s(\kappa), \rho \in \Nu_{0, 2^{\kappa -\chi_J} -\e}^d $ обозначим
через $ H_{\kappa, \rho}^{d,l} $ подпространство в $
L_\infty(Q_{\kappa -\chi_J, \rho}^d), $ состоящее из функций $ h,
$ для каждой из которых существует функция $ g \in \mathfrak
P_\kappa^{d,l} $ такая, что $ h = g \mid_{Q_{\kappa -\chi_J,
\rho}^d}, $ а через $ U_{\kappa, \rho}^{d,l}: \mathfrak
P_\kappa^{d,l} \mapsto H_{\kappa, \rho}^{d,l} $ обозначим линейное
отображение, которое каждому $ g \in \mathfrak P_\kappa^{d,l} $
ставит в соответствие $ h = U_{\kappa, \rho}^{d,l} g = g
\mid_{Q_{\kappa -\chi_J, \rho}^d}, $ обозначим еще через $
H_\kappa^{d,l} $ прямое произведение $ H_\kappa^{d,l} =
\prod_{\rho \in \Nu_{0, 2^{\kappa -\chi_J} -\e}^d} H_{\kappa,
\rho}^{d,l} $ пространств $  H_{\kappa, \rho}^{d,l}, $ а также
через $ U_\kappa^{d,l}: \mathfrak P_\kappa^{d,l} \mapsto
H_\kappa^{d,l} $ -- отображение, которое каждому $ g \in \mathfrak
P_\kappa^{d,l} $ сопоставляет $ h = \{ h_\rho, \rho \in \Nu_{0,
2^{\kappa -\chi_J} -\e}^d\} $ с компонентами $ h_\rho = U_{\kappa,
\rho}^{d,l} g, \rho \in \Nu_{0, 2^{\kappa -\chi_J} -\e}^d. $

Справедлива

Лемма 1.4.7

При $ d \in \N, l, \kappa \in \Z_+^d $ отображение $
U_\kappa^{d,l} $ является изоморфизмом $ \mathfrak P_\kappa^{d,l}
$ на $ H_\kappa^{d,l}. $

Далее, при $ d \in \N $ для $ \kappa \in \Z_+^d $  и $ \nu \in
\Z^d $ обозначим через $ A_{\kappa, \nu}^d  $ отображение, которое
каждой функции $ f  $, заданной на некотором множестве  $ S
\subset \R^d, $ ставит в соответствие функцию $ A_{\kappa,  \nu}^d
f, $ определяемую на множестве $ \{ x \in \R^d: 2^{-\kappa} \nu
+2^{-\kappa} x \in S\} $ равенством $ (A_{\kappa, \nu}^d f)(x) =
f(2^{-\kappa} \nu +2^{-\kappa}  x). $ Так как для $ \kappa \in
\Z_+^d, \nu \in \Z^d $  отображение $  \R^d \ni x \mapsto
2^{-\kappa} \nu +2^{-\kappa} x \in \R^d $ --- взаимно однозначно,
то отображение $ A_{\kappa, \nu}^d $ является биекцией на  себя
множества всех функций с  областью  определения в $ \R^d. $

Лемма 1.4.8

При $ d \in \N, l \in \Z_+^d, \kappa \in \Z_+^d, J = \s(\kappa),
\rho \in \Nu_{0, 2^{\kappa -\chi_J} -\e}^d$ отображение $
A_{\kappa -\chi_J, \rho}^d \mid_{H_{\kappa, \rho}^{d,l}} $
является изоморфизмом $ H_{\kappa, \rho}{d, l} $ на $
H_{\chi_J}^{d, l}. $

При $ d \in \N, l, \kappa \in \Z_+^d $ положим
$$
\mathfrak R_\kappa^{d,l} = \dim \mathfrak P_\kappa^{d,l}.
$$

Лемма 1.4.9

Пусть $ d \in \N, l, \kappa \in \Z_+^d, J = \s(\kappa), 1 \le p
\le \infty. $ Тогда имеет место равенство
\begin{equation*} \tag{1.4.10}
\mathfrak R_\kappa^{d,l} = \mathfrak R_{\chi_J}^{d,l} 2^{(\kappa
-\chi_J, \e)},
\end{equation*}
и можно построить линейный изоморфизм $ \mathfrak I_\kappa^{d,l} $
подпространства $ \mathfrak P_\kappa^{d,l} $ на
пространство $ \R^{\mathfrak R_\kappa^{d,l}}, $ обладающий тем
свойством, что для $ g \in \mathfrak P_\kappa^{d,l} $ выполняются
неравенства
\begin{equation*} \tag{1.4.11}
c_3 \| g \|_{L_p(I^d)} \le 2^{-(\kappa -\chi_J, \e)/p} \|
\mathfrak I_\kappa^{d,l} g \|_{l_p^{\mathfrak R_\kappa^{d,l}}} \le
c_4 \| g \|_{L_p(I^d)}
\end{equation*}
с некоторыми константами $ c_3 >0, c_4 >0, $ зависящими только
от $ d, l, p. $

Лемма 1.4.10

Пусть $ d \in \N, l \in \Z_+^d, J \subset \{1, \ldots,d \}: J \ne
\emptyset, $ и система функций $ \phi_i^{d,l,J} \in L_\infty(\R^d),
i =1, \ldots, \mathfrak R_{\chi_J}^{d, l} $ такова, что
$ \phi_i^{d,l,J}(x) =0 $ при $ x \in \R^d \setminus I^d,
i =1, \ldots, \mathfrak R_{\chi_J}^{d, l}, $
а
$$
\{ \phi_i^{d,l,J} \mid_{I^d} \in \mathfrak P_{\chi_J}^{d, l},
i =1, \ldots, \mathfrak R_{\chi_J}^{d, l} \}
$$
является ортонормированным базисом в $ \mathfrak P_{\chi_J}^{d, l} \cap
L_2(I^d).$
Тогда для любого $ \kappa \in \Z_+^d: \s(\kappa) = J, $ система функций
$$
\{ 2^{( \kappa -\chi_J, \e) /2} \phi_i^{d,l,J}( 2^{ \kappa -\chi_J} x
-\rho), i =1, \ldots, \mathfrak R_{\chi_J}^{d, l}, \rho \in
\Nu_{0, 2^{\kappa -\chi_J} -\e}^d \}
$$
образует ортонормированный базис в
$ \mathfrak P_\kappa^{d, l} \cap L_2(I^d), $
и для $ f \in L_1(I^d) $ почти для всех $ x \in I^d $
выполняется равенство
\begin{multline*} \tag{1.4.12}
( \mathcal E_\kappa^{d,l} f)(x) = \sum_{ i =1, \ldots, \mathfrak
R_{\chi_J}^{d, l}, \rho \in \Nu_{0, 2^{\kappa -\chi_J} -\e}^d }
2^{( \kappa -\chi_J, \e)}\\
\times \biggl(\int_{ I^d} \phi_i^{d,l,J}( 2^{ \kappa -\chi_J} y
-\rho) f(y) dy\biggr) \phi_i^{d,l,J}( 2^{ \kappa -\chi_J} x
-\rho) \\
= \sum_{ i =1, \ldots, \mathfrak R_{\chi_J}^{d, l}, \rho \in
\Nu_{0, 2^{\kappa -\chi_J} -\e}^d } 2^{( \kappa -\chi_J, \e)}\\
\times \biggl(\int_{ Q_{\kappa -\chi_J, \rho}^d } \phi_i^{d,l,J}(
2^{ \kappa -\chi_J} y -\rho) f(y) dy\biggr) \phi_i^{d,l,J}( 2^{
\kappa
-\chi_J} x -\rho).\\
\end{multline*}

Доказательство.

Тот факт, что система функций
$$
\{ 2^{( \kappa -\chi_J, \e) /2} \phi_i^{d,l,J}( 2^{ \kappa -\chi_J} x
-\rho), i =1, \ldots, \mathfrak R_{\chi_J}^{d, l}, \rho \in
\Nu_{0, 2^{\kappa -\chi_J} -\e}^d \}
$$
образует ортонормированный базис в
$ \mathfrak P_\kappa^{d, l} \cap L_2(I^d), $
легко установить, используя леммы 1.4.7 и 1.4.8, а также соотношения
$$
\phi_i^{d,l,J}( 2^{ \kappa -\chi_J} x -\rho) =
((A_{ \kappa -\chi_J, \rho})^{-1} \phi_i^{d,l,J})(x), x \in \R^d;
$$
$$
\phi_i^{d,l,J}( 2^{ \kappa -\chi_J} x -\rho) =0, \text{при}
x \in \R^d \setminus Q_{ \kappa -\chi_J, \rho}^d,
i =1, \ldots, \mathfrak R_{\chi_J}^{d, l}, \rho \in
\Nu_{0, 2^{\kappa -\chi_J} -\e}^d ;
$$
и равенство (1.4.10).

Равенство (1.4.12) вытекает из (1.4.6), (1.4.7), (1.4.8) и разложения
вектора в ортонормированном базисе.
$ \square $

Отметим еще, что в условиях леммы 1.4.10 в силу включения $
\mathfrak P_{\chi_J}^{d, l} \subset  \mathcal P_{\chi_J}^{d, l} $
и с учетом (1.4.1) имеют место соотношения
\begin{equation*} \tag{1.4.13}
(\phi_i^{d,l,J}) \mid_{ Q_{\chi_J, \nu}^d} \in
\mathcal P^{d,l}(Q_{\chi_J, \nu}^d), \nu \in \Nu_{0, 2^{\chi_J} -\e}^d,
i =1, \ldots, \mathfrak R_{\chi_J}^{d, l}.
\end{equation*}
\bigskip

1.5. В этом пункте содержатся утверждения, касающиеся
аппроксимационных свойств операторов проектирования, построенных в п. 1.4.,
которые понадобятся в следующих параграфах.

Из [7], [8] можно извлечь лемму 1.5.1.

Лемма 1.5.1

Пусть $ d \in \N, l \in \N^d, 1 \le p < \infty. $ Тогда существуют
константы $ c_1(d,l) >0, c_2(d) >0, c_3(d, l) >0 $ и $ c_4(d,l) >0 $
такие, что для любой функции $ f \in L_p(I^d) $ при $ \kappa \in \Z_+^d $
справедливы соотношения
\begin{equation*} \tag{1.5.1}
\begin{split}
&\|f -E_\kappa^{d,l -\e}f \|_{L_p(I^d)} \le\\
&\le c_1 \sum_{j=1}^d 2^{\kappa_j /p} \biggl(\int_{ c_2
2^{-\kappa_j} B^1} \int_{(I^d)_{l_j \xi e_j}} |\Delta_{\xi
e_j}^{l_j} f(x)|^p dx
d\xi\biggr)^{1/p} \to 0 \\
&\text{при}\, \mn(\kappa) \to \infty,
\end{split}
\end{equation*}
\begin{equation*}
\| E_\kappa^{d,l -\e}f \|_{L_p(I^d)} \le
 c_3 \| f \|_{L_p(I^d)},
\end{equation*}
\begin{equation*} \tag{1.5.2}
\| \mathcal E_\kappa^{d,l -\e}f \|_{L_p(I^d)} \le
 c_4 \| f \|_{L_p(I^d)}.
\end{equation*}

Из (1.5.1) и леммы 1.2.1 вытекает

Следствие

В условиях леммы 1.5.1 для $ f \in L_p(I^d) $ в пространстве
$ L_p(I^d) $ имеет место равенство
\begin{equation*} \tag{1.5.3}
f = \sum_{ \kappa \in \Z_+^d} \mathcal E_\kappa^{d,l -\e} f.
\end{equation*}

В [5] установлена

Теорема 1.5.2

Пусть $ d \in \N, l \in \N^d. $ Тогда для любой функции $ f \in
L_2(I^d) $ справедливо равенство
\begin{equation*} \tag{1.5.4}
\| f \|_{L_2(I^d)} = (\sum_{ \kappa \in \Z_+^d} \| \mathcal
E_\kappa^{d,l -\e} f \|_{L_2(I^d)}^2)^{1/2}.
\end{equation*}
\bigskip

\centerline{ \S 2. Теорема типа Литтлвуда-Пэли}
\centerline{и оценка сверху нормы функции $ \| f \|_{L_p(I^d)} $ }
\centerline{ через нормы проекций $ \| \mathcal E_\kappa^{d,l} f \|_{L_p(I^d)} $ }
\bigskip

2.1. В этом пункте приведены вспомогательные утверждения,
на которые опирается доказательство занимающей центральное место в работе
леммы 2.2.1. Все утверждения, доказательство которых опущено,
взяты из [9, гл. I].

Для локально суммируемой в $ \R^d$ функции  $ f $ определим её
максимальную функцию $ M_f, $ полагая
$$
M_f(x) = \sup_{r \in \R_+} (1/ \mes B(x,r)) \int_{B(x,r)} |f(y)|
dy, x \in \R^d,
$$
где $ B(x,r) = \{y \in \R^d: |y -x| < r\}, $ а $ |z| = \|z\|_{l_2^d} $
для $ z \in \R^d. $

Лемма 2.1.1

Пусть $ d \in \N $. Тогда существует константа $ c_1(d) >0 $
такая, что для любой функции $ f \in L_1(\R^d) $ при $ \alpha \in
\R_+ $ множество $ \{x \in \R^d: M_f(x) > \alpha\} $ открыто и
выполняется неравенство
\begin{equation*} \tag{2.1.1}
\mes \{x \in \R^d: M_f(x) > \alpha\} \le (c_1 / \alpha)
\int_{\R^d} |f| dx.
\end{equation*}

Лемма 2.1.2
Для $ f \in L_1(\R^d) $ почти для всех $ x \in \R^d $
выполняется равенство
\begin{equation*} \tag{2.1.2}
f(x) = \lim_{ r \to 0} (1/ \mes B(x,r)) \int_{B(x,r)} f(y) dy.
\end{equation*}

Для замкнутого множества $ F \subset \R^d $ и $ x \in \R^d $
положим
$$
\rho(x) = \rho(x,F) = \min_{y \in F} | x -y|.
$$

Лемма 2.1.3

Пусть $ d \in \N. $ Тогда существует константа $ c_2(d) >0 $
такая, что для любого замкнутого множества $ F \subset \R^d, $
для которого $ \mes (\R^d \setminus F) < \infty, $ функция $ M(x), $
определяемая для $ x \in F $ равенством
$$
M(x) = \int_{\R^d} \rho(u) |x -u|^{-(d+1)} du,
$$
суммируема на $ F $ и соблюдается  неравенство
\begin{equation*} \tag{2.1.3}
\int_F M(x) dx \le c_2 \mes (\R^d \setminus F).
\end{equation*}

Лемму 2.1.4 можно рассматривать как некоторый уточняющий вариант
теоремы 3 из гл. I в [9].

Лемма 2.1.4

Пусть $ d \in \N. $ Тогда существуют константы $ c_3(d) >0,
c_4(d) >0 $ такие, что для любого замкнутого множества $ F \subset \R^d, $
у которого для $ W =
\R^d \setminus F $ мера $ \mes W < \infty, $
существует семейство кубов $ \{ Q_r, r \in \N \}, $
обладающих следующими свойствами:

1) для каждого $ r \in \N $ куб $ Q_r $ имеет вид
\begin{equation*} \tag{2.1.4}
Q_r = 2^{-k^r} \nu^r +2^{-k^r} I^d,
\end{equation*}
где $ k^r \in \Z, \nu^r \in \Z^d; $

2) \begin{equation*} \tag{2.1.5}
Q_r \cap Q_s = \emptyset, \text{ для}   r,s \in \N: r \ne s;
\end{equation*}

3) \begin{equation*} \tag{2.1.6}
W = \cup_{ r \in \N} \overline Q_r,
\end{equation*}
где  $ \overline Q_r =
 2^{-k^r} \nu^r +2^{-k^r} \overline I^d, r \in \N; $

4) при $ r \in \N $ справедливы неравенства
\begin{equation*} \tag{2.1.7}
c_3 \diam Q_r < \inf_{ x \in Q_r} \rho(x, F) \le
c_4 \diam Q_r,
\end{equation*}
где  $ \diam Q = \sup_{x,y \in Q} |x -y|, Q \subset \R^d. $

Доказательство.

Поскольку множество $ W $ открыто и мера $ \mes W < \infty, $
то множество
$ \{k \in \Z| \exists \nu \in \Z^d: \inf_{x \in Q_{k \e, \nu}^d}
\rho(x, F) > d^{1/2} 2^{-k} \} $ -- непусто и ограничено снизу.
Поэтому существует  $ k_0 \in \Z, $ для которого соблюдается равенство
$ k_0 = \min \{k \in \Z| \exists \nu \in \Z^d: \inf_{x \in Q_{k \e, \nu}^d}
\rho(x, F) > d^{1/2} 2^{-k} \}. $
Положим
$ \Nu_0 = \{ \nu \in \Z^d: \inf_{x \in Q_{k_0 \e, \nu}^d}
\rho(x, F) > d^{1/2} 2^{-k_0} \}, $ а
$ \mathfrak Q_0 = \{ Q_{k_0 \e, \nu}^d: \nu \in \Nu_0 \}. $
И определим по индукции для $ k \in \N $ множества
\begin{multline*}
\Nu_k = \{ \nu \in \Z^d: \inf_{x \in Q_{(k_0 +k) \e, \nu}^d}
\rho(x, F) > d^{1/2} 2^{-(k_0 +k)}, \\ Q_{(k_0 +k) \e, \nu}^d
\cap Q_{(k_0 +k^\prime) \e, \nu^\prime}^d = \emptyset \forall
\nu^\prime \in \Nu_{k^\prime}, k^\prime =0, \ldots, k -1 \}, \\
\mathfrak Q_k = \{ Q_{(k_0 +k) \e, \nu}^d: \nu \in \Nu_k \}.
\end{multline*}

Проверим, что для кубов из $ \cup_{k \in \Z_+} \mathfrak Q_k $ соблюдаются
условия (2.1.4) -- (2.1.7).

Соотношения (2.1.4),  (2.1.5) и первое неравенство в (2.1.7)
соблюдаются по построению. Покажем, что второе неравенство в (2.1.7) также выполнено.
Пусть $ k \in \Z_+, \nu \in \Nu_k. $ Тогда, учитывая, что
$ \R^d = \cup_{\nu^\prime \in \Z^d} \overline Q_{(k_0 +k -1) \e,
\nu^\prime}^d, $
выберем $ \nu^\prime \in \Z^d, $ для которого $ Q_{(k_0 +k) \e, \nu}^d \cap
\overline Q_{(k_0 +k -1) \e, \nu^\prime}^d \ne \emptyset, $
а, следовательно, и
$ Q_{(k_0 +k) \e, \nu}^d \cap Q_{(k_0 +k -1) \e, \nu^\prime}^d \ne \emptyset. $
При этом, ввиду (1.4.3) имеет место включение
\begin{equation*} \tag{2.1.8}
Q_{(k_0 +k) \e, \nu}^d \subset Q_{(k_0 +k -1) \e, \nu^\prime}^d.
\end{equation*}
Проверим, что
\begin{equation*} \tag{2.1.9}
\inf_{x \in Q_{(k_0 +k -1) \e, \nu^\prime}^d} \rho(x, F) \le
d^{1/2} 2^{-(k_0 +k -1)}.
\end{equation*}
Если $ k=0, $ то (2.1.9) соблюдаеется в силу выбора $ k_0. $
При $ k =1 $ выполнение (2.1.9)  с учётом (2.1.8) следует из
определения множеств $ \Nu_0, \Nu_1. $
А при $ k \in \N: k > 1, $ предположим, что
$$
\inf_{x \in Q_{(k_0 +k -1) \e, \nu^\prime}^d} \rho(x, F) >
d^{1/2} 2^{-(k_0 +k -1)}.
$$
Тогда вследствие (2.1.8),  определения множества $ \Nu_k $ и включения
$ \nu \in \Nu_k $ мультииндекс $\nu^\prime \notin \Nu_{k -1}, $
и, значит, существуют $ j \in \Z_+: j < k -1, $  и $ n \in \Nu_j $ такие, что
$$
Q_{(k_0 +k -1) \e, \nu^\prime}^d \cap
Q_{(k_0 +j) \e, n}^d \ne \emptyset,
$$
что в силу (1.4.2), (2.1.8) влечёт включение
$$
Q_{(k_0 +k) \e, \nu}^d \subset
Q_{(k_0 +j) \e, n}^d,
$$
которое противоречит принадлежности $ \nu \in \Nu_k. $ Полученное
противоречие подтверждает справедливость (2.1.9). С учётом (2.1.9)
выбирая для произвольного $ \epsilon >0 $ точки $ y \in F $ и $
\xi \in Q_{(k_0 +k -1) \e, \nu^\prime}^d $ так, чтобы было $ | y
-\xi | < d^{1/2} 2^{-(k_0 +k -1)} +\epsilon, $ для $ x \in Q_{(k_0
+k) \e, \nu}^d, $ учитывая (2.1.8), имеем

\begin{multline*}
| x -y | \le | x -\xi | +| y -\xi| < \diam Q_{(k_0 +k -1) \e,
\nu^\prime}^d +d^{1/2} 2^{-(k_0 +k -1)} +\epsilon\\ = 2d^{1/2}
2^{-(k_0 +k -1)} +\epsilon = 4d^{1/2} 2^{-(k_0 +k)} +\epsilon = 4
\diam Q_{(k_0 +k) \e, \nu}^d +\epsilon,
\end{multline*}
откуда
$$
\inf_{ x \in Q_{(k_0 +k) \e, \nu}^d} \rho(x, F) =
\inf_{ x \in Q_{(k_0 +k) \e, \nu}^d, y \in F} | x -y| <
4 \diam Q_{(k_0 +k) \e, \nu}^d +\epsilon,
$$
что в силу произвольности $ \epsilon >0 $ даёт второе неравенство в (2.1.7).

Теперь получим (2.1.6). Для $ x_0 \in W, $ благодаря открытости $ W, $
возьмём $ \epsilon >0, $ для которого $ B(x_0, 3 \epsilon) \subset W, $
и найдём $ k \in \Z_+ $ такое, что $ d^{1/2} 2^{-(k_0 +k)} < \epsilon. $
Выбирая $ \nu \in \Z^d, $ для которого
$ x_0 \in \overline Q_{(k_0 +k) \e, \nu}^d, $
видим, что для $ x \in Q_{(k_0 +k) \e, \nu}^d, y \in F $
выполняется неравенство
\begin{multline*}
| x -y| = | x -x_0 +x_0 -y | \ge | x_0 -y | -| x -x_0| \ge 3 \epsilon -
\diam \overline Q_{(k_0 +k) \e, \nu}^d \\ = 3 \epsilon -d^{1/2} 2^{-(k_0 +k)} >
3d^{1/2} 2^{-(k_0 +k)} -d^{1/2} 2^{-(k_0 +k)} = 2d^{1/2} 2^{-(k_0 +k)},
\end{multline*}
и, следовательно,
$$
\inf_{ x \in Q_{(k_0 +k) \e, \nu}^d} \rho(x, F) =
\inf_{ x \in Q_{(k_0 +k) \e, \nu}^d, y \in F} | x -y| \ge
2d^{1/2} 2^{-(k_0 +k)} > d^{1/2} 2^{-(k_0 +k)}.
$$
Принимая во внимание сказанное, получаем, что если $ \nu \in \Nu_k, $
то $ x_0 \in \cup_{ r \in \N} \overline Q_r. $
Если же $ \nu \notin \Nu_k, $ то согласно определению $ \Nu_k $
существуют $ j \in \Z_+: j < k, \nu^\prime \in \Nu_j $ такие, что
$$
Q_{(k_0 +k) \e, \nu}^d \cap Q_{(k_0 +j) \e, \nu^\prime}^d \ne \emptyset,
$$
и, значит, (см. (1.4.2))
$$
Q_{(k_0 +k) \e, \nu}^d \subset Q_{(k_0 +j) \e, \nu^\prime}^d,
$$
а
$$
\overline Q_{(k_0 +k) \e, \nu}^d \subset
\overline Q_{(k_0 +j) \e, \nu^\prime}^d,
$$
т.е. $ x_0 \in \overline Q_{(k_0 +j) \e, \nu^\prime}^d
\subset \cup_{ r \in \N} \overline Q_r. \square $

Предложение 2.1.5

Пусть $ d \in \N. $ Тогда существует константа $ C_5(d) >0 $ такя, что
для любой функции $ f \in L_1(\R^d) $  при любом $ \alpha \in \R_+ $
существуют замкнутое множесво $ F \subset \R^d $ и семейство кубов
$ \{ Q_r, r \in \N \}, $ со следующими свойствами:

1) почти для всех $ x \in F $  выполняется неравнство
\begin{equation*} \tag{2.1.10}
| f(x) | \le \alpha;
\end{equation*}

2) для $ W = \R^d \setminus F $ справедливо неравенство
\begin{equation*} \tag{2.1.11}
\mes W \le (c_1 / \alpha) \int_{\R^d} |f| dx;
\end{equation*}

3) для кубов семейства $ \{Q_r, r \in \N \} $ соблюдаются соотношения
(2.1.4) -- (2.1.7), а также

4) при  $ r \in \N $ имеет место оценка
\begin{equation*} \tag{2.1.12}
(1 / \mes Q_r) \int_{ Q_r} | f(x) | dx \le c_5 \alpha.
\end{equation*}

Схема доказательства предложения 2.1.5 взята из \S 3 гл. I в [9].

Доказательство.

Для $ f \in L_1(\R^d) $ и $ \alpha \in \R_+ $ определим множество
$ F $ равенством
$$
F = \{ x \in \R^d: M_f(x) \le \alpha \}.
$$
ввиду леммы 2.1.1 множество $ W = \R^d \setminus F = \{x \in
\R^d: M_f(x) > \alpha \} $ -- открыто и, значит, множество $ F $
-- замкнуто. А из (2.1.1) следует (2.1.11). Благодаря (2.1.2),
почти для всех $ x \in F $ имеем
$$
| f(x) | = \lim_{ r \to 0} (1/ \mes B(x,r)) \int_{B(x,r)} | f(y) | dy
\le M_f(x) \le \alpha,
$$
т.е. выполняется (2.1.10).

В соответствии с леммой 2.1.4 для множества $ F $ построим
семейство кубов $ \{ Q_r: r \in \N \}, $  удовлетворяющих условиям
(2.1.4) -- (2.1.7). Для проверки (2.1.12), принимая во внимание
(2.1.7), при $ r \in \N $ выберем $ \xi_r \in Q_r $ и $ x_r \in F $ так,
чтобы соблюдалось неравенство
$ | \xi_r -x_r | < 2 c_4 \diam Q_r. $ Тогда для $ x \in Q_r $ справедливо
неравенство
$$
| x -x_r | \le | x -\xi_r | +| \xi_r -x_r | \le \diam Q_r +
2 c_4 \diam Q_r = c_6 \diam Q_r = \delta_r
$$
или $ Q_r \subset B(x_r, \delta_r). $
Поэтому при $ r \in \N $ имеет место неравенство
\begin{multline*}
(1 / \mes Q_r) \int_{ Q_r} | f(x) | dx \le
(1 / \mes Q_r) \int_{ B(x_r, \delta_r)} | f(x) | dx = \\
(\mes B(x_r, \delta_r) / \mes Q_r) (1/ \mes B(x_r, \delta_r))
\int_{ B(x_r, \delta_r)} | f(x) | dx \\ = c_5
(1/ \mes B(x_r, \delta_r)) \int_{ B(x_r, \delta_r)} | f(x) | dx \le
c_5 M_f(x_r) \le c_5 \alpha. \square
\end{multline*}

При $ d \in \N $ через $ L(\R^d) $ обозначим пространство измеримых по
Лебегу функций в $ \R^d. $ Как обычно, при $ 1 \le p_1, p_2 \le
\infty $ под суммой $ L_{p_1}(\R^d) +L_{p_2}(\R^d) $ понимается
подпространство в $ L(\R^d), $ состоящее из всех функций $ f \in L(\R^d), $
для которых существуют функции $ f_1 \in L_{p_1}(\R^d) $ и
$ f_2 \in L_{p_2}(\R^d) $ такие, что $ f = f_1 +f_2.$
Напомним, что при $ 1 \le p_1 \le p \le p_2 < \infty $ справедливо включение
$ L_p(\R^d) \subset L_{p_1}(\R^d) +L_{p_2}(\R^d). $

Теорема 2.1.6

Пусть $ d \in \N, 1 < q < \infty, C_0 \in \R_+, C_1 \in \R_+ $ и $
T: (L_1(\R^d) +L_q(\R^d)) \mapsto L(\R^d) $ -- отображение,
удовлетворяющее следующим условиям:

1) для любых $ f, g \in (L_1(\R^d) +L_q(\R^d)) $ почти для всех
$ x \in \R^d $ выполняется неравенство
  \begin{equation*} \tag{2.1.13}
| (T(f +g))(x) | \le | (T f)(x) | +| (T g)(x) |;
\end{equation*}

2) для $ f \in L_1(\R^d) $ при $ \alpha >0 $ соблюдается неравенство
\begin{equation*} \tag{2.1.14}
\mes\{ x \in \R^d: | (T f)(x) | > \alpha \} \le
(C_0 / \alpha) \| f \|_{L_1(\R^d)};
\end{equation*}

3) для $ f \in L_q(\R^d) $ при $ \alpha >0 $ выполняется неравенство
\begin{equation*} \tag{2.1.15}
\mes\{ x \in \R^d: | (T f)(x) | > \alpha \} \le
((C_1 / \alpha) \| f \|_{L_q(\R^d)})^q.
\end{equation*}
Тогда при $ 1< p < q $ существует константа $ c_7(p, q, C_0, C_1)
>0 $ такая, что для любого отображения $ T: (L_1(\R^d) +L_q(\R^d))
\mapsto L(\R^d), $ подчинённого условиям (2.1.13) -- (2.1.15), для
$ f \in L_p(\R^d) $ имеет место неравенство
\begin{equation*} \tag{2.1.16}
\| T f \|_{L_p(\R^d)} \le c_7 \| f \|_{L_p(\R^d)}.
\end{equation*}
\bigskip

2.2. В этом пункте устанавливается аналог теоремы Литтлвуда-Пэли для
операторов $ \{\mathcal E_\kappa^{d,l}, \kappa \in \Z_+^d\}, d \in \N, l \in \Z_+^d $
(см. теорему 2.2.4), а из него выводится утверждение, содержащее оценку,
объявленную в названии параграфа (см. следствие из теоремы 2.2.4). При этом при
доказательстве теоремы 2.2.4 будем придерживаться того же подхода,
что в [3, п. 1.5.2] в случае теоремы Литтлвуда-Пэли для кратных рядов Фурье.
Убедимся, что имеет место

Лемма 2.2.1

Пусть $  l \in \Z_+, 1 < p < \infty. $ Тогда существует константа
$ c_1(l,p) >0 $ такая, что при любом $ k \in \N $
для любого набора чисел $ \sigma = \{ \sigma_\kappa \in \{-1, 1\}:
\kappa =1, \ldots, k \}, $ для $ f \in L_p(I) $
справедливо неравенство
\begin{equation*} \tag{2.2.1}
\| \sum_{\kappa =1}^k \sigma_\kappa \cdot (\mathcal E_\kappa^{1,l} f)
\|_{L_p(I)} \le c_1 \| f \|_{L_p(I)}.
\end{equation*}

Отметим, что доказательство леммы 2.2.1 проводится по схеме, использованной
в [9] при доказательстве теоремы 1 из гл. II.

Доказательство.

Сначала установим справедливость (2.2.1) при $ 1 < p \le 2. $
Принимая во внимание, что для $ f \in L_1(\R) +L_2(\R) $ имеет место
включение $ f \mid_I \in L_1(I), $ определим при $ k \in \N $
отображение $ T = T_{k,\sigma}: L_1(\R) +L_2(\R) \mapsto L(\R), $ полагая для
$ f \in (L_1(\R) +L_2(\R)) $ значение
$$
( T f)(x) = \begin{cases} 0, \text{при} x \in \R \setminus I; \\
 \sum_{\kappa =1}^k \sigma_\kappa (\mathcal E_\kappa^{1,l}( f \mid_I))(x) ,
\text{при} x \in I.
\end{cases}
$$

При $ k \in \N, \sigma = \{ \sigma_\kappa \in \{-1, 1\}:
\kappa =1, \ldots, k \} $ для $ f, g \in (L_1(\R) +L_2(\R)) $ имеем
\begin{multline*} \tag{2.2.2}
| (T (f +g))(x) | =0 = 0 +0  = | (T f)(x) |
+| (T g)(x) |, x \in \R \setminus I; \\
| (T (f +g))(x) | =
| \sum_{\kappa =1}^k \sigma_\kappa (\mathcal E_\kappa^{1,l}( (f +g)
\mid_I))(x) | = \\
| \sum_{\kappa =1}^k \sigma_\kappa (\mathcal E_\kappa^{1,l}( f \mid_I
+g \mid_I))(x)| = \\
| \sum_{\kappa =1}^k (\sigma_\kappa (\mathcal E_\kappa^{1,l} (f \mid_I))(x)
+\sigma_\kappa (\mathcal E_\kappa^{1,l} (g \mid_I))(x))| = \\
| \sum_{\kappa =1}^k \sigma_\kappa (\mathcal E_\kappa^{1,l}( f \mid_I))(x)
+\sum_{\kappa =1}^k \sigma_\kappa (\mathcal E_\kappa^{1,l}( g \mid_I))(x)|
\le \\
| \sum_{\kappa =1}^k \sigma_\kappa (\mathcal E_\kappa^{1,l}( f \mid_I))(x)|
+| \sum_{\kappa =1}^k \sigma_\kappa (\mathcal E_\kappa^{1,l}( g \mid_I))(x)|
= \\
= | (T f)(x) | +| (T g)(x) |, \text{ почти для всех } x \in  I,
\end{multline*}
т.е. выполняется (2.1.13).

Далее, покажем, что при $ k \in \N, \sigma = \{ \sigma_\kappa \in
\{-1, 1\}: \kappa =1, \ldots, k \}, $ для $ f \in L_2(\R) $ и $ \alpha >0 $
соблюдается неравенство
\begin{equation*} \tag{2.2.3}
\mes \{ x \in \R: | (T f)(x) | > \alpha \} \le
((1 / \alpha) \| f \|_{L_2(\R)})^2.
\end{equation*}

В самом деле, при $ k \in \N, \sigma = \{ \sigma_\kappa \in \{-1, 1\}:
\kappa =1, \ldots, k \}, $ для $ f \in L_2(\R) $ и $ \alpha >0 $ ввиду
(1.4.9), (1.5.4) имеем
\begin{multline*}
\int_\R | (T f)(x)|^2 dx = \int_I ( \sum_{\kappa =1}^k \sigma_\kappa
(\mathcal E_\kappa^{1,l}( f \mid_I))(x))^2 dx = \\
\int_I \sum_{\kappa =1}^k \sum_{\kappa^\prime =1}^k
\sigma_\kappa (\mathcal E_\kappa^{1,l}( f \mid_I))(x)
\sigma_{\kappa^\prime} (\mathcal E_{\kappa^\prime}^{1,l}( f \mid_I))(x) dx = \\
\sum_{\kappa =1}^k \sum_{\kappa^\prime =1}^k
\sigma_\kappa \sigma_{\kappa^\prime}
\int_I (\mathcal E_\kappa^{1,l}( f \mid_I))(x)
(\mathcal E_{\kappa^\prime}^{1,l}( f \mid_I))(x) dx \\
= \sum_{\kappa =1}^k \int_I
((\mathcal E_\kappa^{1,l}( f \mid_I))(x))^2 dx \\
= \sum_{ \kappa =1}^k \| \mathcal E_\kappa^{1,l} f \mid_I \|_{L_2(I)}^2
\le \sum_{\kappa \in \Z_+} \| \mathcal E_\kappa^{1,l} f \mid_I \|_{L_2(I)}^2
= \| f \mid_I \|_{L_2(I)}^2 \le \| f \|_{L_2(\R)}^2,
\end{multline*} откуда, как обычно, получаем
\begin{multline*}
\alpha^2 \mes \{ x \in \R: | (T f)(x) | > \alpha \} = \int_{ \{
x \in \R: | (T f)(x) | > \alpha \} } \alpha^2 dx \\
\le \int_{ \{
x \in \R: | (T f)(x) | > \alpha \} } | (T f )(x)|^2 dx \le
\int_\R | (T f)(x)|^2 dx \le \| f \|_{L_2(\R)}^2,
\end{multline*} и, значит, верно (2.2.3).

Теперь установим, что существует константа $ C_0(l) >0 $ такая,
что при $ k \in \N, \sigma = \{ \sigma_\kappa \in \{-1, 1\}:
\kappa =1, \ldots, k \} $ для $ f \in L_1(\R) $ и $ \alpha >0 $
соблюдается неравенство
\begin{equation*} \tag{2.2.4}
\mes \{ x \in \R: | (T f)(x) | > \alpha \} \le
(C_0 / \alpha) \| f \|_{L_1(\R)}.
\end{equation*}

Пусть $ k \in \N, \sigma = \{ \sigma_\kappa \in \{-1, 1\}:
\kappa =1, \ldots, k \}, f \in L_1(\R) $ и $ \alpha >0. $
Для функции $ f $ и числа $ \alpha $ построим замкнутое множество
$ F,$ множество $ W = \R \setminus F $ и семейство интервалов
$ \{ Q_r, r \in \N \}, $ для которых соблюдаются условия
(2.1.4) -- (2.1.7) и (2.1.10) -- (2.1.12) при $ d =1. $
Обозначая через $ \chi_A $ характеристическую функцию
множества $ A \subset \R, $ определим функции $ g \in L_1(\R) \cap
L_2(\R) $ и $ h \in L_1(\R), $
полагая
$$
g(x) = f(x) \chi_F(x) +\sum_{r=1}^\infty (1/ \mes Q_r) (\int_{Q_r} f(y) dy)
\chi_{Q_r}(x), x \in \R,
$$
и $ h = f -g. $

Из (2.1.6), (2.1.5) с учетом того, что $ \mes (W \setminus (\cup_{
r \in \N } Q_r)) =0 $ (ибо $ (W \setminus (\cup_{ r \in \N } Q_r))
\subset \cup_{ r \in \N }( \overline Q_r \setminus Q_r)), $
вытекает, что почти для всех $ x \in \R $ имеет место равенство $
\chi_W (x) = \sum_{r =1}^\infty \chi_{Q_r} (x). $ Поэтому почти
для всех $ x \in \R $ получаем
\begin{multline*}
h(x) = f(x) -g(x) = f(x) \chi_F (x) +f(x) \chi_W (x) -g(x) \\
=f(x) \chi_F (x) +f(x) (\sum_{r =1}^\infty \chi_{Q_r} (x)) -g(x)
\\ = f(x) (\sum_{r =1}^\infty \chi_{Q_r} (x)) -\sum_{r=1}^\infty (1/
\mes Q_r) (\int_{Q_r} f(y) dy) \chi_{Q_r}(x)\\
 = \sum_{r=1}^\infty
( f(x) -(1/ \mes Q_r) \int_{Q_r} f(y) dy) \chi_{Q_r}(x) =
\sum_{r=1}^\infty h_r (x), \end{multline*}
где $ h_r (x) = ( f(x)
-(1/ \mes Q_r) \int_{Q_r} f(y) dy) \chi_{Q_r}(x). $

 Учитывая (2.1.5), (2.1.6), на основании (2.1.10) и (2.1.12) заключаем,
что почти для всех $ x \in \R $ выполняется неравенство
$ | g(x) | \le c_2 \alpha, $ из которого вытекает оценка
\begin{multline*} \tag{2.2.5}
\| g \|_{L_2(\R)}^2 = \int_\R | g(x) |^2 dx \le \int_\R c_2 \alpha
| g(x) | dx = \\
c_2 \alpha \int_\R | f(x) \chi_F(x) +\sum_{r=1}^\infty
(1/ \mes Q_r) (\int_{Q_r} f(y) dy) \chi_{Q_r}(x) | dx \le \\
c_2 \alpha \int_\R | f(x) | \chi_F(x) +\sum_{r=1}^\infty
(1/ \mes Q_r) | \int_{Q_r} f(y) dy | \chi_{Q_r}(x)  dx = \\
c_2 \alpha (\int_\R | f(x) | \chi_F(x) dx +\sum_{r=1}^\infty \int_\R
(1/ \mes Q_r) | \int_{Q_r} f(y) dy | \chi_{Q_r}(x)  dx ) = \\
c_2 \alpha (\int_F | f(x) | dx +\sum_{r=1}^\infty
(1/ \mes Q_r) | \int_{Q_r} f(y) dy | \int_\R \chi_{Q_r}(x) dx ) \le \\
c_2 \alpha (\int_F | f(x) | dx +\sum_{r=1}^\infty
(1/ \mes Q_r) (\int_{Q_r} | f(y) | dy) \mes Q_r ) = \\
c_2 \alpha (\int_F | f(x) | dx +\sum_{r=1}^\infty
\int_{Q_r} | f(x) | dx ) = \\
c_2 \alpha \int_{F \cup (\cup_{r=1}^\infty Q_r)} | f(x) | dx  =
c_2 \alpha \int_\R | f(x) | dx = c_2 \alpha \| f \|_{L_1(\R)}.
\end{multline*}

Для получения (2.2.4), фиксируя множество $ A \subset \R:
\mes A =0 $ и для $ x \in \R \setminus A $ ввиду (2.2.2)
имеет место неравенство
$$
| (T f)(x) | = | (T (g +h))(x) | \le | (T g)(x) | +| (T h)(x) |,
$$
видим, что
\begin{multline*}
(\{ x \in \R: | (T f)(x) | > \alpha \} \setminus A) \subset \{ x
\in \R: | (T g)(x)| +| (T h)(x) | > \alpha \} \\
\subset \{ x \in \R: | (T g)(x) | > \alpha /2 \} \cup \{ x \in
\R: | (T h)(x) |
> \alpha /2 \},
\end{multline*} и, значит,
\begin{multline*} \tag{2.2.6}
\mes \{ x \in \R: | (T f)(x) | > \alpha \} = \mes (\{ x \in \R:
| (T f)(x) | > \alpha \} \setminus A ) \\
\le \mes \{ x \in \R: | (T g)(x) | > \alpha /2 \} + \mes \{ x
\in \R: | (T h)(x) | > \alpha /2 \}. \end{multline*}

Из (2.2.3) и (2.2.5) выводим
\begin{multline*} \tag{2.2.7}
\mes \{ x \in \R: | (T g)(x) | > \alpha /2 \} \le ((2 / \alpha)
\| g \|_{L_2(\R)})^2 = (2 / \alpha)^2 \| g \|_{L_2(\R)}^2 \\
\le (2 / \alpha)^2 c_2 \alpha \| f \|_{L_1(\R)} = (c_3 / \alpha)
\| f \|_{L_1(\R)}.
\end{multline*}

Для оценки второго слагаемого в правой части (2.2.6) имеем
\begin{multline*}
\{ x \in \R: | (T h)(x) | > \alpha /2 \} \\
= \{ x \in F: | (T
h)(x) | > \alpha /2 \} \cup \{ x \in W: | (T h)(x) | > \alpha /2
\},
\end{multline*}
а, следовательно,
\begin{multline*} \tag{2.2.8}
\mes \{ x \in \R: | (T h)(x) | > \alpha /2 \} \\
= \mes \{ x \in
F: | (T h)(x) | > \alpha /2 \} + \mes \{ x \in W: | (T h)(x) |
> \alpha /2 \}.
\end{multline*}

Второе слагаемое в правой части (2.2.8) в силу (2.1.11) удовлетворяет
неравенству
\begin{equation*} \tag{2.2.9}
\mes \{ x \in W: | (T h)(x) | > \alpha /2 \} \le
\mes W \le (c_4 / \alpha) \int_\R |f| dx = (c_4 / \alpha) \| f \|_{L_1(\R)}.
\end{equation*}

Учитывая, что $ (T h)(x) =0 $ для $ x \in \R \setminus I, $ находим, что
$$
\{ x \in F: | (T h)(x) | > \alpha /2 \} =
\{ x \in (F \cap I): | (T h)(x) | > \alpha /2 \},
$$
а, значит,
\begin{multline*} \tag{2.2.10}
\mes \{ x \in F: | (T h)(x) | > \alpha /2 \} \\
= \mes \{ x \in
(F \cap I): | (T h)(x) | > \alpha /2 \} \le (2 / \alpha) \| T
h \|_{L_1(F \cap I)}.
\end{multline*}

Для проведения оценки правой части (2.2.10) определим при $ m \in \N $
функцию $ h_m^\prime $ равенством
$$
h_m^\prime = h -\sum_{r=1}^m h_r
$$
и заметим, что вследствие (2.2.2) при $ m \in \N $ почти для всех
$ x \in (F \cap I) $ справедливо неравенство
$$
| (T h)(x) | = | (T (\sum_{r=1}^m h_r +h_m^\prime))(x) | \le
\sum_{r=1}^m | (T h_r)(x) | +| (T h_m^\prime)(x) |,
$$
которое влечет оценку
\begin{multline*} \tag{2.2.11}
\| T h \|_{L_1(F \cap I)} = \int_{F \cap I} | (T h)(x) | dx
\le \int_{F \cap I} (\sum_{r=1}^m | (T h_r)(x) |) +| (T
h_m^\prime)(x) | dx\\
= \sum_{r=1}^m \int_{F \cap I} | (T h_r)(x) | dx +\int_{F \cap
I} | (T h_m^\prime)(x) | dx.
\end{multline*}

 При $ r \in \N $ оценим сверху значения $ | (T h_r)(x) | $
для $ x \in F \cap I. $

Если $ r \in \N $ таково, что $ Q_r \cap I = \emptyset, $ то
$ (T h_r )(x) =0 $ почти для всех $ x \in \R. $

Пусть $ F \cap I \ne \emptyset $ и $ r \in \N: Q_r \cap I \ne \emptyset. $
Тогда ввиду (2.1.4), (1.4.2) либо $ Q_r \subset I, $ либо
$ I \subset Q_r. $
Второе включение невозможно, поскольку если оно верно, то $ F \cap I
\subset I \subset Q_r, $ т.е. $ F \cap Q_r \ne \emptyset, $ что
противоречит (2.1.6) и определению $ W. $ Поэтому в рассматриваемой
ситуации $ Q_r \subset I. $
Отметим, что
$$
I = ( \cup_{ n \in \Nu_{0, 2^k -1}^1} Q_{k, n}^1)  \cup
\{ 2^{-k} \nu:  \nu =1, \ldots, 2^k -1 \},
$$
и, значит,
\begin{multline*}
F \cap I = ( \cup_{ n \in \Nu_{0, 2^k -1}^1} ( F \cap Q_{k, n}^1))
\cup ( F \cap \{ 2^{-k} \nu:  \nu =1, \ldots, 2^k -1 \}) \\
= (\cup_{ n \in \Nu_{0, 2^k -1}^1:  F \cap Q_{k, n}^1 \ne
\emptyset } ( F \cap Q_{k, n}^1))  \cup ( F \cap \{ 2^{-k} \nu:
\nu =1, \ldots, 2^k -1 \}). \end{multline*}

Учитывая это замечание, проведем оценку сверху $ | (T h_r)(x) | $
при $ r \in \N: Q_r \cap I \ne \emptyset, n \in \Nu_{0, 2^k -1}^1:
F \cap Q_{k, n}^1 \ne \emptyset $ для $ x \in F \cap Q_{k, n}^1. $
Фиксируя систему функций
$$
\{ \phi_i = \phi_i^{1,l,\{1\} }, i =1, \ldots, \mathfrak R_1^{1, l} \},
$$
удовлетворяющую условиям леммы 1.4.10 при $ d =1, J = \{1\}, $
согласно (1.4.12), при $ r \in \N: Q_r \cap I \ne \emptyset, n \in
\Nu_{0, 2^k -1}^1: F \cap Q_{k, n}^1 \ne \emptyset, $ почти для всех $ x \in F \cap
Q_{k, n}^1 $ имеем
\begin{multline*}
| (T h_r)(x) | = \biggl| \sum_{\kappa =1}^k \sigma_\kappa
(\mathcal E_\kappa^{1,l}( H_r \mid_I))(x)\biggr|
\\
= \biggl| \sum_{\kappa =1}^k \sigma_\kappa \biggl(\sum_{
\rho_{\kappa -1} \in \Nu_{0, 2^{\kappa -1} -1}^1 }
 \sum_{ i =1, \ldots, \mathfrak R_1^{1, l}} 2^{ \kappa -1} \\
 \times\biggl(\int_{ Q_{\kappa -1, \rho_{ \kappa -1}}^1 } \phi_i( 2^{
\kappa -1} y -\rho_{\kappa -1}) h_r(y) dy\biggr) \phi_i( 2^{
\kappa -1} x
-\rho_{\kappa -1}) \biggr) \biggr| \\
= \biggl| \sum_{\kappa =0}^{k -1} \sigma_{\kappa +1} \biggl(
\sum_{ \rho_\kappa \in \Nu_{0, 2^\kappa -1}^1 }
\sum_{ i =1, \ldots, \mathfrak R_1^{1, l}} 2^\kappa \\
\times\biggl(\int_{ Q_{\kappa , \rho_\kappa}^1 } \phi_i( 2^\kappa
y-\rho_\kappa) h_r(y) dy\biggr) \phi_i( 2^\kappa x -\rho_\kappa) \biggr) \biggr| \\
= \biggl| \sum_{\kappa =0}^{k -1} \sigma_{\kappa +1} \biggl(
\sum_{ \rho_\kappa \in \Nu_{0, 2^\kappa -1}^1: Q_{\kappa,
\rho_\kappa}^1 \cap Q_{k, n}^1 \ne \emptyset }
\sum_{ i =1, \ldots, \mathfrak R_1^{1, l}} 2^\kappa \\
\times\biggl(\int_{ Q_{\kappa , \rho_\kappa}^1 } \phi_i( 2^\kappa
y -\rho_\kappa) h_r(y) dy\biggr) \phi_i( 2^\kappa x -\rho_\kappa)
\biggr) \biggr|.
\end{multline*}

Заметим, что при $ \kappa =0, \ldots, k $ множество
$ \{ \rho_\kappa \in \Nu_{0, 2^\kappa -1}^1:
Q_{\kappa, \rho_\kappa}^1 \cap Q_{k, n}^1 \ne \emptyset \} $
состоит из единственного элемента, который обозначим $ \rho_\kappa(n). $
Непустота этого множества следует из включения $ Q_{k, n}^1 \subset I
\subset \cup_{ \rho_\kappa \in \Nu_{0, 2^\kappa -1}^1}
\overline Q_{\kappa, \rho_\kappa}^1, $
а единственность, в силу (1.4.2), вытекает
из включения $ Q_{k, n}^1 \subset Q_{\kappa, \rho_\kappa}^1 $
для $ \rho_\kappa \in \Nu_{0, 2^\kappa -1}^1:
Q_{\kappa, \rho_\kappa}^1 \cap Q_{k, n}^1 \ne \emptyset. $

Принимая во внимание сказанное, получаем, что при $ r \in \N:
Q_r \cap I \ne \emptyset, n \in \Nu_{0, 2^k -1}^1: F \cap Q_{k, n}^1 \ne
\emptyset $ почти для всех $ x \in F \cap Q_{k, n}^1 $
выполняется соотношение
\begin{multline*} \tag{2.2.12}
| (T h_r)(x) | = \biggl| \sum_{\kappa =0}^{k -1} \sigma_{\kappa
+1} \cdot
( \sum_{ i =1, \ldots, \mathfrak R_1^{1, l}} 2^\kappa \\
\times \biggl(\int_{ Q_{\kappa, \rho_\kappa(n)}^1 } \phi_i(
2^\kappa y -\rho_\kappa(n)) h_r(y) dy\biggr) \phi_i( 2^\kappa x
-\rho_\kappa(n)) ) \biggr| \\
\le \sum_{\kappa =0}^{k -1} | \sigma_{\kappa +1} |
( \sum_{ i =1, \ldots, \mathfrak R_1^{1, l}} | 2^\kappa \\
\times \biggl(\int_{ Q_{\kappa, \rho_\kappa(n)}^1 } \phi_i(
2^\kappa y -\rho_\kappa(n)) h_r(y) dy\biggr) \phi_i( 2^\kappa x
-\rho_\kappa(n)) | ) \\
\le  \sum_{\kappa =0}^{k -1}
\sum_{ i =1, \ldots, \mathfrak R_1^{1, l}}
2^{\kappa} \\
\times \biggl| \int_{ Q_{\kappa, \rho_\kappa(n)}^1 } \phi_i(
2^\kappa y -\rho_\kappa(n)) h_r(y) dy \biggr|
\cdot \| \phi_i \|_{L_\infty(I)} \\
\le c_5(l) \sum_{\kappa =0}^{k -1} \sum_{ i =1, \ldots, \mathfrak
R_1^{1, l}} 2^{\kappa} \biggl| \int_{ Q_{\kappa, \rho_\kappa(n)}^1}
\phi_i(2^\kappa y -\rho_\kappa(n)) h_r(y) dy \biggr| \\
= c_5  \sum_{ i =1, \ldots, \mathfrak R_1^{1, l}} \sum_{\kappa
=0}^{k -1} 2^{\kappa} \biggl| \int_{ Q_{\kappa, \rho_\kappa(n)}^1
} \phi_i( 2^\kappa y -\rho_\kappa(n)) h_r(y) dy \biggr| \\
=  c_5 \sum_{ i =1, \ldots, \mathfrak R_1^{1, l}} \sum_{\kappa
=0}^{k -1} 2^{\kappa} \biggl| \int_{ Q_{\kappa, \rho_\kappa(n)}^1
\cap Q_r} \phi_i( 2^\kappa y -\rho_\kappa(n)) h_r(y) dy \biggr| \\
= c_5 \sum_{ i =1, \ldots, \mathfrak R_1^{1, l}} \sum_{\kappa =0,
\ldots, k -1: Q_{\kappa, \rho_\kappa(n)}^1 \cap Q_r \ne
\emptyset} 2^\kappa \\
\times \biggl| \int_{ Q_{\kappa, \rho_\kappa(n)}^1 \cap Q_r}
\phi_i( 2^\kappa y -\rho_\kappa(n)) h_r(y) dy \biggr|.
\end{multline*}

Для оценки правой части (2.2.12) установим справедливость леммы 2.2.2.

Лемма 2.2.2

Пусть $ r \in \N: Q_r \cap I \ne \emptyset,
n \in \Nu_{0, 2^k -1}^1: F \cap Q_{k, n}^1 \ne \emptyset. $ Положим
$$
\iota = \iota(r,n) = \max \{ \kappa =0, \ldots, k -1:
Q_{\kappa, \rho_\kappa(n)}^1 \cap Q_r \ne \emptyset \}.
$$
Тогда для $ \kappa =0, \ldots, \iota $ существует $ \nu_\kappa =
\nu_\kappa(r, n) \in \{0, 1\} $ такое, что справедливо включение
\begin{equation*} \tag{2.2.13}
Q_r \subset (2^{-(\kappa +1)} (2 \rho_\kappa(n) +\nu_\kappa)
+2^{-(\kappa +1)} I).
\end{equation*}

Доказательство.

Прежде всего отметим, что при $ \kappa =0, \ldots, k -1 $ в силу (1.4.2)
$ Q_{\kappa +1, \rho_{\kappa +1}(n)}^1 \subset
Q_{\kappa, \rho_\kappa(n)}^1, $
и, следовательно,
$$
2^{-\kappa} \rho_\kappa(n) \le
2^{-(\kappa +1)} \rho_{\kappa +1}(n) <
2^{-(\kappa +1)} \rho_{\kappa +1}(n) +2^{-(\kappa +1)} \le
2^{-\kappa} \rho_\kappa(n) +2^{-\kappa},
$$
или
$$
2 \rho_\kappa(n) \le \rho_{\kappa +1}(n) < 2 \rho_\kappa(n) +2,
$$
т.е.
$$
2 \rho_\kappa(n) \le \rho_{\kappa +1}(n) \le 2 \rho_\kappa(n) +1,
$$
а, значит, существует $ \epsilon_\kappa(n) \in \{0, 1\} $
такое, что соблюдается равенство
\begin{equation*} \tag{2.2.14}
\rho_{\kappa +1}(n) = 2 \rho_\kappa(n) +\epsilon_\kappa(n).
\end{equation*}
ввиду (2.2.14) при $ \kappa =0, \ldots, k -1 $ имеем
\begin{multline*} \tag{2.2.15}
Q_{\kappa, \rho_\kappa(n)}^1 = 2^{-\kappa} \rho_\kappa(n) +2^{-\kappa} I = \\
2^{-\kappa} \rho_\kappa(n) +2^{-\kappa} ((2^{-1} I) \cup (2^{-1} +2^{-1} I )
\cup \{2^{-1} \}) = \\
2^{-\kappa} \rho_\kappa(n) +((2^{-(\kappa +1)} I) \cup (2^{-(\kappa +1)}
+2^{-(\kappa +1)} I ) \cup \{2^{-(\kappa +1)} \}) = \\
(2^{-(\kappa +1)} 2 \rho_\kappa(n)
+2^{-(\kappa +1)} I) \cup
(2^{-(\kappa +1)} (2 \rho_\kappa(n) +1)
+2^{-(\kappa +1)} I ) \cup \{2^{-(\kappa +1)} (2 \rho_\kappa(n) +1)\} = \\
(2^{-(\kappa +1)} (2 \rho_\kappa(n) +\epsilon_\kappa(n))
+2^{-(\kappa +1)} I) \cup\\
\cup (2^{-(\kappa +1)} (2
\rho_\kappa(n) +1 -\epsilon_\kappa(n))
+2^{-(\kappa +1)} I ) \cup \{2^{-(\kappa +1)} (2 \rho_\kappa(n) +1)\} = \\
Q_{\kappa +1, \rho_{\kappa +1}(n)}^1 \cup
(2^{-(\kappa +1)} (2 \rho_\kappa(n) +1 -\epsilon_\kappa(n))
+2^{-(\kappa +1)} I ) \cup \{2^{-(\kappa +1)} (2 \rho_\kappa(n) +1)\}.
\end{multline*}

Далее, при $ r \in \N: Q_r \cap I \ne \emptyset, n \in \Nu_{0, 2^k
-1}^1: F \cap Q_{k, n}^1 \ne \emptyset, \kappa =0, \ldots, \iota $
с учетом (1.4.2) имеем
$$
( Q_{\kappa, \rho_\kappa(n)}^1 \cap Q_r) \supset
( Q_{\iota, \rho_\iota(n)}^1 \cap Q_r) \ne \emptyset
$$
и, следовательно (см. (2.1.4), (1.4.2)), либо
\begin{equation*} \tag{2.2.16}
Q_r \subset Q_{ \kappa, \rho_\kappa(n)}^1,
\end{equation*}
либо
$$
Q_{\kappa, \rho_\kappa(n)}^1 \subset Q_r.
$$
Последнее включение невозможно, ибо если оно справедливо, то
$$
(F \cap Q_r) \supset (F \cap Q_{ \kappa, \rho_\kappa(n)}^1)
\supset (F \cap Q_{k, n}^1)
\ne \emptyset,
$$
т.е. $ F \cap Q_r \ne \emptyset, $ что неверно.
Итак, при $ r \in \N: Q_r \cap I \ne \emptyset,
n \in \Nu_{0, 2^k -1}^1: F \cap Q_{k, n}^1 \ne \emptyset,
\kappa =0, \ldots, \iota $ соблюдается (2.2.16).

Теперь при $ r \in \N: Q_r \cap I \ne \emptyset,
n \in \Nu_{0, 2^k -1}^1: F \cap Q_{k, n}^1 \ne \emptyset,
\kappa =0, \ldots, \iota -1 $ из (2.2.16) с $ \kappa +1 $ вмество
$ \kappa $ и (2.2.14) следует (2.2.13).

Далее, при $ r \in \N: Q_r \cap I \ne \emptyset, n \in \Nu_{0, 2^k
-1}^1: F \cap Q_{k, n}^1 \ne \emptyset, \kappa = \iota < k -1 $
справедливо соотношение $ Q_{\kappa +1, \rho_{\kappa +1}(n)}^1
\cap Q_r = \emptyset, $ что в соединении с (2.2.16), (2.2.15) дает
включение
\begin{equation*} \tag{2.2.17}
Q_r \subset (2^{-(\kappa +1)} (2 \rho_\kappa(n) +1 -\epsilon_\kappa(n))
+2^{-(\kappa +1)} I ) \cup \{2^{-(\kappa +1)} (2 \rho_\kappa(n) +1)\},
 \end{equation*}
которое влечет (2.2.13), поскольку $ (2^{-(\kappa +1)} (2
\rho_\kappa(n) +1)) \notin Q_r, $ ибо если $ (2^{-(\kappa +1)} (2
\rho_\kappa(n) +1)) \in Q_r, $ то в силу открытости $ Q_r $ с
учетом (2.2.14) будет $ Q_{\kappa +1, \rho_{\kappa +1}(n)}^1 \cap
Q_r \ne \emptyset, $ что неверно.

Наконец, при $ r \in \N: Q_r \cap I \ne \emptyset, n \in \Nu_{0,
2^k -1}^1: F \cap Q_{k, n}^1 \ne \emptyset, \kappa = \iota = k -1,
$ если $ Q_{\kappa +1, \rho_{\kappa +1}(n)}^1 \cap Q_r = Q_{k,
n}^1 \cap Q_r \ne \emptyset, $ то, как и при выводе (2.2.16),
получаем, что $ Q_r \subset Q_{ \kappa +1, \rho_{\kappa +1}(n)}^1,
$ т.е. ввиду (2.2.14) имеет место (2.2.13). Если же при $ r \in
\N: Q_r \cap I \ne \emptyset, n \in \Nu_{0, 2^k -1}^1: F \cap
Q_{k, n}^1 \ne \emptyset, \kappa = \iota = k -1 $ пересечение $
Q_{\kappa +1, \rho_{\kappa +1}(n)}^1 \cap Q_r = \emptyset, $ то из
(2.2.16) и (2.2.15) вытекает (2.2.17), которое, как и выше, влечет
(2.2.13). $ \square $

Отметим еще, что при $ r \in \N $ выполняется равенство
\begin{multline*} \tag{2.2.18}
\int_{Q_r} h_r (x) dx = \int_{Q_r} ( f(x) -(1/ \mes Q_r)
\int_{Q_r} f(y) dy) \chi_{Q_r}(x) dx\\ = \int_{Q_r} ( f(x) -(1/
\mes Q_r) \int_{Q_r} f(y) dy) dx \\
= \int_{Q_r} f(x) dx -(1/ \mes Q_r) (\int_{Q_r} f(y) dy)
\int_{Q_r} dx \\
= \int_{Q_r} f(x) dx -\int_{Q_r} f(y) dy =0.
\end{multline*}

Возвращаясь к (2.2.12), благодаря (2.2.16), при $ r \in \N: Q_r
\cap I \ne \emptyset, n \in \Nu_{0, 2^k -1}^1: F \cap Q_{k, n}^1 \ne
\emptyset, $ почти для всех
$ x \in F \cap Q_{k, n}^1 $ получаем неравенство
\begin{equation*} \tag{2.2.19}
| (T h_r)(x) | \le
c_5 \sum_{ i =1, \ldots, \mathfrak R_1^{1, l}}
\sum_{\kappa =0}^\iota 2^\kappa | \int_{ Q_r}
\phi_i( 2^\kappa y -\rho_\kappa(n)) h_r(y) dy |.
\end{equation*}

Фиксируя для каждого $ r \in \N $ точку $ y_r \in Q_r, $ при $ r
\in \N: Q_r \cap I \ne \emptyset, n \in \Nu_{0, 2^k -1}^1: F \cap
Q_{k, n}^1 \ne \emptyset, \kappa =0, \ldots, \iota, i =1, \ldots,
\mathfrak R_1^{1, l}, $ с учетом (2.2.18) выводим
\begin{multline*} \tag{2.2.20}
| \int_{ Q_r} \phi_i( 2^\kappa y -\rho_\kappa(n)) h_r(y) dy | = \\
| \int_{ Q_r} (\phi_i( 2^\kappa y -\rho_\kappa(n)) -
\phi_i( 2^\kappa y_r -\rho_\kappa(n)))   h_r(y) +
\phi_i( 2^\kappa y_r -\rho_\kappa(n)) h_r(y) dy | = \\
| \int_{ Q_r} (\phi_i( 2^\kappa y -\rho_\kappa(n)) -
\phi_i( 2^\kappa y_r -\rho_\kappa(n)))   h_r(y) dy
+\phi_i( 2^\kappa y_r -\rho_\kappa(n)) \int_{Q_r} h_r(y) dy | = \\
| \int_{ Q_r} (\phi_i( 2^\kappa y -\rho_\kappa(n)) -
\phi_i( 2^\kappa y_r -\rho_\kappa(n)))   h_r(y) dy | \le \\
\int_{ Q_r} | \phi_i( 2^\kappa y -\rho_\kappa(n)) -
\phi_i( 2^\kappa y_r -\rho_\kappa(n))|
\cdot | h_r(y) | dy.
\end{multline*}

Принимая во внимание (2.2.13), (1.4.13), при $ r \in \N: Q_r \cap I \ne
\emptyset, n \in \Nu_{0, 2^k -1}^1: F \cap Q_{k, n}^1 \ne \emptyset,
\kappa =0, \ldots, \iota, i =1, \ldots, \mathfrak R_1^{1, l}, $ для
$ y \in Q_r $ имеем
\begin{multline*} \tag{2.2.21}
| \phi_i( 2^\kappa y -\rho_\kappa(n)) -
\phi_i( 2^\kappa y_r -\rho_\kappa(n))| = \\
| (y -y_r) \frac{d} {dz} (\phi_i( 2^\kappa z -\rho_\kappa(n)))
\mid_{z =\theta y +(1 -\theta) y_r} | \le \\
| y -y_r | \sup_{z \in Q_r} | \frac{d} {dz} (\phi_i( 2^\kappa z
-\rho_\kappa(n))) | = \\
| y -y_r | \sup_{z \in Q_r} | 2^\kappa \frac{d \phi_i} {du} ( 2^\kappa z
-\rho_\kappa(n)) | \le \\
2^\kappa ( \diam Q_r ) \sup_{z \in Q_r} | \frac{d \phi_i} {du} ( 2^\kappa z
-\rho_\kappa(n)) | \le \\
2^\kappa ( \diam Q_r ) \sup_{z \in (2^{-(\kappa +1)} (2 \rho_\kappa(n)
+\nu_\kappa) +2^{-(\kappa +1)} I)} | \frac{d \phi_i} {du} ( 2^\kappa z
-\rho_\kappa(n)) | = \\
2^\kappa ( \diam Q_r ) \sup_{u \in (2^{-1} \nu_\kappa +2^{-1} I)}
| \frac{d \phi_i} {du} (u) | = \\
2^\kappa ( \diam Q_r ) \sup_{u \in Q_{1, \nu_\kappa}^1} | \frac{d
\phi_i} {du} (u) | \le  c_6(l) 2^\kappa ( \diam Q_r ).
\end{multline*}

Соединяя (2.2.20) и (2.2.21), находим, что
при $ r \in \N: Q_r \cap I \ne \emptyset,
n \in \Nu_{0, 2^k -1}^1: F \cap Q_{k, n}^1 \ne \emptyset,
\kappa =0, \ldots, \iota, i =1, \ldots, \mathfrak R_1^{1, l} $ выполняется неравенство
\begin{multline*}
| \int_{ Q_r} \phi_i( 2^\kappa y -\rho_\kappa(n)) h_r(y) dy |
\\ \le \int_{ Q_r} c_6 2^\kappa ( \diam Q_r ) | h_r(y) | dy = c_6
2^\kappa ( \diam Q_r ) \int_{ Q_r} | h_r(y) | dy.
\end{multline*}

Подставляя последнюю оценку в (2.2.19), получаем, что
при $ r \in \N: Q_r \cap I \ne \emptyset,
n \in \Nu_{0, 2^k -1}^1: F \cap Q_{k, n}^1 \ne \emptyset, $ почти для всех
$ x \in F \cap Q_{k, n}^1 $ соблюдается неравенство
\begin{multline*} \tag{2.2.22}
| (T h_r)(x) | \le
c_5 \sum_{ i =1, \ldots, \mathfrak R_1^{1, l}}
\sum_{\kappa =0}^\iota 2^\kappa c_6 2^\kappa ( \diam Q_r )
\int_{ Q_r} | h_r(y) | dy = \\
c_7(l) ( \diam Q_r )
(\int_{ Q_r} | h_r(y) | dy) \sum_{\kappa =0}^\iota 2^{2 \kappa} = \\
c_7 ( \diam Q_r ) (\int_{ Q_r} | h_r(y) | dy) 2^{2 \iota}
\sum_{\kappa =0}^\iota 2^{-2(\iota -\kappa)} \le \\
c_7 ( \diam Q_r ) (\int_{ Q_r} | h_r(y) | dy) 2^{2 \iota}
\sum_{s=0}^\infty 2^{-2s} =  c_8 ( \diam Q_r ) 2^{2 \iota} \int_{
Q_r} | h_r(y) | dy.
\end{multline*}

Заметим, что при $ r \in \N $ в силу (2.1.12) верно неравенство
\begin{multline*} \tag{2.2.23}
\int_{Q_r} | h_r (y) | dy = \int_{Q_r} | ( f(y) -(1/ \mes Q_r)
\int_{Q_r} f(z) dz) \chi_{Q_r}(y) | dy = \\
\int_{Q_r} | f(y) -(1/ \mes Q_r)
\int_{Q_r} f(z) dz | dy \le \\
\int_{Q_r} | f(y) | dy +(1/ \mes Q_r)
| \int_{Q_r} f(z) dz | \int_{Q_r} dy
\le \\
\int_{Q_r} | f(y) | dy +\int_{Q_r} | f(z) | dz = 2 \int_{Q_r} | f(y) | dy
\le c_9 \alpha \mes Q_r.
\end{multline*}

Кроме того, благодаря (2.1.7), при $ r \in \N, $ для $ y \in Q_r $
справедливо неравенство
\begin{equation*} \tag{2.2.24}
\diam Q_r \le c_{10} \rho(y, F).
\end{equation*}

А также при $ r \in \N: Q_r \cap I \ne \emptyset,
n \in \Nu_{0, 2^k -1}^1: F \cap Q_{k, n}^1 \ne \emptyset, $ для
$ x \in F \cap Q_{k, n}^1 \subset Q_{\iota, \rho_\iota(n)}^1 $
и $ y \in Q_r \subset Q_{\iota, \rho_\iota(n)}^1 $ (см. (2.2.16))
соблюдается неравенство
$ | x -y| \le 2^{-\iota} $ или
\begin{equation*} \tag{2.2.25}
2^\iota \le | x -y |^{-1}.
\end{equation*}

Подставляя (2.2.23) -- (2.2.25) в (2.2.22), заключаем, что
при $ r \in \N: Q_r \cap I \ne \emptyset,
n \in \Nu_{0, 2^k -1}^1: F \cap Q_{k, n}^1 \ne \emptyset, $ почти для всех
$ x \in F \cap Q_{k, n}^1 $ имеет место неравенство
\begin{multline*} \tag{2.2.26}
| (T h_r)(x) | \le c_{11} ( \diam Q_r ) 2^{2 \iota} \alpha \mes
Q_r = c_{11} \alpha ( \diam Q_r ) 2^{2 \iota} \int_{ Q_r } dy \\
= c_{11} \alpha \int_{ Q_r } ( \diam Q_r ) 2^{2 \iota} dy \le
c_{12} \alpha \int_{ Q_r } \rho(y, F) | x -y |^{-2} dy.
\end{multline*}

Из сказанного ясно, что при $ r \in \N $
неравенство (2.2.26) выполняется почти для всех
$ x \in F \cap I, $
а, следовательно,
\begin{multline*} \tag{2.2.27}
\int_{ F \cap I } | (T h_r)(x) | dx \le \int_{ F \cap I } c_{12}
\alpha (\int_{ Q_r } \rho(y, F) | x -y |^{-2} dy ) dx\\
\le c_{12} \alpha \int_F \int_{ Q_r } \rho(y, F) | x -y |^{-2} dy
dx.
\end{multline*}

Подставляя (2.2.27) в (2.2.11) и с учетом (2.1.5), (2.1.6)
применяя (2.1.3), а затем (2.1.11), приходим к выводу, что при $ m
\in \N $ справедливо неравенство
\begin{multline*} \tag{2.2.28}
\| T h \|_{L_1(F \cap I)} \le
c_{12} \alpha \sum_{r=1}^m \int_F \int_{ Q_r } \rho(y, F) | x -y |^{-2} dy dx
+ \\ +\int_{F \cap I} | (T h_m^\prime)(x) | dx \le
c_{12} \alpha \int_F (\sum_{r=1}^m \int_{ Q_r } \rho(y, F) | x -y |^{-2} dy) dx
 + \\ +\int_{I} | (T h_m^\prime)(x) | dx =
c_{12} \alpha \int_F \int_{ \cup_{r=1}^m Q_r } \rho(y, F) | x -y |^{-2} dy dx
+ \\ +\int_{I} | (T h_m^\prime)(x) | dx \le
c_{12} \alpha \int_F \int_W \rho(y, F) | x -y |^{-2} dy dx + \\
+\int_{I} | (T h_m^\prime)(x) | dx \le
c_{13} \alpha \mes (\R \setminus F)
+\int_{I} | (T h_m^\prime)(x) | dx = \\
c_{13} \alpha \mes W +\int_{I} | (T h_m^\prime)(x) | dx
\le c_{14} \int_{\R} |f(x)| dx
+\int_{I} | (T h_m^\prime)(x) | dx = \\
c_{14}  \| f \|_{L_1(\R)}
+\int_I | (T h_m^\prime)(x) | dx.
\end{multline*}

Принимая во внимание, что при $ m \in \N $ в силу
оценок (1.5.2) и (2.2.23), условия (2.1.5) и включения
$ f \in L_1(\R) $ имеет место соотношение
\begin{multline*}
\int_I | (T h_m^\prime)(x) | dx \\
= \int_I | \sum_{\kappa =1}^k \sigma_\kappa \cdot (\mathcal E_\kappa^{1,l}
( h_m^\prime \mid_I))(x) | dx \le \int_I \sum_{\kappa =1}^k |
(\mathcal E_\kappa^{1,l}( h_m^\prime \mid_I))(x) | dx \\=
\sum_{\kappa =1}^k \int_I | (\mathcal E_\kappa^{1,l}( h_m^\prime
\mid_I))(x) | dx = \sum_{\kappa =1}^k \| \mathcal E_\kappa^{1,l}(
h_m^\prime
\mid_I) \|_{L_1(I)} \\
\le \sum_{\kappa =1}^k c_{15} \| h_m^\prime \|_{L_1(I)} = c_{15} k
\| h_m^\prime \|_{L_1(I)} \le
c_{15} k \| h_m^\prime \|_{L_1(\R)} = \\
c_{15} k \| \sum_{r=m+1}^\infty h_r \|_{L_1(\R)} \le
c_{15} k \sum_{r=m+1}^\infty \| h_r \|_{L_1(\R)} = \\
c_{15} k \sum_{r=m+1}^\infty \int_\R | h_r (x) | dx =
c_{15} k \sum_{r=m+1}^\infty \int_{ Q_r} | h_r (x) | dx \le \\
c_{15} k \sum_{r=m+1}^\infty 2 \int_{ Q_r} | f(x) | dx \to 0
\text{ при } m \to \infty,
\end{multline*}
из (2.2.28) следует, что выполняется неравенство
\begin{equation*} \tag{2.2.29}
\| T h \|_{L_1(F \cap I)} \le c_{14}  \| f \|_{L_1(\R)}.
\end{equation*}

Подставляя (2.2.29) в (2.2.10), выводим неравенство
$$
\mes \{ x \in F: | (T h)(x) | > \alpha /2 \} \le
(c_{16} / \alpha) \| f \|_{L_1(\R)},
$$
которое в соединении с (2.2.8), (2.2.9) дает оценку
\begin{equation*} \tag{2.2.30}
\mes \{ x \in \R: | (T h)(x) | > \alpha /2 \} \le
(c_{17} / \alpha) \| f \|_{L_1(\R)}.
\end{equation*}

Объединяя (2.2.6), (2.2.7) и (2.2.30), приходим к (2.2.4).
Сопоставляя (2.2.2) -- (2.2.4) с (2.1.13) -- (2.1.15), заключаем, что
при $ k \in \N, \sigma = \{ \sigma_\kappa \in \{-1, 1\}:
\kappa =1, \ldots, k \}, 1 < p < 2 $ для $ f \in L_p(\R) $ согласно (2.1.16)
верно неравенство
$$
\| T f \|_{L_p(\R)} \le c_1 \| f \|_{L_p(\R)},
$$
из которого следует (2.2.1) при $ 1 < p < 2. $ Справедливость (2.2.1)
при $ p =2 $ установлена при выводе (2.2.3).

Теперь проверим соблюдение (2.2.1)
при $ 2 < p < \infty. $
В самом деле, при $ k \in \N, \sigma = \{ \sigma_\kappa \in \{-1, 1\}:
\kappa =1, \ldots, k \}, 2 < p < \infty $ для $ f \in L_p(I), $ используя
 неравенство Гельдера, (1.4.8) и (2.2.1) при $ p^\prime $ вместо $ p $
( $ p^\prime = p/(p-1) \in (1,2)), $ имеем
\begin{multline*}
\| \sum_{\kappa =1}^k \sigma_\kappa \cdot (\mathcal E_\kappa^{1,l}
f) \|_{L_p(I)} = \sup_{ g \in B(L_{p^\prime}(I))} \int_I
(\sum_{\kappa =1}^k \sigma_\kappa (\mathcal E_\kappa^{1,l} f ))
\cdot g dx \\
= \sup_{ g \in B(L_{p^\prime}(I))} \int_I \sum_{\kappa =1}^k
(\sigma_\kappa (\mathcal E_\kappa^{1,l} f ) \cdot g) dx = \sup_{ g
\in B(L_{p^\prime}(I))} \sum_{\kappa =1}^k \sigma_\kappa \int_I
(\mathcal E_\kappa^{1,l} f ) \cdot g dx \\
= \sup_{ g \in B(L_{p^\prime}(I))} \sum_{\kappa =1}^k
\sigma_\kappa \int_I f \cdot (\mathcal E_\kappa^{1,l} g) dx =
\sup_{ g \in B(L_{p^\prime}(I))} \int_I \sum_{\kappa =1}^k
(\sigma_\kappa f \cdot (\mathcal E_\kappa^{1,l} g)) dx \\
= \sup_{ g \in B(L_{p^\prime}(I))} \int_I f \cdot (\sum_{\kappa
=1}^k \sigma_\kappa (\mathcal E_\kappa^{1,l} g)) dx \le \sup_{ g
\in B(L_{p^\prime}(I))} \| f \|_{L_p(I)} \cdot \| \sum_{\kappa
=1}^k \sigma_\kappa (\mathcal E_\kappa^{1,l} g)
\|_{L_{p^\prime}(I)} \\
= \| f \|_{L_p(I)} \cdot \sup_{ g \in B(L_{p^\prime}(I))} \|
\sum_{\kappa =1}^k \sigma_\kappa (\mathcal E_\kappa^{1,l} g)
\|_{L_{p^\prime}(I)} \\
\le  \| f \|_{L_p(I)} \cdot \sup_{ g \in B(L_{p^\prime}(I))}
c_1(l,p^\prime) \| g \|_{L_{p^\prime}(I)} =
c_1 \| f \|_{L_p(I)}.   \square
\end{multline*}

Следствие

В условиях леммы 2.2.1 существует константа $ c_{18}(l,p) >0 $
такая, что при любом $ k \in \Z_+ $
для любого набора чисел $ \{ \sigma_\kappa \in \{-1, 1\}:
\kappa =0, \ldots, k \}, $ для $ f \in L_p(I) $
имеет место неравенство
\begin{equation*} \tag{2.2.31}
\| \sum_{\kappa =0}^k \sigma_\kappa (\mathcal E_\kappa^{1,l} f) \|_{L_p(I)}
\le c_{18} \| f \|_{L_p(I)}.
\end{equation*}

Доказательство.

В самом деле, в условиях леммы 2.2.1, используя (1.5.2) и (2.2.1), выводим
\begin{multline*}
\| \sum_{\kappa =0}^k \sigma_\kappa (\mathcal E_\kappa^{1,l} f)
\|_{L_p(I)} \le \| \sigma_0 (\mathcal E_0^{1,l} f) \|_{L_p(I)} +
\| \sum_{\kappa =1}^k \sigma_\kappa (\mathcal E_\kappa^{1,l} f)
\|_{L_p(I)} \\
= \| \mathcal E_0^{1,l} f \|_{L_p(I)} + \|
\sum_{\kappa =1}^k \sigma_\kappa (\mathcal E_\kappa^{1,l} f)
\|_{L_p(I)} \le c_{18} \| f \|_{L_p(I)}. \square
\end{multline*}

На основании леммы 2.2.1 устанавливается

Теорема 2.2.3

Пусть $ d \in \N, l \in \Z_+^d, 1 < p < \infty. $ Тогда существует константа
$ c_{19}(d,l,p) >0 $ такая, что для любого семейства чисел
$ \{ \sigma_\kappa: \kappa \in \Z_+^d \} $ вида $ \sigma_\kappa = \prod_{j=1}^d
\sigma^j_{\kappa_j}, $ где $ \sigma^j_{\kappa_j} \in \{-1, 1\}, \kappa_j
\in \Z_+, j =1,\ldots,d, $ для $ f \in L_p(I^d) $
справедливо неравенство
\begin{equation*} \tag{2.2.32}
\| \sum_{\kappa \in \Z_+^d} \sigma_\kappa \cdot
(\mathcal E_\kappa^{d,l} f ) \|_{L_p(I^d)}
\le c_{19} \| f \|_{L_p(I^d)}.
\end{equation*}

Доказательство.

Сначала покажем, что в условиях теоремы для любого непустого
множества $ J \subset \{1, \ldots, d \} $ при любом $ k^J \in (\Z_+^d)^J $
и любых наборах чисел $ \{ \sigma^j_{\kappa_j} \in \{-1, 1\}, \kappa_j
= 0,\ldots, k_j \}, j \in J, $
для $ f \in L_p(I^d) $ имеет место неравенство
\begin{multline*} \tag{2.2.33}
\| \sum_{\kappa^J \in \Z_+^m(k^J)}  (\prod_{j \in J} (\sigma^j_{\kappa_j}
V_j (\mathcal E_{\kappa_j}^{1,l_j})))f \|_{L_p(I^d)}
\le (\prod_{j \in J} c_{18}(l_j,p)) \cdot
\| f \|_{L_p(I^d)},
\end{multline*}
где $ m = \card J. $

Доказательство (2.2.33) проведем по индукции относительно $ m. $
При $ m =1, $ т.е. для $ j =1, \ldots, d, $ используя п. 2) леммы
1.3.1, теорему Фубини, (1.3.1), (2.2.31),  имеем
\begin{multline*}
\| \sum_{\kappa_j =0}^{k_j} \sigma^j_{\kappa_j} \cdot (V_j
(\mathcal E_{\kappa_j}^{1,l_j})) f \|_{L_p(I^d)}^p = \| (V_j
(\sum_{\kappa_j =0}^{k_j} \sigma^j_{\kappa_j} \mathcal
E_{\kappa_j}^{1,l_j})) f \|_{L_p(I^d)}^p \\
= \int_{I^d} | (V_j (\sum_{\kappa_j =0}^{k_j} \sigma^j_{\kappa_j}
\mathcal E_{\kappa_j}^{1,l_j})) f |^p dx\\ = \int_{I^{d -1}}
\int_I | ((V_j (\sum_{\kappa_j =0}^{k_j} \sigma^j_{\kappa_j}
\mathcal E_{\kappa_j}^{1,l_j})) f) (x_1,\ldots, x_{j -1}, x_j,
x_{j +1}, \ldots, x_d)|^p\\
\times dx_j dx_1 \ldots dx_{j -1} dx_{j +1} \ldots dx_d
\\
= \int_{I^{d -1}} \int_I | ( (\sum_{\kappa_j =0}^{k_j}
\sigma^j_{\kappa_j} \mathcal E_{\kappa_j}^{1,l_j}) f (x_1,\ldots,
x_{j -1}, \cdot, x_{j +1}, \ldots, x_d))(x_j) |^p\\
\times dx_j dx_1 \ldots
dx_{j -1} dx_{j +1} \ldots dx_d\\
 = \int_{I^{d -1}} \|
(\sum_{\kappa_j =0}^{k_j} \sigma^j_{\kappa_j} \mathcal
E_{\kappa_j}^{1,l_j}) f (x_1,\ldots, x_{j -1}, \cdot, x_{j +1},
\ldots, x_d) \|_{L_p(I)}^p \\
\times dx_1 \ldots dx_{j -1} dx_{j +1} \ldots
dx_d \\
\le \int_{I^{d -1}} ( c_{18}(l_j,p) \| f (x_1,\ldots, x_{j -1},
\cdot, x_{j +1}, \ldots, x_d) \|_{L_p(I)})^p\\
\times dx_1 \ldots dx_{j -1} dx_{j +1} \ldots dx_d\\ = (
c_{18}(l_j,p))^p \int_{I^{d -1}} \int_I | f (x_1,\ldots, x_{j -1},
x_j, x_{j +1}, \ldots, x_d) |^p\\
\times dx_j
dx_1 \ldots dx_{j -1} dx_{j +1} \ldots dx_d \\
= ( c_{18}(l_j,p))^p
\int_{I^d} | f(x)|^p dx = ( c_{18}(l_j,p) \| f \|_{L_p(I^d)})^p,
\end{multline*}
откуда
\begin{equation*} \tag{2.2.34}
\| \sum_{\kappa_j =0}^{k_j} \sigma^j_{\kappa_j} (V_j (\mathcal
E_{\kappa_j}^{1,l_j})) f \|_{L_p(I^d)} \le c_{18}(l_j,p) \| f
\|_{L_p(I^d)},
\end{equation*}
что совпадает с (2.2.33) при $ m =1, J = \{j\}. $

Предположим теперь, что при некотором $ m: 1 \le m \le d -1, $
оценка (2.2.33) имеет место для любого множества $ J \subset \{1,
\ldots, d \}: \card J = m, $ при любом $ k^J \in (\Z_+^d)^J, $
любых наборах чисел $ \{ \sigma^j_{\kappa_j} \in \{-1, 1\},
\kappa_j = 0, \ldots, k_j \}, j \in J, $ и любой функции $ f \in
L_p(I^d).$ Покажем, что тогда неравенство (2.2.33) справедливо при
$ m+1 $ вместо $ m $ для любого множества $ \mathcal J \subset
\{1, \ldots, d \} $ вместо $ J, $ у которого $ \card \mathcal J =
m +1, $ при любом $ k^{ \mathcal J} \in (\Z_+^d)^{ \mathcal J}, $
любых наборах чисел $ \{ \sigma^j_{\kappa_j} \in \{-1, 1\},
\kappa_j = 0, \ldots, k_j \}, j \in \mathcal J, $ и любой функции
$ f \in L_p(I^d). $ Представляя множество $ \mathcal J \subset
\{1, \ldots, d \}: \card \mathcal J = m +1, $ в виде $ \mathcal J
= J \cup \{i\}, i \notin J, $ с учетом (1.3.2) в силу (2.2.34) и
предположения индукции получаем
\begin{multline*}
\| \sum_{\kappa^{\mathcal J} \in \Z_+^{m+1}(k^{\mathcal J})}
(\prod_{j \in \mathcal J} (\sigma^j_{\kappa_j} V_j (\mathcal
E_{\kappa_j}^{1,l_j})))f \|_{L_p(I^d)} \\
= \| \sum_{(\kappa_i, \kappa^J): \kappa_i =0, \ldots, k_i,
\kappa^J \in \Z_+^m(k^J)} \sigma^i_{\kappa_i} (V_i (\mathcal
E_{\kappa_i}^{1,l_i})) ((\prod_{j \in J} (\sigma^j_{\kappa_j} V_j
(\mathcal E_{\kappa_j}^{1,l_j})) )f) \|_{L_p(I^d)}\\ = \|
\sum_{\kappa_i =0}^{k_i} \sigma^i_{\kappa_i} (V_i (\mathcal
E_{\kappa_i}^{1,l_i})) (\sum_{\kappa^J \in \Z_+^m(k^J)}  (\prod_{j
\in J} (\sigma^j_{\kappa_j} V_j (\mathcal E_{\kappa_j}^{1,l_j}))
)f) \|_{L_p(I^d)} \\
\le c_{18}(l_i,p) \| \sum_{\kappa^J \in \Z_+^m(k^J)}  (\prod_{j
\in J} (\sigma^j_{\kappa_j} V_j (\mathcal E_{\kappa_j}^{1,l_j})
))f \|_{L_p(I^d)} \\\le c_{18}(l_i,p) (\prod_{j \in J}
c_{18}(l_j,p)) \cdot \| f \|_{L_p(I^d)} = (\prod_{j \in \mathcal
J} c_{18}(l_j,p)) \cdot \| f \|_{L_p(I^d)},
\end{multline*}
что завершает вывод (2.2.33).

В частности, из (2.2.33) при $ m = d $ ввиду (1.4.4) получаем, что
в условиях теоремы при любом $ k \in \Z_+^d $ соблюдается неравенство
\begin{multline*} \tag{2.2.35}
\| \sum_{\kappa \in \Z_+^d(k)} \sigma_\kappa \cdot (\mathcal
E_\kappa^{d,l} f ) \|_{L_p(I^d)} \le c_{19} \| f \|_{L_p(I^d)},
\sigma_\kappa = \prod_{j=1}^d \sigma^j_{\kappa_j}, \\ \text{ где }
 \sigma^j_{\kappa_j} \in \{-1, 1\}, \kappa_j
\in \Z_+, j =1,\ldots,d,  f \in L_p(I^d),
c_{19} = \prod_{j=1}^d c_{18}(l_j,p).
\end{multline*}

Теперь убедимся в справедливости (2.2.32).
Для произвольного семейства чисел $ \{ \sigma_\kappa = \prod_{j=1}^d
\sigma^j_{\kappa_j}: \sigma^j_{\kappa_j} \in \{-1, 1\}, j =1,\ldots,d,
\kappa \in \Z_+^d \}, $ функции $ f \in L_p(I^d) $ рассмотрим последовательность
$$
\{ (\sum_{\kappa \in \Z_+^d(k \e)} \sigma_\kappa \cdot
(\mathcal E_\kappa^{d,l} f )) \in L_p(I^d), k \in \Z_+ \}
$$
и, принимая во внимание (2.2.35), секвенциальную компактность
шара $ B(L_p(I^d)) $ относительно $*$-слабой топологии в пространсте
$ L_p(I^d) = (L_{p^\prime}(I^d))^*, $ выберем
подпоследовательность
$$
\{ \sum_{\kappa \in \Z_+^d(k_n \e)} \sigma_\kappa \cdot
(\mathcal E_\kappa^{d,l} f ): k_n < k_{n+1}, n \in \N \}
$$
и функцию $ F \in L_p(I^d), $ обладающие тем свойством, что для
любой функции $ g \in L_{p^\prime}(I^d) $ выполняется равенство
\begin{equation*} \tag{2.2.36}
\lim_{ n \to \infty} \int_{I^d} (\sum_{\kappa \in \Z_+^d(k_n \e)}
\sigma_\kappa \cdot \\ (\mathcal E_\kappa^{d,l} f )) g dx =
\int_{I^d} F g dx.
\end{equation*}
Заметим, что при любом $ \kappa \in \Z_+^d, $ благодаря (1.4.8),
(2.2.36), (1.4.6), для  $ g \in L_{p^\prime}(I^d) $ имеет место равенство
\begin{multline*}
\int_{I^d} (\mathcal E_\kappa^{d,l} F) \cdot g dx = \int_{I^d} F
\cdot (\mathcal E_\kappa^{d,l} g) dx = \lim_{ n \to \infty}
\int_{I^d} (\sum_{\kappa^\prime \in \Z_+^d(k_n \e)}
\sigma_{\kappa^\prime} \cdot (\mathcal E_{\kappa^\prime}^{d,l} f
)) (\mathcal E_\kappa^{d,l} g) dx \\
= \lim_{ n \to \infty}
\int_{I^d} \mathcal E_\kappa^{d,l}
 (\sum_{\kappa^\prime \in \Z_+^d(k_n \e)}
\sigma_{\kappa^\prime} \cdot (\mathcal E_{\kappa^\prime}^{d,l} f
)) \cdot g dx \\= \lim_{ n \to \infty} \int_{I^d}
(\sum_{\kappa^\prime \in \Z_+^d(k_n \e)} \sigma_{\kappa^\prime}
\cdot \mathcal E_\kappa^{d,l} (\mathcal E_{\kappa^\prime}^{d,l} f
)) g dx  = \int_{I^d} \sigma_{\kappa} \cdot (\mathcal
E_\kappa^{d,l} f) g dx,
\end{multline*}
и, значит,
\begin{equation*} \tag{2.2.37}
\mathcal E_\kappa^{d,l} F =
\sigma_{\kappa} \cdot (\mathcal E_\kappa^{d,l} f).
\end{equation*}
Учитывая (2.2.37), (1.2.1), (1.5.1), заключаем, что
\begin{equation*}
\sum_{\kappa \in \Z_+^d(k)} \sigma_\kappa \cdot (\mathcal
E_\kappa^{d,l} f ) = \sum_{\kappa \in \Z_+^d(k)} (\mathcal
E_\kappa^{d,l} F) = E_k^{d,l} F
\end{equation*}
сходится к $ F $ в $ L_p(I^d) $ при $ \mn(k) \to \infty. $
Поэтому, переходя к пределу при $ \mn(k) \to \infty $ в
неравенстве (2.2.35), приходим к (2.2.32). $ \square $

Следствие

В условиях теоремы 2.2.3 для любого семейства чисел $ \{ \sigma_\kappa:
\kappa \in \Z_+^d \} $ вида $ \sigma_\kappa = \prod_{j=1}^d
\sigma^j_{\kappa_j}, $ где $ \sigma^j_{\kappa_j} \in \{-1, 1\}, \kappa_j
\in \Z_+, j =1,\ldots,d, $ и любой функции $ f \in L_p(I^d) $
соблюдается неравенство
\begin{multline*} \tag{2.2.38}
\| f \|_{L_p(I^d)} \le c_{19}
\| \sum_{\kappa \in \Z_+^d} \sigma_\kappa \cdot
(\mathcal E_\kappa^{d,l} f ) \|_{L_p(I^d)}.
\end{multline*}

Доказательство.

Сначала покажем, что при любом $ k \in \Z_+^d, $ любом наборе чисел
$ \{ \sigma_\kappa = \prod_{j=1}^d \sigma^j_{\kappa_j}: \sigma^j_{\kappa_j} \in \{-1, 1\}, j =1,\ldots,d,
\kappa \in \Z_+^d(k) \}, $ для $ f \in L_p(I^d) $
справедливо неравенство
\begin{multline*} \tag{2.2.39}
\| E_k^{d,l} f \|_{L_p(I^d)} \le c_{19}
\| \sum_{\kappa \in \Z_+^d(k)}
\sigma_\kappa \cdot (\mathcal E_\kappa^{d,l} f ) \|_{L_p(I^d)}.
\end{multline*}

В самом деле, ввиду (1.4.6), (1.2.1) имеем
\begin{multline*}
(\sum_{\kappa \in \Z_+^d(k)} \sigma_\kappa \cdot \mathcal
E_\kappa^{d,l})^2 f = (\sum_{\kappa, \kappa^\prime \in \Z_+^d(k)}
\sigma_\kappa \sigma_{\kappa^\prime} \cdot \mathcal E_\kappa^{d,l}
\mathcal E_{\kappa^\prime}^{d,l}) f\\
 = (\sum_{\kappa \in
\Z_+^d(k)} \sigma_\kappa^2 \cdot \mathcal E_\kappa^{d,l}) f =
\sum_{\kappa \in \Z_+^d(k)} \mathcal E_\kappa^{d,l} f = E_k^{d,l}
f.
\end{multline*}
Откуда, применяя (2.2.35), выводим
\begin{multline*}
\| E_k^{d,l} f \|_{L_p(I^d)} = \| (\sum_{\kappa \in \Z_+^d(k)}
\sigma_\kappa \cdot \mathcal E_\kappa^{d,l})^2 f \|_{L_p(I^d)}\\
 =\| \sum_{\kappa \in \Z_+^d(k)} \sigma_\kappa \cdot \mathcal
E_\kappa^{d,l} (\sum_{\kappa^\prime \in \Z_+^d(k)}
\sigma_{\kappa^\prime} \cdot (\mathcal E_{\kappa^\prime}^{d,l} f))
\|_{L_p(I^d)} \\
\le c_{19} \| \sum_{\kappa^\prime \in \Z_+^d(k)}
\sigma_{\kappa^\prime} \cdot (\mathcal E_{\kappa^\prime}^{d,l} f)
\|_{L_p(I^d)}.
\end{multline*}

Как видно из (2.2.32) и (1.5.1), в неравенстве (2.2.39) можно
перейти к пределу при $ \mn(k) \to \infty $, в результате чего
получим (2.2.38). $ \square $

С помощью теоремы 2.2.3 и следствия из нее, опираясь на схему доказательства
теоремы Литтлвуда-Пэли, приведенную в [3] для
операторов взятия частных сумм кратных рядов Фурье, устанавливается

Теорема 2.2.4

В условиях теоремы 2.2.3 существуют константы $ c_{20}(d,l,p) >0,
c_{21}(d,l,p) >0 $ такие, что для любой функции $ f \in L_p(I^d) $
выполняются неравенства
\begin{equation*} \tag{2.2.40}
c_{20} \| f \|_{L_p(I^d)} \le
( \int_{I^d} (\sum_{\kappa \in \Z_+^d} ((\mathcal E_\kappa^{d,l} f)(x))^2 )^{p/2}
dx)^{1/p} \le c_{21} \| f \|_{L_p(I^d)}.
\end{equation*}

Доказательство.

Рассмотрим систему Радемахера, состоящую из функций
$$
\omega_\kappa(t) = \operatorname{sign} \sin( 2^{\kappa +1} \pi t),
t \in I, \kappa \in \Z_+,
$$
и определим семейство функций $ \omega_\kappa^d, \kappa \in \Z_+^d, $
полагая
$$
\omega_\kappa^d(t) = \prod_{j=1}^d \omega_{\kappa_j}(t_j), t \in I^d.
$$
Как известно (см., например, [3]), существуют константы $ c_{22}(d,p) >0,
c_{23}(d,p) >0 $ такие, что при любом $ k \in \Z_+^d $ для любого набора
чисел $ \{ a_\kappa \in \R, \kappa \in \Z_+^d(k) \} $
имеет место неравенство
\begin{equation*} \tag{2.2.41}
c_{22} ( \sum_{ \kappa \in \Z_+^d(k)} a_\kappa^2 )^{1/2}
\le
\| \sum_{ \kappa \in \Z_+^d(k)} a_\kappa \omega_\kappa^d(\cdot) \|_{L_p(I^d)}
\le
c_{23} ( \sum_{ \kappa \in \Z_+^d(k)} a_\kappa^2 )^{1/2}.
\end{equation*}

Для $ f \in L_p(I^d) $ при $ k \in \Z_+^d, $ используя (2.2.39), теорему
Фубини, (2.2.41), выводим
\begin{multline*} \tag{2.2.42}
\| E_k^{d,l} f \|_{L_p(I^d)}^p = \int_{I^d} \| E_k^{d,l} f
\|_{L_p(I^d)}^p dt \le \int_{I^d} (c_{19} \| \sum_{\kappa \in
\Z_+^d(k)} \omega_\kappa^d(t) \cdot (\mathcal E_\kappa^{d,l} f )
\|_{L_p(I^d)})^p dt\\
 = (c_{19})^p \int_{I^d} \int_{I^d} |
\sum_{\kappa \in \Z_+^d(k)} \omega_\kappa^d(t) \cdot (\mathcal
E_\kappa^{d,l} f )(x) |^p dx dt \\
= (c_{19})^p \int_{I^d}
\int_{I^d} | \sum_{\kappa \in \Z_+^d(k)} \omega_\kappa^d(t) \cdot
(\mathcal
E_\kappa^{d,l} f )(x) |^p dt dx\\
 = (c_{19})^p \int_{I^d} \|
\sum_{\kappa \in \Z_+^d(k)} (\mathcal E_\kappa^{d,l} f )(x) \cdot
\omega_\kappa^d(\cdot) \|_{L_p(I^d)}^p dx \\
\le (c_{19})^p \int_{I^d} (c_{23} (\sum_{\kappa \in \Z_+^d(k)}
((\mathcal E_\kappa^{d,l} f )(x))^2 )^{1/2})^p dx\\
 = (c_{24})^p \int_{I^d}
(\sum_{\kappa \in \Z_+^d(k)} ((\mathcal E_\kappa^{d,l} f )(x))^2
)^{p/2} dx,
\end{multline*}
и, пользуясь (2.2.41), теоремой Фубини, (2.2.35), получаем
\begin{multline*}
\int_{I^d} (\sum_{\kappa \in \Z_+^d(k)} ((\mathcal E_\kappa^{d,l}
f )(x))^2 )^{p/2} dx \le \int_{I^d} (c_{25} \| \sum_{\kappa \in
\Z_+^d(k)} (\mathcal E_\kappa^{d,l} f )(x) \cdot
\omega_\kappa^d(\cdot) \|_{L_p(I^d)})^p dx\\
 = (c_{25})^p \int_{I^d} \int_{I^d} | \sum_{\kappa \in \Z_+^d(k)}
\omega_\kappa^d(t) \cdot (\mathcal E_\kappa^{d,l} f )(x) |^p dt
dx\\
 = (c_{25})^p \int_{I^d} \int_{I^d} | \sum_{\kappa \in
\Z_+^d(k)} \omega_\kappa^d(t) \cdot (\mathcal E_\kappa^{d,l} f
)(x) |^p dx
dt\\
 = (c_{25})^p \int_{I^d} \| \sum_{\kappa \in \Z_+^d(k)}
\omega_\kappa^d(t) \cdot (\mathcal E_\kappa^{d,l} f )
\|_{L_p(I^d)}^p dt \\
\le (c_{25})^p \int_{I^d} (c_{19} \| f \|_{L_p(I^d)})^p dt =
(c_{21})^p \| f \|_{L_p(I^d)}^p,
\end{multline*}
откуда, в частности, имеем
\begin{equation*} \tag{2.2.43}
\int_{I^d} (\sum_{\kappa \in \Z_+^d(n \e)} ((\mathcal
E_\kappa^{d,l} f )(x))^2 )^{p/2} dx \le (c_{21})^p \| f
\|_{L_p(I^d)}^p, n \in \Z_+.
\end{equation*}

Для получения второго неравенства в (2.2.40) достаточно применить
теорему Леви о предельном переходе под знаком интеграла к монотонно
возрастающей последовательности функций
$ \{ (\sum_{\kappa \Z_+^d(n \e)} ((\mathcal E_\kappa^{d,l} f)(x))^2 )^{p/2},
n \in \Z_+\}, $ принимая во внимание (2.2.43) и учитывая, что
почти для всех $ x \in I^d $ предел
\begin{multline*}
\lim_{n \to \infty} (\sum_{\kappa \in \Z_+^d(n \e)} ((\mathcal
E_\kappa^{d,l} f)(x))^2 )^{p/2}\\ =
 (\lim_{n \to \infty} \sum_{\kappa \in \Z_+^d(n \e)}
((\mathcal E_\kappa^{d,l} f)(x))^2 )^{p/2} = (\sum_{\kappa
\in \Z_+^d} ((\mathcal E_\kappa^{d,l} f)(x))^2 )^{p/2}.
\end{multline*}
Последнее равенство справедливо в силу того, что при $ k \in
\Z_+^d $ имеет место соотношение
$$
\sum_{\kappa \in \Z_+^d(\mn(k) \e)} ((\mathcal E_\kappa^{d,l} f)(x))^2 \le
\sum_{\kappa \in \Z_+^d(k)} ((\mathcal E_\kappa^{d,l} f)(x))^2 \le
\sum_{\kappa \in \Z_+^d(\mx(k) \e)} ((\mathcal E_\kappa^{d,l} f)(x))^2, x \in I^d.
$$

ввиду сказанного и на основании (1.5.1) переходя к пределу в (2.2.42) при
$ \mn(k) \to \infty, $ приходим
к первому неравенству в (2.2.40). $ \square $

Из теоремы 2.2.4, используя соображения, приведенные в [2] для аналогичного
утверждения об операторах взятия частных сумм ряда Фурье, выводится

Следствие

Пусть $ d \in \N, l \in \Z_+^d, 1 \le p < \infty. $ Тогда существует
константа $ c_{26}(d,l,p) >0 $ такая, что если для $ f \in L_1(I^d) $ ряд
$ \sum_{\kappa \in \Z_+^d } \| \mathcal E_\kappa^{d,l} f \|_{L_p(I^d)}^{\p} $
сходится, где $ \p = \min(2,p), $ то
$ f \in L_p(I^d) $ и выполняется неравенство
\begin{equation*} \tag{2.2.44}
\| f \|_{L_p(I^d)} \le c_{26} (\sum_{\kappa \in \Z_+^d }
\| \mathcal E_\kappa^{d,l} f \|_{L_p(I^d)}^{\p} )^{1/\p}.
\end{equation*}

Доказательство.

В условиях теоремы, с помощью неравенства треугольника для нормы в
пространстве $ L_{p/2}(I^d) $ при $ p > 2, $ и неравенства
Гельдера с показателем $ p/2 $ при $ 1 < p \le 2, $ из (2.2.42)
выводим
\begin{multline*} \tag{2.2.45}
\| E_k^{d,l} f \|_{L_p(I^d)} \le c_{24} \cdot ( \| \sum_{\kappa
\in \Z_+^d(k)} (\mathcal E_\kappa^{d,l} f )^2 \|_{L_{p/2}(I^d)}
)^{1/2} \le \\
c_{24} \cdot ( \sum_{\kappa \in \Z_+^d(k)} \| (\mathcal
E_\kappa^{d,l} f )^2 \|_{L_{p/2}(I^d)} )^{1/2} = c_{24} \cdot (
\sum_{\kappa \in \Z_+^d(k)} \| \mathcal
E_\kappa^{d,l} f \|_{L_p(I^d)}^2 )^{1/2}, \\
f \in L_p(I^d), k \in
\Z_+^d, p > 2,
\end{multline*}
и
\begin{multline*} \tag{2.2.46}
\| E_k^{d,l} f \|_{L_p(I^d)} \le c_{24} (\int_{I^d} (\sum_{\kappa
\in \Z_+^d(k)} ((\mathcal E_\kappa^{d,l} f )(x))^2 )^{p/2}
dx)^{1/p} \le \\
c_{24} (\int_{I^d} \sum_{\kappa \in \Z_+^d(k)} |(\mathcal
E_\kappa^{d,l} f )(x)|^p dx)^{1/p} = c_{24} ( \sum_{\kappa \in
\Z_+^d(k)} \| \mathcal E_\kappa^{d,l} f
\|_{L_p(I^d)}^p )^{1/p}, \\
f \in L_p(I^d), k \in \Z_+^d, 1 < p \le
2.
\end{multline*}

Заметим, что на самом деле неравенства (2.2.45) и (2.2.46)
имеют место для любой функции $ f \in L_1(I^d).$

Действительно, если $ k \in \Z_+^d $ и $ f \in L_1(I^d), $
то применяя (2.2.45) в случае $ p > 2, $ соответственно (2.2.46) в случае
$ 1 < p \le 2 $, к функции $ E_k^{d,l}f \in L_\infty(I^d) \subset L_p(I^d) $ и
учитывая (1.4.5) и то обстоятельство, что вследствие (1.2.1) и (1.4.6)
соблюдаются равенства
\begin{multline*}
\mathcal E_\kappa^{d,l} E_k^{d,l} f = \mathcal E_\kappa^{d,l}
(\sum_{ \kappa^\prime \in \Z_+^d(k)} \mathcal
E_{\kappa^\prime}^{d,l} f) = \sum_{ \kappa^\prime \in \Z_+^d(k)}
\mathcal E_\kappa^{d,l} \mathcal E_{\kappa^\prime}^{d,l} f \\
= \mathcal E_\kappa^{d,l} f, \text{ при } \kappa \in \Z_+^d(k), f \in
L_1(I^d), k \in \Z_+^d,
\end{multline*}
получаем, что (2.2.45), соответственно (2.2.46), справедливо для $ f \in L_1(I^d). $

Отметим еще, что для любого семейства неотрицательных чисел $
\{a_\kappa \in \R: a_\kappa \ge 0, \kappa \in \Z_+^d \} $
существует предел (конечный или бесконечный)
$$
\lim_{ \mn(k) \to \infty} \sum_{ \kappa \in \Z_+^d(k)} a_\kappa =
\sup_{k \in \Z_+^d} \sum_{ \kappa \in \Z_+^d(k)} a_\kappa,
$$
поскольку для любых $ k, k^\prime \in \Z_+^d: \mn(k) \ge \mx(k^\prime), $ соблюдается неравенство
$$
\sum_{ \kappa \in \Z_+^d(k)} a_\kappa \ge
\sum_{ \kappa \in \Z_+^d(k^\prime)} a_\kappa.
$$

Пусть теперь $ f \in L_1(I^d) $ и ряд $ \sum_{\kappa \in \Z_+^d }
\| \mathcal E_\kappa^{d,l} f \|_{L_p(I^d)}^{\p} $ сходится. Тогда,
с учетом замечаний, на основании (2.2.45), (2.2.46) заключаем, что
при $ k \in \Z_+ $ справедливо неравенство
\begin{equation*} \tag{2.2.47}
\| E_{k \e}^{d,l} f \|_{L_p(I^d)} \le c_{26} (\sum_{\kappa \in \Z_+^d}
\| \mathcal E_\kappa^{d,l} f \|_{L_p(I^d)}^{\p} )^{1/\p} < \infty.
\end{equation*}
Причем, в силу (1.5.1) последовательность функций $ \{ E_{k
\e}^{d,l}f, k \in \Z_+ \} $ сходится к $ f $ в $ L_1(I^d), $ а,
следовательно, $ \{ E_{k \e}^{d,l}f \} $ сходится к $ f $ по мере.
Поэтому существует подпоследовательность $ \{ E_{k_n \e}^{d,l}f, n
\in \N \}, $ которая сходится к $ f $ почти всюду на $ I^d. $
Применяя теорему Фату о предельном переходе под знаком  интеграла
к последовательности функций $ \{ | E_{k_n \e}^{d,l}f |^p, n \in
\N \}, $ почти всюду на $ I^d $ сходящейся к $ | f |^p, $ с учетом
(2.2.47), приходим к выводу, что $ f \in L_p(I^d) $ и соблюдается
(2.2.44) при $ 1 < p < \infty. $ Справедливость (2.2.44) при $ p
=1 $ вытекает из (1.5.3) и неравенства треугольника для нормы в $
L_1(I^d). \square $

В заключение отметим, что утверждения аналогичные тем, что
установлены в этом пункте, получены автором для ортопроекторов на
взаимно ортогональные подпространства, соответствующие кратно-масштабному
анализу типа Добеши.
\bigskip

\centerline{ \S 3. Колмогоровский $n$-поперечник классов функций,}
\centerline{удовлетворяющих смешанным условиям Гельдера }
\bigskip

3.1. В этом пункте определим классы функций, рассматриваемые в
этом параграфе, и приведем сведения, необходимые для вывода
оценки их поперечников.

Для области $ D \subset \R^d $ и векторов $ h \in \R^d $ и $ l \in
\Z_+^d $ через $ D_h^l $ обозначим множество
\begin{equation*}
\begin{split}
&D_h^l = (\ldots (D_{l_d h_d e_d})_{l_{d-1} h_{d-1} e_{d-1}}
\ldots)_{l_1 h_1 e_1}
= \{ x \in D: x +tlh \in D \forall t \in \overline I^d\} =\\
&=\{ x \in D: (x +\sum_{j \in \s(l)} t_j l_j h_j e_j) \in D
\forall t^{\s(l)} \in (\overline I^d)^{\s(l)} \}.
\end{split}
\end{equation*}

Пусть $ d \in \N, D $ -- область в $ \R^d $ и $ 1 \le p \le
\infty. $ Тогда для $ f \in L_p(D), h \in \R^d $ и $ l \in \Z_+^d
$ определим в $ D_h^l $ смешанную разность  функции $ f $ порядка
$ l, $ соответствующую вектору $ h, $ равенством
\begin{equation*}
\begin{split}
&(\Delta_h^l f)(x) = ((\prod_{j=1}^d \Delta_{h_j e_j}^{l_j}) f)(x)
=
((\prod_{j \in \s(l)} \Delta_{h_j e_j}^{l_j}) f)(x) =\\
&= \sum_{k \in \Z_+^d(l)} (-\e)^{l-k} C_l^k f(x+kh),\ x \in D_h^l,
\end{split}
\end{equation*}

где $ C_l^k = \prod_{j=1}^d C_{l_j}^{k_j}, k \in \Z_+^d(l). $

Имея в виду, что для $ f \in L_p(D), l \in \Z_+^d $ и векторов $
h,h^\prime \in \R^d: h^{\s(l)} = (h^\prime)^{\s(l)}, $ соблюдается
соотношение
$$
\| \Delta_h^l f\|_{L_p(D_h^l)} = \| \Delta_{h^\prime}^l
f\|_{L_p(D_{h^\prime}^l)},
$$
определим для функции $ f $ смешанный модуль непрерывности в $
L_p(D) $ порядка $ l $ равенством
$$
\Omega^l (f,t^{\s(l)})_{L_p(D)} = \supvrai_{ \{ h \in \R^d:
h^{\s(l)} \in t^{\s(l)} (B^d)^{\s(l)} \}} \| \Delta_h^l
f\|_{L_p(D_h^l)}, t^{\s(l)} \in (\R_+^d)^{\s(l)}.
$$

Теперь определим классы функций, изучаемые в этом параграфе (см. [10], [11]).

Пусть $ \alpha \in \R_+^d, 1 \le p \le \infty $ и $ D $ --
область в $ \R^d. $ Тогда зададим вектор $ l = l(\alpha) \in \N^d,
$ полагая $ l_j = \min \{m \in \N: \alpha_j < m \}, j =1,\ldots,d,
$ и обозначим через $ S_p^\alpha H(D) (\mathcal S_p^\alpha
\mathcal H(D)) $ множество всех функций $ f \in L_p(D), $
обладающих тем свойством, что для любого непустого множества $ J
\subset \{1,\ldots,d\} $ выполняется неравенство
$$
\sup_{t^J \in
(\R_+^d)^J}\!(t^J)^{-\alpha^J}\!\Omega^{l\chi_J}(f,t^J)_{L_p(D)}
\!=\!\sup_{t^J \in (\R_+^d)^J}\!(\prod_{j \in J} t_j^{-\alpha_j})
\Omega^{l\chi_J}(f, t^{\s(l\chi_J)})_{L_p(D)}\!<\!\infty(\le\!1).
$$

Пусть $ \alpha,p,D $ и $ l =l(\alpha) $ -- те же, что и выше, и $
\theta \in \R: 1 \le \theta < \infty. $ Тогда обозначим через $
S_{p,\theta}^\alpha B(D) (\mathcal S_{p,\theta}^\alpha \mathcal
B(D)) $ множество всех функций $ f \in L_p(D), $ которые для
любого непустого множества $ J \subset \{1,\ldots,d\} $
удовлетворяют условию
\begin{equation*}
\begin{split}
&(\int_{(\R_+^d)^J} (t^J)^{-\e^J -\theta \alpha^J}
(\Omega^{l \chi_J}(f, t^J)_{L_p(D)})^\theta dt^J)^{1/\theta} =\\
&(\int_{(\R_+^d)^J} (\prod_{j \in J} t_j^{-1 -\theta \alpha_j})
(\Omega^{l \chi_J}(f, t^{\s(l \chi_J)})_{L_p(D)})^\theta \prod_{j
\in J} dt_j)^{1/\theta} < \infty (\le 1).
\end{split}
\end{equation*}

При $ \theta = \infty $ положим $ S_{p,\infty}^\alpha B(D) =
S_p^\alpha H(D), \mathcal S_{p,\infty}^\alpha \mathcal B(D) =
\mathcal S_p^\alpha \mathcal H(D). $

Как известно (см., например, [8]), имеет место включение
\begin{equation*} \tag{3.1.1}
\mathcal S_{p, \theta}^\alpha \mathcal B(D) \subset c_1(\alpha)
\mathcal S_p^\alpha \mathcal H(D),
\end{equation*}
где $ c_1(\alpha) = \prod_{j=1}^d 2^{1+\alpha_j}. $

В [5], [8] установлена справедливость леммы 3.1.1 и теоремы 3.1.2.

Лемма 3.1.1

Пусть $ d \in \N, l \in \N^d, 1 \le p < \infty, 1 \le q \le
\infty. $ Тогда существуют константы $ c_2(d,l,p,q) >0 $ и $ c_3(d) >0 $
такие, что для любой функции $ f \in L_p(I^d) $ при $ \kappa \in \Z_+^d $
выполняется неравенство
\begin{equation*} \tag{3.1.2}
\begin{split}
&\| \mathcal E_\kappa^{d,l-\e}f \|_{L_q(I^d)} \le c_2 (\prod_{j
\in \s(\kappa)} 2^{\kappa_j (p^{-1} +(p^{-1} -q^{-1})_+)})
\biggl(\int_{ (c_3 2^{-\kappa} B^d)^{\s(\kappa)}}\\
&\int_{ (I^d)_\xi^{l \chi_{\s(\kappa)}}} |\Delta_\xi^{l
\chi_{\s(\kappa)}} f(x)|^p dx d\xi^{\s(\kappa)} \biggr)^{1/p},
\end{split}
\end{equation*}
где $ t_+ = \frac{1}{2} (t +|t|), t \in \R. $

Теорема 3.1.2

Пусть $ d \in \N, \alpha \in \R_+^d, l = l(\alpha), 1 \le p <
\infty, 1 \le q \le \infty. $
Тогда существует константа $ c_4(d,\alpha,p,q) >0 $ такая, что для
любой функции $ f \in \mathcal S_p^\alpha \mathcal H(I^d) $ при $
\kappa \in \Z_+^d \setminus \{0\} $ соблюдается неравенство
\begin{equation*} \tag{3.1.3}
\| \mathcal E_\kappa^{d,l -\e}f \|_{L_q(I^d)} \le c_4 2^{-(\kappa,
\alpha -(p^{-1} -q^{-1})_+ \e)}
\end{equation*}
и при выполнении условия
\begin{equation*} \tag{3.1.4}
\alpha -(p^{-1} -q^{-1})_+ \e >0
\end{equation*}
в $ L_q(I^d) $ имеет место равенство (1.5.3).
\bigskip

3.2. В этом пункте напомним некоторые сведения о поперечниках (см. [12]).

Пусть $ C $ -- подмножество банахова пространства $X$ и
$n \in \Z_+$. Тогда $n$-поперечни\-ком по Колмогорову множества $C$
в  пространстве  $X$ называется величина
$$
d_n(C,X)=\inf_{M \in \mathcal M_n(X)}\sup_{x \in C}\inf_{y \in
M}\|x-y\|_X,
$$
где $ \mathcal M_n(X) $ -- совокупность всех плоскостей  $ M $ в $
X, $ у которых  $ \dim  M   \le   n. $

Отметим, что для  симметричного (относительно нуля)  множества $ C
$  значение величины $ d_n(C,X) $ не изменится, если вместо
совокупности  всех  плоскостей размерности  не  больше  $ n $ в
определении  $  d_n(C,X)  $ рассматривать лишь совокупность   тех
из   них,   которые  проходят  через $ 0, $ т.е.  линейных
подпространств.

Предложение 3.2.1

Пусть  $  U:  X  \mapsto  Y$ -- непрерывное линейное  отображение
банахова пространства $ X $ в банахово пространство $Y $ и $ C
\subset X$ -- некоторое множество. Тогда  при $ n \in \Z_+ $ справедлива оценка
\begin{equation*} \tag{3.2.1}
d_n(U(C),Y) \le \|U\|_{\mathcal B(X,Y)} d_n(C,X).
\end{equation*}

Напомним, что при  $ 1\le  p,q \le \infty, n \in \N $ для $ x  \in \R^n  $ справедливо неравенство
\begin{equation*} \tag{3.2.2}
\|x\|_{l_q^n} \le n^{(1/q -1/p)_+ } \|x\|_{l_p^n}.
\end{equation*}

Как показано в [13], существует константа $ c_1 >0 $  такая, что для $ n,m \in \N: n \le m, $
имеет место оценка
\begin{equation*} \tag{3.2.3}
d_n(B(l_2^m), l_{\infty}^m) \le c_1 n^{-1/2}
(1+\log(m/n))^{3/2}.
\end{equation*}
\bigskip

3.3. В этом пункте с использованием следствия из теоремы 2.2.4 проводится
оценка сверху колмогоровского $n$-поперечника класса
$ \mathcal S_{p,\theta}^\alpha \mathcal B(I^d) $  в пространстве
$ L_q(I^d). $ Отметим, что вывод этой оценки имеет как сходство, так и
отличия от вывода подобной оценки для периодических функций в [2].

Лемма 3.3.1

Пусть $ d \in \N, \alpha \in \R_+^d, l = l(\alpha), 1 \le \theta
\le \infty, 1 \le p < \infty, 1 \le q < \infty $ и выполнено
условие (3.1.4). Пусть еще $ J = \{j =1,\ldots,d: \alpha_j =
\mn(\alpha)\}, \cmn = \cmn(\alpha) = \card J, $ а вектор $ \beta
\in \R_+^d $ удовлетворяет условиям
\begin{multline*} \tag{3.3.1}
\beta_j =1, \text{для} j \in J; \beta_j > 1, \beta_j^{-1}
(\alpha_j -(p^{-1} -q^{-1})_+) > \mn(\alpha -(p^{-1} -q^{-1})_+
\e)\\
\text{ для } j \in J^\prime = \{1,\ldots,d\} \setminus J.
\end{multline*}
Тогда существует константа $ c_1(d,\alpha,\theta,p,q,\beta) >0$ такая,
что для $ f \in \mathcal S_{p,\theta}^\alpha \mathcal B(I^d) $ при любом
$ r \in \N $ имеет место неравенство
\begin{equation*} \tag{3.3.2}
\| \sum_{\kappa \in \Z_+^d: (\kappa, \beta) > r }
\mathcal E_\kappa^{d,l-\e} f \|_{L_q(I^d)} \le
c_1 \begin{cases} 2^{-\mn(\alpha -(1/p -1/q) \e) r}
r^{(\cmn(\alpha) -1)(1/\q -1/\theta)_+}, \text{ при } p \le q; \\
2^{-\mn(\alpha) r} r^{(\cmn(\alpha) -1)(1/\p -1/\theta)_+}, \text{
при } q < p, \end{cases}
\end{equation*}
где $ \q = \min(2,q),  \p = \min(2,p). $

Доказательство.

Сначала установим справедливость (3.3.2) при $ p \le q. $
Принимая во внимание (3.1.1), (3.1.3), (3.1.4),
в соответствии с (2.2.44), (1.4.6)
для $ f \in \mathcal S_{p,\theta}^\alpha \mathcal B(I^d) $ при $ r \in \N $
имеем неравенство
\begin{multline*} \tag{3.3.3}
\| \sum_{\kappa \in \Z_+^d: (\kappa, \beta) > r }
\mathcal E_\kappa^{d,l-\e} f \|_{L_q(I^d)} \le \\
c_2 ( \sum_{\kappa^\prime \in \Z_+^d} \| \mathcal E_{\kappa^\prime}^{d,l-\e}
( \sum_{\kappa \in \Z_+^d: (\kappa, \beta) > r }
\mathcal E_\kappa^{d,l-\e} f) \|_{L_q(I^d)}^{\q} )^{1/\q} = \\
c_2 ( \sum_{\kappa^\prime \in \Z_+^d} \|
 \sum_{\kappa \in \Z_+^d: (\kappa, \beta) > r }
\mathcal E_{\kappa^\prime}^{d,l-\e}
\mathcal E_\kappa^{d,l-\e} f \|_{L_q(I^d)}^{\q} )^{1/\q} = \\
c_2 (\sum_{\kappa \in \Z_+^d: (\kappa, \beta) > r }
\| \mathcal E_\kappa^{d,l-\e} f \|_{L_q(I^d)}^{\q} )^{1/\q}.
\end{multline*}

Далее, для $ f \in \mathcal S_{p,\theta}^\alpha \mathcal B(I^d) $
при $ \kappa \in \Z_+^d, $ используя (3.1.2), выводим
\begin{equation*} \tag{3.3.4}
\begin{split}
\| \mathcal E_\kappa^{d,l-\e}f \|_{L_q(I^d)} &\le c_3 (\prod_{j
\in \s(\kappa)} 2^{\kappa_j (p^{-1} +(p^{-1} -q^{-1})_+)})
(\int_{(c_4 2^{-\kappa} B^d)^{\s(\kappa)}}
\int_{ (I^d)_\xi^{l \chi_{\s(\kappa)}}}\\
|\Delta_\xi^{l \chi_{\s(\kappa)}} f(x)|^p dx d\xi^{\s(\kappa)}
)^{1/p} &\le c_5 2^{(\kappa, (1/p -1/q)_+ \e)} \Omega^{l
\chi_{\s(\kappa)}}(f, c_4 (2^{-\kappa})^{\s(\kappa)})_{L_p(I^d)}.
\end{split}
\end{equation*}

Учитывая (3.3.4), для $ f \in \mathcal S_{p,\theta}^\alpha
\mathcal B(I^d) $ получаем
\begin{multline*} \tag{3.3.5}
( \sum_{\kappa \in \Z_+^d: (\kappa, \beta) > r} \!\| \mathcal
E_\kappa^{d,l-\e} f \|_{L_q(I^d)}^{\q} )^{1/\q} \\
\le (\sum_{\kappa \in \Z_+^d: (\kappa, \beta) > r} \!(c_5
2^{(\kappa, (1/p -1/q)_+ \e)} \Omega^{l \chi_{\s(\kappa)}} (f, c_4
(2^{-\kappa})^{\s(\kappa)})_{L_p(I^d)})^{\q} )^{1/\q} \\
= c_5 (\sum_{\kappa \in \Z_+^d: (\kappa, \beta) > r} \!2^{-\q
(\kappa, \alpha -(1/p -1/q)_+ \e)} (2^{(\kappa, \alpha)}
\Omega^{l \chi_{\s(\kappa)}}(f, c_4
(2^{-\kappa})^{\s(\kappa)})_{L_p(I^d)})^{\q} )^{1/\q}.
\end{multline*}

При оценке правой части (3.3.5) рассмотрим два случая. В первом
случае, когда $ \theta > \q, $ применяя неравенство Гельдера с
показателем $ \theta /\q > 1,$ для $ f \in \mathcal
S_{p,\theta}^\alpha \mathcal B(I^d) $ находим, что
\begin{multline*} \tag{3.3.6}
( \sum_{\kappa \in \Z_+^d: (\kappa, \beta) > r} \!2^{-\q(\kappa,
\alpha -(1/p -1/q)_+ \e)} (2^{(\kappa, \alpha)} \Omega^{l
\chi_{\s(\kappa)}}(f,
c_4 (2^{-\kappa})^{\s(\kappa)})_{L_p(I^d)})^{\q} )^{1/\q} \\
\le ( \sum_{\kappa \in \Z_+^d: (\kappa, \beta) > r}
\!2^{-\q \theta (\kappa, \alpha -(1/p -1/q)_+ \e) /(\theta -\q)})^{1/\q -1/\theta} \\
\times ( \sum_{\kappa \in \Z_+^d: (\kappa, \beta) > r}
\!(2^{(\kappa, \alpha)} \Omega^{l \chi_{\s(\kappa)}}(f, c_4
(2^{-\kappa})^{\s(\kappa)})_{L_p(I^d)})^\theta )^{1/\theta}.
\end{multline*}

Согласно (1.2.3), с учетом (3.3.1) получаем, что
\begin{multline*} \tag{3.3.7}
( \sum_{\kappa \in \Z_+^d: (\kappa, \beta) > r}
2^{-\q \theta (\kappa, \alpha -(1/p -1/q)_+ \e) /(\theta -\q)})^{1/\q -1/\theta} \\
\le (c_6 2^{-\! \mn(\q \theta (\theta \!-\! \q)^{\!-\! 1}
\beta^{\!-\! 1} (\alpha - (1/p - 1/q)_+ \!\e)) r} r^{\cmn(\q
\theta (\theta \!-\! \q)^{\!-1} \beta^{\!-1}
(\alpha - (1/p - 1/q)_+ \!\e)) - 1})^{1/\q - 1/\theta} \\
= (c_6  2^{-\q \theta (\theta -\q)^{-1} \mn(\alpha -(1/p -1/q)_+
\e) r}
r^{\cmn(\alpha) -1})^{1/\q -1/\theta} \\
= c_7 2^{-\mn(\alpha -(1/p -1/q)_+ \e) r} r^{(\cmn(\alpha)
-1)(1/\q -1/\theta)}.
\end{multline*}

Кроме того, для $ f \in \mathcal S_{p,\theta}^\alpha \mathcal B(I^d) $
выводим
\begin{multline*} \tag{3.3.8}
( \sum_{\kappa \in \Z_+^d: (\kappa, \beta) > r} (2^{(\kappa,
\alpha)} \Omega^{l \chi_{\s(\kappa)}}(f,
c_4 (2^{-\kappa})^{\s(\kappa)})_{L_p(I^d)})^\theta )^{1/\theta} \\
\le ( \sum_{\kappa \in \Z_+^d \setminus \{0\}} (2^{(\kappa,
\alpha)} \Omega^{l \chi_{\s(\kappa)}}(f,
c_4 (2^{-\kappa})^{\s(\kappa)})_{L_p(I^d)})^\theta )^{1/\theta} \\
= ( \sum_{ \J \subset \Nu_{1,d}^1: \J \ne \emptyset} \sum_{\kappa
\in \Z_+^d: \s(\kappa) = \J } (2^{(\kappa, \alpha)} \Omega^{l
\chi_{\J}}(f,
c_4 (2^{-\kappa})^{\J})_{L_p(I^d)})^\theta)^{1/\theta} \\
= ( \sum_{ \J \subset \Nu_{1,d}^1: \J \ne \emptyset} \sum_{\kappa
\in \Z_+^d: \s(\kappa) = \J } (2^{(\kappa^{\J}, \alpha^{\J})} \Omega^{l
\chi_{\J}}(f, c_4 2^{-\kappa^{\J}})_{L_p(I^d)})^\theta)^{1/\theta} \\
=(\!\sum_{ \J \subset \Nu_{1,d}^1: \J \ne \emptyset} \sum_{\kappa^{\J}
\in (\N^d)^{\J}} \!\int_{2^{-\kappa^{\J}} \!+\! 2^{-\kappa^{\J}} (I^d)^{\J}}
\!2^{(\kappa^{\J} \!, \e^{\J}) +\theta (\kappa^{\J} \!, \alpha^{\J})}
\!(\Omega^{l \chi_{\J}}(f, c_4 2^{-\kappa^{\J}})_{L_p(I^d)})^\theta
dt^{\J})^{1/\theta} \\
\le (\!\sum_{ \J \subset \Nu_{1,d}^1: \J \ne \emptyset}
\sum_{\kappa^{\J} \in (\N^d)^{\J}} \!\int_{2^{-\kappa^{\J}} \!+\!
2^{-\kappa^{\J}} (I^d)^{\J}} \! c_8 (t^{\J})^{-\e^{\J} -\theta \alpha^{\J}}
(\Omega^{l \chi_{\J}}(f, c_4 t^{\J})_{L_p(I^d)})^\theta dt^{\J})^{1/\theta} \\
\le ( \sum_{ \J \subset \Nu_{1,d}^1: \J \ne \emptyset} \int_{
(I^d)^{\J}} c_8 (t^{\J})^{-\e^{\J} -\theta \alpha^{\J}}
(\Omega^{l \chi_{\J}}(f, c_4 t^{\J})_{L_p(I^d)})^\theta dt^{\J})^{1/\theta} \\
\le (c_8 \sum_{ \J \subset \Nu_{1,d}^1: \J \ne \emptyset} \int_{
(\R_+^d)^{\J}} (t^{\J})^{-\e^{\J} -\theta \alpha^{\J}} (\Omega^{l \chi_{\J}}(f,
c_4 t^{\J})_{L_p(I^d)})^\theta dt^{\J})^{1/\theta} \le c_9.
\end{multline*}

Соединяя (3.3.5) -- (3.3.8), для $ f \in \mathcal S_{p,\theta}^\alpha
\mathcal B(I^d) $ при $ \theta > \q $ приходим к неравенству
\begin{equation*} \tag{3.3.9}
( \sum_{\kappa \in \Z_+^d: (\kappa, \beta) > r} \!\| \mathcal
E_\kappa^{d,l-\e} f \|_{L_q(I^d)}^{\q} )^{\frac{1}{\q}} \!\le c_{10}
2^{-\!\mn(\alpha -(1/p -1/q)_+ \e) r} r^{(\cmn(\alpha) \!-\!
1)(1/\q -1/\theta)_+}.
\end{equation*}

Оценивая правую часть (3.3.5) при $ \theta \le \q, $ для $ f \in
\mathcal S_{p,\theta}^\alpha \mathcal B(I^d) $ в силу (3.3.1) и
неравенства Гельдера с показателем $ \theta /\q \le 1, $ имеем
\begin{multline*} \tag{3.3.10}
( \sum_{\kappa \in \Z_+^d: (\kappa, \beta) > r} 2^{-\q(\kappa,
\alpha -(1/p -1/q)_+ \e)} (2^{(\kappa, \alpha)} \Omega^{l
\chi_{\s(\kappa)}}(f, c_4
(2^{-\kappa})^{\s(\kappa)})_{L_p(I^d)})^{\q} )^{1/\q} \\
=( \sum_{\kappa \in \Z_+^d: (\kappa, \beta) > r} 2^{-\q(\kappa,
\beta \beta^{-1} (\alpha -(1/p -1/q)_+ \e))} (2^{(\kappa, \alpha)}
\Omega^{l \chi_{\s(\kappa)}}(f, c_4
(2^{-\kappa})^{\s(\kappa)})_{L_p(I^d)})^{\q} )^{1/\q} \\
\le ( \sum_{\kappa \in \Z_+^d: (\kappa, \beta) > r} 2^{-\q
\mn(\beta^{-1} (\alpha -(1/p -1/q)_+ \e))(\kappa, \beta)}
(2^{(\kappa, \alpha)} \Omega^{l \chi_{\s(\kappa)}}(f, c_4
(2^{-\kappa})^{\s(\kappa)})_{L_p(I^d)})^{\q} )^{1/\q} \\
\le ( \sum_{\kappa \in \Z_+^d: (\kappa, \beta) > r} 2^{-\q
\mn(\alpha -(1/p -1/q)_+ \e) r} (2^{(\kappa, \alpha)} \Omega^{l
\chi_{\s(\kappa)}}(f, c_4
(2^{-\kappa})^{\s(\kappa)})_{L_p(I^d)})^{\q} )^{1/\q} \\
= 2^{-\mn(\alpha -(1/p -1/q)_+ \e) r} ( \sum_{\kappa \in \Z_+^d:
(\kappa, \beta) > r} (2^{(\kappa, \alpha)} \Omega^{l
\chi_{\s(\kappa)}}(f, c_4
(2^{-\kappa})^{\s(\kappa)})_{L_p(I^d)})^{\q} )^{1/\q} \\
\le 2^{-\mn(\alpha -(1/p -1/q)_+ \e) r} ( \sum_{\kappa \in \Z_+^d:
(\kappa, \beta) > r} (2^{(\kappa, \alpha)} \Omega^{l
\chi_{\s(\kappa)}}(f, c_4
(2^{-\kappa})^{\s(\kappa)})_{L_p(I^d)})^\theta )^{1/ \theta}.
\end{multline*}

Подставляя (3.3.8) в (3.3.10) и соединяя с (3.3.5), для $ f \in
\mathcal S_{p,\theta}^\alpha \mathcal B(I^d) $ получаем
(3.3.9) при $ \theta \le \q. $

Объединяя (3.3.3) с (3.3.9), приходим к (3.3.2) при $ p \le q. $

Для получения (3.3.2) в случае $ q < p $ достаточно заметить, что в этом
случае справедлива оценка
\begin{equation*}
\| \sum_{\kappa \in \Z_+^d: (\kappa, \beta) > r } \mathcal
E_\kappa^{d,l-\e} f \|_{L_q(I^d)} \le \| \sum_{\kappa \in \Z_+^d:
(\kappa, \beta) > r } \mathcal E_\kappa^{d,l-\e} f \|_{L_p(I^d)},
f \in \mathcal S_{p,\theta}^\alpha \mathcal B(I^d),
\end{equation*}
и применить (3.3.2) при $ q = p. \square$

Лемма 3.3.2

Пусть выполнены условия леммы 3.3.1 и $ p \le q. $ Тогда для $ C =
\mathcal S_{p,\theta}^\alpha \mathcal B(I^d) $ и $ X = L_q(I^d) $
существует константа $ c_{11}(d,\alpha,\theta,p,q,\beta) >0$
такая, что для любых $ n,r \in \N, j_0 \in \Z_+ $ и любого набора
чисел $ \{n_\kappa \in \Z_+: \kappa \in \Z_+^d, r < (\kappa,
\beta) \le r +j_0\}, $ подчиненных условию
\begin{equation*} \tag{3.3.11}
n  \ge \sum_{\kappa \in \Z_+^d: (\kappa, \beta) \le r}
\mathfrak R_\kappa^{d,l-\e} +\sum_{\kappa \in \Z_+^d:
r < (\kappa, \beta) \le r +j_0} n_\kappa,
\end{equation*}
можно построить линейное подпространство $ M \subset X, $ размерность
которого $ \dim M \le n, $ а для $ f \in C $ имеет место неравенство
\begin{multline*} \tag{3.3.12}
\inf_{ g \in M} \| f -g \|_{L_q(I^d)} \le  c_{11} \biggl(
\biggl(\sum_{\kappa \in \Z_+^d: r < (\kappa, \beta) \le r+j_0,
n_\kappa < \mathfrak R_\kappa^{d,l-\e}}
(2^{(\kappa, (p^{-1}  -q^{-1}) \e)} \\
\times d_{n_\kappa} (B(l_p^{\mathfrak R_{\kappa}^{d,l-\e}}),
l_q^{\mathfrak R_{\kappa}^{d,l-\e}} ) \Omega^{l
\chi_{\s(\kappa)}}(f,
c_4 (2^{-\kappa})^{\s(\kappa)})_{L_p(I^d)})^{\q} \biggr)^{1/\q} \\
+2^{-(r+j_0) \mn(\alpha -(p^{-1}
-q^{-1})_+ \e)} (r +j_0)^{(\cmn -1)(1/\q -1/\theta)_+}\biggr).
\end{multline*}

Доказательство.

Пусть $ n,r \in \N, j_0 \in \Z_+ $ и набор чисел
$ \{n_\kappa \in \Z_+: \kappa \in \Z_+^d, r < (\kappa, \beta)
\le r +j_0\} $ удовлетворяют условию (3.3.11).
Для каждого $ \kappa \in \Z_+^d: r < (\kappa, \beta)
\le r +j_0, n_\kappa < \mathfrak R_\kappa^{d,l-\e} $ фиксируем линейное
подпространство $ M_\kappa \subset \mathfrak P_\kappa^{d,l-\e}, $
для которого $ \dim M_\kappa \le n_\kappa $ и
\begin{equation*}
\sup_{f \in B(L_p(I^d)) \cap \mathfrak P_\kappa^{d,l-\e}}
\inf_{g \in M_\kappa} \|f -g\|_{L_q(I^d)} <
2d_{n_\kappa}(B(L_p(I^d)) \cap \mathfrak P_\kappa^{d,l-\e},
\mathfrak P_\kappa^{d,l-\e} \cap L_q(I^d)) >0,
\end{equation*}
а, следовательно, для каждого $ f \in B(L_p(I^d)) \cap
\mathfrak P_\kappa^{d,l-\e} $ выполняется неравенство
\begin{equation*}
\inf_{g \in M_\kappa} \|f -g\|_{L_q(I^d)} <
2d_{n_\kappa}(B(L_p(I^d)) \cap \mathfrak P_\kappa^{d,l-\e},
\mathfrak P_\kappa^{d,l-\e} \cap L_q(I^d)).
\end{equation*}
А для $ \kappa \in \Z_+^d: r < (\kappa, \beta) \le r +j_0, n_\kappa \ge
\mathfrak R_\kappa^{d,l-\e}, $ положим $ M_\kappa =
\mathfrak P_\kappa^{d,l-\e}. $

Учитывая сказанное, для $ f \in C $ при $ \kappa \in \Z_+^d: r < (\kappa, \beta)
\le r +j_0, n_\kappa < \mathfrak R_\kappa^{d,l-\e}, $ выберем функцию
$ g_\kappa = g_\kappa(f) \in M_\kappa $ такую, что
\begin{equation*}
\begin{split}
\| (1 / \| \mathcal E_\kappa^{d,l-\e} f \|_{L_p(I^d)}) \mathcal
E_\kappa^{d,l-\e} f -g_\kappa \|_{L_q(I^d)}\\ <
2d_{n_\kappa}(B(L_p(I^d)) \cap \mathfrak P_\kappa^{d,l-\e},
\mathfrak P_\kappa^{d,l-\e} \cap L_q(I^d)),
\text{ если } \mathcal E_\kappa^{d,l-\e} f \ne 0; \\
g_\kappa =0, \text{ если } \mathcal E_\kappa^{d,l-\e} f =0,
\end{split}
\end{equation*}
и, значит,
\begin{multline*} \tag{3.3.13}
\| \mathcal E_\kappa^{d,l-\e} f - \| \mathcal E_\kappa^{d,l-\e} f
\|_{L_p(I^d)} g_\kappa \|_{L_q(I^d)} \\
\le \| \mathcal
E_\kappa^{d,l-\e} f \|_{L_p(I^d)} 2d_{n_\kappa}(B(L_p(I^d)) \cap
\mathfrak P_\kappa^{d,l-\e}, \mathfrak P_\kappa^{d,l-\e} \cap
L_q(I^d)).
\end{multline*}
Тогда получаем, что для $ f \in C $ и подпространства
$$
M = \sum_{\kappa \in \Z_+^d: (\kappa, \beta) \le r} \mathfrak P_\kappa^{d,l-\e}
+\sum_{\kappa \in \Z_+^d: r < (\kappa, \beta) \le r +j_0 } M_\kappa
$$
ввиду (3.1.1), (3.1.4), (1.5.3) соблюдается неравенство
\begin{multline*} \tag{3.3.14}
\inf_{ g \in M} \| f -g \|_{L_q(I^d)} \le \biggl\| f -\sum_{\kappa
\in \Z_+^d: (\kappa, \beta) \le r} \mathcal E_\kappa^{d,l-\e} f
-\sum_{\kappa \in \Z_+^d: r < (\kappa, \beta) \le r +j_0, n_\kappa
\ge \mathfrak R_\kappa^{d,l-\e}} \mathcal E_\kappa^{d,l-\e} f \\
-\sum_{\kappa \in \Z_+^d: r < (\kappa, \beta) \le r +j_0, n_\kappa
< \mathfrak R_\kappa^{d,l-\e}} \| \mathcal E_\kappa^{d,l-\e} f
\|_{L_p(I^d)} g_\kappa \biggr\|_{L_q(I^d)} \\
= \biggl\| \sum_{\kappa \in \Z_+^d: r < (\kappa, \beta) \le r
+j_0, n_\kappa < \mathfrak R_\kappa^{d,l-\e}} (\mathcal
E_\kappa^{d,l-\e} f - \| \mathcal E_\kappa^{d,l-\e} f
\|_{L_p(I^d)} g_\kappa) +\sum_{\kappa \in \Z_+^d: (\kappa, \beta)
> r +j_0} \mathcal
E_\kappa^{d,l-\e} f \biggr\|_{L_q(I^d)} \\
\le \biggl\| \sum_{\kappa \in \Z_+^d: r < (\kappa, \beta) \le r
+j_0, n_\kappa < \mathfrak R_\kappa^{d,l-\e}} (\mathcal
E_\kappa^{d,l-\e} f - \| \mathcal E_\kappa^{d,l-\e} f
\|_{L_p(I^d)} g_\kappa) \biggr\|_{L_q(I^d)} \\
+\| \sum_{\kappa \in \Z_+^d: (\kappa, \beta) > r +j_0} \mathcal
E_\kappa^{d,l-\e} f \|_{L_q(I^d)}.
\end{multline*}

Оценивая первое слагаемое в правой части (3.3.14), на основании (2.2.44),
(1.4.7), (1.4.6), (3.3.13) заключаем, что для $ f \in C $
соблюдается неравенство
\begin{multline*} \tag{3.3.15}
\| \sum_{\kappa \in \Z_+^d: r < (\kappa, \beta) \le r +j_0,
n_\kappa < \mathfrak R_\kappa^{d,l-\e}} (\mathcal
E_\kappa^{d,l-\e} f - \| \mathcal E_\kappa^{d,l-\e} f
\|_{L_p(I^d)} g_\kappa) \|_{L_q(I^d)}\\
 \le c_{12} (
\sum_{\kappa^\prime \in \Z_+^d} \| \mathcal
E_{\kappa^\prime}^{d,l-\e} ( \sum_{\kappa \in \Z_+^d: r < (\kappa,
\beta) \le r +j_0, n_\kappa < \mathfrak R_\kappa^{d,l-\e}}
(\mathcal E_\kappa^{d,l-\e} f - \| \mathcal E_\kappa^{d,l-\e} f
\|_{L_p(I^d)} g_\kappa))
\|_{L_q(I^d)}^{\q} )^{1/\q} \\
=c_{12} ( \sum_{\kappa^\prime \in \Z_+^d} \| \sum_{\kappa \in
\Z_+^d: r < (\kappa, \beta) \le r +j_0, n_\kappa < \mathfrak
R_\kappa^{d,l-\e}} \mathcal E_{\kappa^\prime}^{d,l-\e} (\mathcal
E_\kappa^{d,l-\e} f - \| \mathcal E_\kappa^{d,l-\e} f
\|_{L_p(I^d)} g_\kappa)
\|_{L_q(I^d)}^{\q} )^{1/\q} \\
=c_{12} ( \sum_{\kappa \in \Z_+^d: r < (\kappa, \beta) \le r +j_0,
n_\kappa < \mathfrak R_\kappa^{d,l-\e}} \| \mathcal
E_\kappa^{d,l-\e} f - \| \mathcal E_\kappa^{d,l-\e} f
\|_{L_p(I^d)}
g_\kappa \|_{L_q(I^d)}^{\q} )^{1/\q} \\
\le c_{12} ( \sum_{\kappa \in \Z_+^d: r < (\kappa, \beta) \le r
+j_0, n_\kappa < \mathfrak R_\kappa^{d,l-\e}} \| \mathcal
E_\kappa^{d,l-\e} f \|_{L_p(I^d)}^{\q}\\ \times
(2d_{n_\kappa}(B(L_p(I^d)) \cap \mathfrak P_\kappa^{d,l-\e},
\mathfrak P_\kappa^{d,l-\e} \cap L_q(I^d)) )^{\q})^{1/\q}  \\
\le c_{13} ( \sum_{\kappa \in \Z_+^d: r < (\kappa, \beta) \le r
+j_0, n_\kappa < \mathfrak R_\kappa^{d,l-\e}} (\| \mathcal
E_\kappa^{d,l-\e} f \|_{L_p(I^d)}\\ \times
d_{n_\kappa}(B(L_p(I^d)) \cap \mathfrak P_\kappa^{d,l-\e},
\mathfrak P_\kappa^{d,l-\e} \cap L_q(I^d)) )^{\q})^{1/\q}.
\end{multline*}

Заметим, что при $ \kappa \in \Z_+^d: r < (\kappa, \beta) \le r +j_0,
n_\kappa < \mathfrak R_\kappa^{d,l-\e}, $
в силу (3.2.1), (1.4.11) справедлива оценка
\begin{multline*} \tag{3.3.16}
d_{n_\kappa}(B(L_p(I^d)) \cap \mathfrak P_\kappa^{d,l-\e},
\mathfrak P_\kappa^{d,l-\e} \cap L_q(I^d)) \\
=d_{n_\kappa}((\mathfrak I_\kappa^{d,l-\e})^{-1} \mathfrak
I_\kappa^{d,l-\e}(B(L_p(I^d)) \cap \mathfrak P_\kappa^{d,l-\e}),
\mathfrak P_\kappa^{d,l-\e} \cap L_q(I^d)) \\
\le c_{14} 2^{-(\kappa, \e)/q} d_{n_\kappa}( \mathfrak
I_\kappa^{d,l-\e}(B(L_p(I^d)) \cap \mathfrak P_\kappa^{d,l-\e}),
l_q^{\mathfrak R_\kappa^{d,l-\e}}) \\
\le c_{14} 2^{-(\kappa,
\e)/q} d_{n_\kappa}(c_{15} 2^{(\kappa, \e)/p} B(l_p^{\mathfrak
R_\kappa^{d,l-\e}}),
l_q^{\mathfrak R_\kappa^{d,l-\e}}) \\
=c_{16} 2^{(\kappa, (1/p -1/q) \e)} d_{n_\kappa}(B(l_p^{\mathfrak
R_\kappa^{d,l-\e}}), l_q^{\mathfrak R_\kappa^{d,l-\e}}).
\end{multline*}

Поэтому, подставляя (3.3.16) и (3.3.4) при $ q=p $ в (3.3.15), приходим к
неравенству
\begin{multline*} \tag{3.3.17}
\| \sum_{\kappa \in \Z_+^d: r < (\kappa, \beta) \le r +j_0,
n_\kappa < \mathfrak R_\kappa^{d,l-\e}} (\mathcal
E_\kappa^{d,l-\e} f - \| \mathcal E_\kappa^{d,l-\e} f
\|_{L_p(I^d)} g_\kappa) \|_{L_q(I^d)}\\
 \le c_{17} (\sum_{\kappa \in \Z_+^d:
r < (\kappa, \beta) \le r+j_0,
n_\kappa < \mathfrak R_\kappa^{d,l-\e}}
(2^{(\kappa, (p^{-1} -q^{-1}) \e)} \\
\times d_{n_\kappa} (B(l_p^{\mathfrak R_{\kappa}^{d,l-\e}}),
l_q^{\mathfrak R_{\kappa}^{d,l-\e}} ) \Omega^{l
\chi_{\s(\kappa)}}(f, c_4
(2^{-\kappa})^{\s(\kappa)})_{L_p(I^d)})^{\q} )^{1/\q}, f \in C.
\end{multline*}

Применение (3.3.2) ко второму слагаемому в правой части (3.3.14)
дает оценку
\begin{multline*} \tag{3.3.18}
\| \sum_{\kappa \in \Z_+^d: (\kappa, \beta) > r +j_0} \mathcal
E_\kappa^{d,l-\e} f \|_{L_q(I^d)}\\ \le c_{18} 2^{-(r+j_0)
\mn(\alpha -(p^{-1} -q^{-1})_+ \e)} (r +j_0)^{(\cmn -1)(1/\q
-1/\theta)_+}, f \in C.
\end{multline*}

Соединяя (3.3.14), (3.3.17), (3.3.18), получаем (3.3.12).

Причем, по построению, ввиду (3.3.11) соблюдается неравенство
\begin{multline*}
\dim M \le \sum_{\kappa \in \Z_+^d: (\kappa, \beta) \le r}
\dim \mathfrak P_\kappa^{d,l-\e}
+\sum_{\kappa \in \Z_+^d: r < (\kappa, \beta) \le r +j_0 } \dim M_\kappa = \\
\sum_{\kappa \in \Z_+^d: (\kappa, \beta) \le r}
\dim \mathfrak P_\kappa^{d,l-\e}
+\sum_{\kappa \in \Z_+^d: r < (\kappa, \beta) \le r +j_0, n_\kappa \ge
\mathfrak R_\kappa^{d,l-\e}}
\dim \mathfrak P_\kappa^{d,l-\e} \\
 +\sum_{\kappa \in \Z_+^d: r < (\kappa, \beta) \le r +j_0, n_\kappa <
\mathfrak R_\kappa^{d,l-\e}} \dim M_\kappa =
\sum_{\kappa \in \Z_+^d: (\kappa, \beta) \le r}
\mathfrak R_\kappa^{d,l-\e} \\
+\sum_{\kappa \in \Z_+^d: r < (\kappa, \beta) \le r +j_0, n_\kappa \ge
\mathfrak R_\kappa^{d,l-\e}}
\mathfrak R_\kappa^{d,l-\e}
+\sum_{\kappa \in \Z_+^d: r < (\kappa, \beta) \le r +j_0, n_\kappa <
\mathfrak R_\kappa^{d,l-\e}} \dim M_\kappa \le \\
\sum_{\kappa \in \Z_+^d: (\kappa, \beta) \le r}
\mathfrak R_\kappa^{d,l-\e}
+\sum_{\kappa \in \Z_+^d: r < (\kappa, \beta) \le r +j_0} n_\kappa \le n. \square
\end{multline*}

Теорема 3.3.3

Пусть $ d \in \N, \alpha \in \R_+^d, 1 \le p < \infty, 1 \le q <
\infty, 1 \le \theta \le \infty. $ Пусть еще $ J = \{j \in \{1,
\ldots, d\}: \alpha_j = \mn(\alpha)\}, \cmn = \card J. $ Тогда для
$ C = \mathcal S_{p,\theta}^\alpha \mathcal B(I^d) $ и $ X =
L_q(I^d) $ существуют константы $ c_{19}(C,X) >0$ и $ c_{20}(C,X)
>0 $ такие, что при достаточно больших $ n \in \N $ имеют место
неравенства
\begin{equation*} \tag{3.3.19}
d_n(C,X) \le c_{19}
\begin{cases}
n^{-\mn(\alpha -(p^{-1} -q^{-1})_+ \e)}
(\log n)^{(\mn(\alpha -(p^{-1} -q^{-1})_+ \e) +(1/\pq -1/\theta)_+)
(\cmn -1)}, \\
\text{ при }  q \le p \text{ или } (p < q \le 2 \text{ и соблюдении условия (3.1.4)}); \\
n^{-\mn(\alpha -(p^{-1} -1/2)_+ \e)} (\log n)^{(\mn(\alpha
-(p^{-1} -1/2)_+ \e) +(1/2 -1/\theta)_+)
(\cmn -1)}, \\
\text{ при } q \ge \max(2, p) \text{ и выполнении условия }\\
\alpha -p^{-1} \e -(1/2 -p^{-1})_+ \e
>0,
\end{cases}
\end{equation*}

где $ \pq = \min(2, \max(p, q)), $
\begin{multline*} \tag{3.3.19'}
d_n(C,X) \ge c_{20}
n^{-\mn(\alpha -(p^{-1} -1/2)_+ \e)} (\log n)^{(\mn(\alpha
-(p^{-1} -1/2)_+ \e) +(1/2 -1/\theta)_+)
(\cmn -1)}, \\
\text{ при } q \ge 2 \text{ и выполнении условия (3.1.4)}.
\end{multline*}

Доказательство.

В первом случае, когда $  q \le p \text{или} (p < q \le 2 \text{и
соблюдается условие (3.1.4)}), $ полагая $ l =l(\alpha) $ и
фиксируя вектор $ \beta \in \R_+^d, $ подчиненный условиям
(3.3.1), для достаточно  большого $ n \in \N $ выберем число $ r
\in \N $ так, чтобы удовлетворить неравенству
\begin{equation*} \tag{3.3.20}
\sum_{\kappa \in \Z_+^d: (\kappa, \beta) \le r}
\mathfrak R_\kappa^{d,l-\e} \le n <
\sum_{\kappa \in \Z_+^d: (\kappa, \beta) \le r +1}
\mathfrak R_\kappa^{d,l-\e},
\end{equation*}
и рассмотрим линейное подпространство $ M \subset X, $ определяемое равенством
$$
M = \sum_{\kappa \in \Z_+^d: (\kappa, \beta) \le r}
\mathfrak P_\kappa^{d,l-\e}.
$$

Тогда, учитывая (3.3.20), размерность
$$
\dim M = \sum_{\kappa \in \Z_+^d: (\kappa, \beta) \le r}
\dim \mathfrak P_\kappa^{d,l-\e} =
\sum_{\kappa \in \Z_+^d: (\kappa, \beta) \le r}
\mathfrak R_\kappa^{d,l-\e} \le n,
$$
и для $ f \in C, $ благодаря (1.5.3), (3.3.2), соблюдается неравенство
\begin{multline*}
\inf_{ g \in M} \| f -g \|_{L_q(I^d)}  \le \| f -\sum_{\kappa \in
\Z_+^d: (\kappa, \beta) \le r }
\mathcal E_\kappa^{d,l-\e} f \|_{L_q(I^d)} \\
=\| \sum_{\kappa \in \Z_+^d: (\kappa, \beta) > r } \mathcal
E_\kappa^{d,l-\e} f \|_{L_q(I^d)} \le c_{21} 2^{-r \mn(\alpha
-(p^{-1} -q^{-1})_+ \e)} r^{(\cmn -1)(1/\pq -1/\theta)_+}.
\end{multline*}
Отсюда следует, что
\begin{equation*} \tag{3.3.21}
d_n(C,X) \le \sup_{ f \in C} \inf_{ g \in M} \| f -g \|_X \le c_{21}
2^{-r \mn(\alpha -(p^{-1} -q^{-1})_+ \e)} r^{(\cmn -1)(1/\pq -1/\theta)_+}.
\end{equation*}

Замечая, что вследствие (1.4.10), (1.2.2), (3.3.1) для $ r \in \N $ справедливо неравенство
\begin{equation*} \tag{3.3.22}
c_{22}(d,l,\beta) 2^r r^{\cmn -1} \le  \sum_{\kappa \in \Z_+^d:
(\kappa, \beta) \le r} \mathfrak R_\kappa^{d,l-\e} \le c_{23}(d,l,\beta)
2^r r^{\cmn -1},
\end{equation*}
из (3.3.21), (3.3.20) и (3.3.22) получаем (3.3.19)
в первом случае.

Теперь рассмотрим второй случай, когда $ q \ge \max(2,p) $ и
выполняется условие $ \alpha -p^{-1} \e -(1/2 -p^{-1})_+ \e >0, $
которое влечет (3.1.4).

В этом случае, полагая $ l = l(\alpha) $ и обозначая $ J^\prime =
\{1,\ldots,d\} \setminus J, $ фиксируем вектор $ \beta \in \R_+^d, $
для которого выполняются условия (3.3.1) и для $ j \in J^\prime $
верно неравенство
$$
\beta_j^{-1} (\alpha_j -p^{-1} -(1/2 -p^{-1})_+) >
\mn(\alpha -p^{-1} \e -(1/2 -p^{-1})_+ \e).
$$

Далее, фиксируем $ \epsilon > 0 $ такое,
что для $ j \in J^\prime $
соблюдается неравенство
\begin{multline*}
\beta_j^{-1} (\alpha_j -p^{-1} -(1/2 -p^{-1})_+) >
\mu +\epsilon, \\
\text{ где } \mu = \mn(\alpha -p^{-1} \e -(1/2 -p^{-1})_+ \e),
\end{multline*}
а также выберем $ \gamma >0 $ и $ \gamma^\prime >0 $ так, чтобы выполнялись условия:
\begin{multline*}
\gamma < 1/3, \mu -\gamma /2 >0,\\ (1 +\frac{1}{3\gamma})
\mn(\alpha -(p^{-1} -q^{-1})_+ \e) \ge
\mn(\alpha -(p^{-1} -1/2)_+ \e), \\
\gamma^\prime < \gamma, \epsilon -\gamma^\prime /2 >0.
\end{multline*}

Теперь для достаточно большого $ n \in \N, $ подобрав $ r \in \N $
так, чтобы удовлетворить неравенству
\begin{equation*} \tag{3.3.23}
c_{20} 2^r r^{\cmn -1} \le n < c_{20} 2^{r +1} (r +1)^{\cmn -1},
\end{equation*}
где константа $ c_{20} > 0 $ будет указана ниже,
положим число $ j_0 = j_0(r) =
[\frac{r}{3\gamma}] $ и зададим набор чисел
$ \{ n_\kappa \in \N:
\kappa \in \Z_+^d, r < (\kappa, \beta) \le
r +j_0\} $ равенством
$$
n_\kappa = \min([c_0 2^{r -\gamma j -\gamma^\prime
(\kappa^{J^\prime}, \beta^{J^\prime})}] +1, \mathfrak R_\kappa^{d,l-\e}),
$$
для $ \kappa \in \Z_+^d: r +j -1 < (\kappa, \beta) \le r +j,
j=1, \ldots, j_0, $ где $ c_0 = \max_{ \mathcal J \subset \{1, \ldots, d\}}
\mathfrak R_{\chi_{\mathcal J}}^{d,l-\e} $ (см. (1.4.10)).

Тогда в силу (1.2.4) справедливо неравенство
\begin{multline*} \tag{3.3.24}
\sum_{ \kappa \in \Z_+^d: r < (\kappa, \beta) \le r +j_0} n_\kappa
= \sum_{j=1}^{j_0} \sum_{ \kappa \in \Z_+^d: r +j -1 < (\kappa,
\beta) \le r +j} n_\kappa \\
\le \sum_{j=1}^{j_0} \sum_{ \kappa \in
\Z_+^d: r +j -1 < (\kappa, \beta) \le r +j} ([c_0 2^{r -\gamma
j -\gamma^\prime (\kappa^{J^\prime}, \beta^{J^\prime})}] +1) \\
\le \sum_{j=1}^{j_0} \sum_{ \kappa \in \Z_+^d: r +j -1 < (\kappa,
\beta) \le r +j} 2 c_0 2^{r -\gamma j -\gamma^\prime
(\kappa^{J^\prime}, \beta^{J^\prime})} \\
\le 2 c_0 2^r \sum_{j=1}^{j_0} 2^{-\gamma j} \sum_{ \kappa \in
\Z_+^d: r +j -1 < (\kappa, \beta) \le r +j} 2^{-\gamma^\prime
(\kappa^{J^\prime}, \beta^{J^\prime})} \\
\le c_{27} 2^r
\sum_{j=1}^{j_0} 2^{-\gamma j} (r +j)^{\cmn -1} \le c_{27} 2^r
r^{\cmn -1} \sum_{j=1}^\infty 2^{-\gamma j} (j +1)^{\cmn -1} =
c_{28} 2^r r^{\cmn -1}.
\end{multline*}

Из (3.3.22) и (3.3.24) вытекает, что для $ r, j_0 \in \N, \{ n_\kappa \in \N:
\kappa \in \Z_+^d, r < (\kappa, \beta) \le
r +j_0\} $ соблюдается неравенство
\begin{equation*} \tag{3.3.25}
\sum_{\kappa \in \Z_+^d: (\kappa, \beta) \le r} \mathfrak R_\kappa^{d,l-\e}
+\sum_{\kappa \in \Z_+^d: r < (\kappa, \beta) \le r +j_0} n_\kappa
\le c_{20} 2^r r^{\cmn -1},
\end{equation*}
где $ c_{20} = c_{19} +c_{21}. $
Тогда, сопоставляя (3.3.25), (3.3.23) с (3.3.11), в соответствии с леммой 3.3.2
построим подпространство $ M \subset X, $ размерность которого
$ \dim M \le n, $ а для $ f \in C $ имеет место оценка (3.3.12).
Из (3.3.12), пользуясь тем, что в силу (3.2.1), (1.4.10), (3.2.2) и (3.2.3)
справедливо соотношение
\begin{multline*}
d_{n_\kappa} \bigl(B(l_p^{\mathfrak R_{\kappa}^{d,l-\e}}),
l_q^{\mathfrak R_{\kappa}^{d,l-\e}}\bigr)\\
\le c_{22} 2^{(\kappa, (1/2
-p^{-1})_+ \e +q^{-1} \e)} n_\kappa^{-1/2} (1 +\log
\frac{\mathfrak R_\kappa^{d,l-\e}} {n_\kappa})^{3/2}, \\
\text{ при } \kappa \in \Z_+^d: r < (\kappa, \beta) \le r +j_0,
n_\kappa < \mathfrak R_\kappa^{d,l-\e}, \text{(см. http://)}
\end{multline*}
получаем, что для $ f \in C $ выполняется неравенство
\begin{multline*} \tag{3.3.26}
\inf_{ g \in M} \| f -g \|_X  \le  c_{23} \biggl( (\sum_{\kappa \in \Z_+^d: r <
(\kappa, \beta) \le r +j_0, n_\kappa < \mathfrak R_\kappa^{d,l-\e}}
(2^{-(\kappa, \alpha -p^{-1} \e -(1/2 -p^{-1})_+ \e)}
n_\kappa^{-1/2} \times\\
(1 +\log \frac{\mathfrak R_\kappa^{d,l-\e}}
{n_\kappa})^{3/2}
2^{(\kappa, \alpha)} \Omega^{l \chi_{\s(\kappa)}}(f,
c_4 (2^{-\kappa})^{\s(\kappa)})_{L_p(I^d)})^2 )^{1/2} \\
+2^{-(r+j_0) \mn(\alpha -(p^{-1}
-q^{-1})_+ \e)} (r +j_0)^{(\cmn -1)(1/2 -1/\theta)_+}\biggr).
\end{multline*}

Оценим правую часть (3.3.26).
Для этого заметим, что
при $ j=1, \ldots, j_0 $ для $ \kappa \in \Z_+^d:
r +j -1 < (\kappa, \beta) \le r +j, n_\kappa <
\mathfrak R_\kappa^{d,l-\e}, $
соблюдаются неравенства
\begin{multline*} \tag{3.3.27}
(\kappa, \alpha -p^{-1} \e -(1/2 -p^{-1})_+ \e) = \sum_{i \in J}
\beta_i^{-1} ( \alpha_i -p^{-1}  -(1/2 -p^{-1})_+) \beta_i
\kappa_i \\
+\sum_{i \in J^\prime} \beta_i^{-1} ( \alpha_i -p^{-1} -(1/2
-p^{-1})_+) \beta_i \kappa_i \ge \mn( \alpha -p^{-1} \e -(1/2
-p^{-1})_+ \e) (\kappa^J, \beta^J) \\
+ (\mn( \alpha -p^{-1} \e -(1/2 -p^{-1})_+ \e) +\epsilon)
(\kappa^{J^\prime}, \beta^{J^\prime})\\
= \mu (\kappa, \beta) +\epsilon (\kappa^{J^\prime},
\beta^{J^\prime}) > \mu (r +j -1) +\epsilon (\kappa^{J^\prime},
\beta^{J^\prime}),
\end{multline*}
$$
n_\kappa = [c_0 2^{r -\gamma j -\gamma^\prime (\kappa^{J^\prime}, \beta^{J^\prime})}] +1 \ge
c_0 2^{r -\gamma j -\gamma^\prime (\kappa^{J^\prime}, \beta^{J^\prime})},
$$
\begin{multline*}
\log \frac{\mathfrak R_\kappa^{d,l-\e}} {n_\kappa} \le \log \frac{c_0
2^{(\kappa, \e)}} {c_0 2^{r -\gamma j -\gamma^\prime
(\kappa^{J^\prime}, \beta^{J^\prime})}} \le \log \frac{
2^{(\kappa, \beta)}} { 2^{r -\gamma j -\gamma^\prime
(\kappa^{J^\prime}, \beta^{J^\prime})}} \\
\le \log \frac{ 2^{r +j}} { 2^{r -\gamma j -\gamma^\prime
(\kappa^{J^\prime}, \beta^{J^\prime})}} = (\log 2) (j +\gamma j
+\gamma^\prime (\kappa^{J^\prime}, \beta^{J^\prime}))\ (\text{
см.} (1.4.10)),
\end{multline*}
и рассмотрим два случая. В первом случае, когда $ \theta > 2, $
используя неравенство Гельдера с показателем $ \theta /2 > 1, $
неравенства (3.3.8), (3.3.27), (1.2.4), для $ f \in C $ имеем
\begin{multline*}
(\sum_{\kappa \in \Z_+^d: r <
(\kappa, \beta) \le r +j_0, n_\kappa < \mathfrak R_\kappa^{d,l-\e}}
(2^{-(\kappa, \alpha -p^{-1} \e -(1/2 -p^{-1})_+ \e)}
n_\kappa^{-1/2} \times\\
(1 +\log \frac{\mathfrak R_\kappa^{d,l-\e}}
{n_\kappa})^{3/2}
2^{(\kappa, \alpha)} \Omega^{l \chi_{\s(\kappa)}}(f,
c_4 (2^{-\kappa})^{\s(\kappa)})_{L_p(I^d)})^2 )^{1/2} \le \\
(\sum_{\kappa \in \Z_+^d: r <
(\kappa, \beta) \le r +j_0, n_\kappa < \mathfrak R_\kappa^{d,l-\e}}
(2^{-(\kappa, \alpha -p^{-1} \e -(1/2 -p^{-1})_+ \e)}
n_\kappa^{-1/2} \times\\
(1 +\log \frac{\mathfrak R_\kappa^{d,l-\e}}
{n_\kappa})^{3/2} )^{2\theta /(\theta -2)} )^{1/2 -1/\theta} \times \\
(\sum_{\kappa \in \Z_+^d: r <
(\kappa, \beta) \le r +j_0, n_\kappa < \mathfrak R_\kappa^{d,l-\e}}
(2^{(\kappa, \alpha)} \Omega^{l \chi_{\s(\kappa)}}(f,
c_4 (2^{-\kappa})^{\s(\kappa)})_{L_p(I^d)})^\theta )^{1/\theta} \le \\
(\sum_{\kappa \in \Z_+^d: r <
(\kappa, \beta) \le r +j_0, n_\kappa < \mathfrak R_\kappa^{d,l-\e}}
(2^{-(\kappa, \alpha -p^{-1} \e -(1/2 -p^{-1})_+ \e)}
n_\kappa^{-1/2} \times\\
(1 +\log \frac{\mathfrak R_\kappa^{d,l-\e}}
{n_\kappa})^{3/2} )^{2\theta /(\theta -2)} )^{1/2 -1/\theta} \times \\
(\sum_{\kappa \in \Z_+^d: (\kappa, \beta) > r}
(2^{(\kappa, \alpha)} \Omega^{l \chi_{\s(\kappa)}}(f,
c_4 (2^{-\kappa})^{\s(\kappa)})_{L_p(I^d)})^\theta )^{1/\theta} \le \\
\le c_{30} (\sum_{j=1}^{j_0} \sum_{\kappa \in \Z_+^d: r +j -1 <
(\kappa, \beta) \le r +j, n_\kappa < \mathfrak R_\kappa^{d,l-\e}} (2^{-\mu (r
+j -1) -\epsilon (\kappa^{J^\prime}, \beta^{J^\prime})}
\times\\
(c_0 2^{r -\gamma j -\gamma^\prime (\kappa^{J^\prime},
\beta^{J^\prime})})^{-1/2} (j +\gamma j +\gamma^\prime
(\kappa^{J^\prime}, \beta^{J^\prime}))^{3/2} )^{2\theta /(\theta
-2)} )^{1/2 -1/\theta}\le
\end{multline*}
\begin{multline*} \tag{3.3.28} \le c_{31} ( \sum_{j=1}^{j_0} \sum_{\kappa \in
\Z_+^d: r +j -1 < (\kappa, \beta) \le r +j} (2^{-(\mu +1/2) r}
2^{-(\mu -\gamma /2)
j} j^{3/2} (1 +(\kappa^{J^\prime}, \beta^{J^\prime}))^{3/2}\\
\times 2^{-(\epsilon -\gamma^\prime /2) (\kappa^{J^\prime},
\beta^{J^\prime})} )^{2\theta/(\theta -2)} )^{1/2 -1/\theta} \\
= c_{31} 2^{-(\mu +1/2) r} (\sum_{j=1}^{j_0} (2^{-(\mu -\gamma /2)
j} j^{3/2})^{2\theta /(\theta -2)}\\
\times \sum_{\kappa \in \Z_+^d: r +j -1 < (\kappa, \beta) \le r
+j} (2^{-(\epsilon -\gamma^\prime /2) (\kappa^{J^\prime},
\beta^{J^\prime})} (1 +(\kappa^{J^\prime},
\beta^{J^\prime}))^{3/2})^{2\theta /(\theta -2)} )^{1/2 -1/\theta} \\
\le c_{32} 2^{-(\mu +1/2) r} (\sum_{j=1}^{j_0} (2^{-(\mu -\gamma /2)
j} j^{3/2})^{2\theta /(\theta -2)} (r +j)^{\cmn -1} )^{1/2 -1/\theta} \\
\le c_{32} 2^{-(\mu +1/2) r} (\sum_{j=1}^{j_0} (2^{-(\mu -\gamma /2)
j} j^{3/2})^{2\theta /(\theta -2)} r^{\cmn -1}
(1 +j)^{\cmn -1} )^{1/2 -1/\theta} \\
\le c_{32} 2^{-(\mu +1/2) r} r^{(\cmn -1)(1/2 -1/\theta)}
(\sum_{j=1}^\infty (2^{-(\mu -\gamma /2)
j} j^{3/2})^{2\theta /(\theta -2)}
(1 +j)^{\cmn -1} )^{1/2 -1/\theta} \\
 = c_{33}
2^{-\mn(\alpha -p^{-1} \e -(1/2 -p^{-1})_+ \e +(1/2) \e) r}
r^{(\cmn -1)(1/2 -1/\theta)} \\
= c_{33} 2^{-\mn(\alpha -(p^{-1}
-1/2)_+ \e ) r} r^{(\cmn -1)(1/2 -1/\theta)} = c_{33}
2^{-\mn(\alpha -(p^{-1} -1/2)_+ \e ) r} r^{(\cmn -1)(1/2
-1/\theta)_+}.
\end{multline*}

Во втором случае, когда $ \theta \le 2, $ применяя неравенство
Гельдера с показетелем $ \theta /2 \le 1, $ а также (3.3.27),
(3.3.8), для $ f \in C $ получаем
\begin{multline*} \tag{3.3.29}
(\sum_{\kappa \in \Z_+^d: r <
(\kappa, \beta) \le r +j_0, n_\kappa < \mathfrak R_\kappa^{d,l-\e}}
(2^{-(\kappa, \alpha -p^{-1} \e -(1/2 -p^{-1})_+ \e)}
n_\kappa^{-1/2} \times\\
(1 +\log \frac{\mathfrak R_\kappa^{d,l-\e}}
{n_\kappa})^{3/2}
2^{(\kappa, \alpha)} \Omega^{l \chi_{\s(\kappa)}}(f,
c_4 (2^{-\kappa})^{\s(\kappa)})_{L_p(I^d)})^2 )^{1/2} \le \\
(\sum_{\kappa \in \Z_+^d: r <
(\kappa, \beta) \le r +j_0, n_\kappa < \mathfrak R_\kappa^{d,l-\e}}
(2^{-(\kappa, \alpha -p^{-1} \e -(1/2 -p^{-1})_+ \e)}
n_\kappa^{-1/2} \times\\
(1 +\log \frac{\mathfrak R_\kappa^{d,l-\e}}
{n_\kappa})^{3/2} \times \\
2^{(\kappa, \alpha)} \Omega^{l \chi_{\s(\kappa)}}(f,
c_4 (2^{-\kappa})^{\s(\kappa)})_{L_p(I^d)})^\theta )^{1/\theta} \le \\
\le c_{34} (\sum_{j=1}^{j_0} \sum_{\kappa \in \Z_+^d: r +j -1 <
(\kappa, \beta) \le r +j, n_\kappa < \mathfrak R_\kappa^{d,l-\e}} (2^{-\mu (r
+j -1) -\epsilon (\kappa^{J^\prime}, \beta^{J^\prime})} \times \\
(c_0 2^{r -\gamma j -\gamma^\prime (\kappa^{J^\prime},
\beta^{J^\prime})})^{-1/2} (j +\gamma j +\gamma^\prime
(\kappa^{J^\prime}, \beta^{J^\prime}))^{3/2}
2^{(\kappa, \alpha)} \Omega^{l \chi_{\s(\kappa)}}(f,
c_4 (2^{-\kappa})^{\s(\kappa)})_{L_p(I^d)})^\theta )^{1/\theta} \le \\
\le c_{35} ( \sum_{j=1}^{j_0} \sum_{\kappa \in \Z_+^d: r +j -1 <
(\kappa, \beta) \le r +j} (2^{-(\mu +1/2) r} 2^{-(\mu -\gamma /2)
j} j^{3/2} (1 +(\kappa^{J^\prime}, \beta^{J^\prime}))^{3/2}\\
\times 2^{-(\epsilon -\gamma^\prime /2) (\kappa^{J^\prime},
\beta^{J^\prime})} 2^{(\kappa, \alpha)} \Omega^{l
\chi_{\s(\kappa)}}(f,
c_4 (2^{-\kappa})^{\s(\kappa)})_{L_p(I^d)})^\theta )^{1/\theta} = \\
= c_{35} 2^{-(\mu +1/2) r} (\sum_{j=1}^{j_0} (2^{-(\mu -\gamma /2)
j} j^{3/2})^\theta\\ \times \sum_{\kappa \in \Z_+^d: r +j -1 <
(\kappa, \beta) \le r +j} (2^{-(\epsilon -\gamma^\prime /2)
(\kappa^{J^\prime}, \beta^{J^\prime})} (1 +(\kappa^{J^\prime},
\beta^{J^\prime}))^{3/2}\\
\times 2^{(\kappa, \alpha)} \Omega^{l \chi_{\s(\kappa)}}(f,
c_4 (2^{-\kappa})^{\s(\kappa)})_{L_p(I^d)})^\theta )^{1/\theta} \le \\
c_{36} 2^{-(\mu +1/2) r} (\sum_{j=1}^{j_0}
\sum_{\kappa \in \Z_+^d: r +j -1 < (\kappa, \beta) \le r +j}
(2^{(\kappa, \alpha)} \Omega^{l \chi_{\s(\kappa)}}(f,
c_4 (2^{-\kappa})^{\s(\kappa)})_{L_p(I^d)})^\theta )^{1/\theta} \le \\
c_{36} 2^{-(\mu +1/2) r}
(\sum_{\kappa \in \Z_+^d: (\kappa, \beta) > r }
(2^{(\kappa, \alpha)} \Omega^{l \chi_{\s(\kappa)}}(f,
c_4 (2^{-\kappa})^{\s(\kappa)})_{L_p(I^d)})^\theta )^{1/\theta} \le \\
c_{37} 2^{-(\mu +1/2) r} = c_{37} 2^{-\mn(\alpha -p^{-1} \e -(1/2
-p^{-1})_+ \e +(1/2) \e) r} = c_{37} 2^{-\mn(\alpha -(p^{-1}
-1/2)_+ \e ) r}\\ = c_{37} 2^{-\mn(\alpha -(p^{-1} -1/2)_+ \e ) r}
r^{(\cmn -1)(1/2 -1/\theta)_+}.
\end{multline*}

Оценивая далее правую часть (3.3.26), выводим
\begin{multline*} \tag{3.3.30}
2^{-(r+j_0) \mn(\alpha -(p^{-1} -q^{-1})_+ \e)} (r +j_0)^{(\cmn
-1)(1/2 -1/\theta)_+}\\
 \le 2^{-(r +[\frac{r}{3\gamma}] )
\mn(\alpha -(p^{-1} -q^{-1})_+ \e)}
(r +[\frac{r}{3\gamma}] )^{(\cmn -1)(1/2 -1/\theta)_+} \\
\le 2^{-(r +\frac{r}{3\gamma} -1) \mn(\alpha -(p^{-1} -q^{-1})_+
\e)} (r +\frac{r}{3\gamma} )^{(\cmn -1)(1/2 -1/\theta)_+}\\
 \le c_{38} 2^{-\mn(\alpha -(p^{-1} -1/2)_+ \e) r} r^{(\cmn -1)(1/2
-1/\theta)_+}.
\end{multline*}

Соединяя (3.3.26), (3.3.28), (3.3.29),  (3.3.30) и учитывая, что
$ \dim M \le n, $ приходим к неравенству
$$
d_n(C,X) \le \sup_{ f \in C} \inf{ g \in M} \| f -g\|_X \le
c_{39} 2^{-\mn(\alpha -(p^{-1} -1/2)_+ \e) r}
r^{(\cmn -1)(1/2 -1/\theta)_+}.
$$
Отсюда и из (3.3.23) вытекает (3.3.19) во втором случае.

Для получения (3.3.19') заметим, что при
$ \alpha -(1/p -1/2)_+ \e >0 $ выполняется неравенство (см. [5])
\begin{multline*}
d_n(B( S_{p,\theta}^\alpha B(I^d)), L_2(I^d))\\ \ge c_{40}
n^{-\mn(\alpha -(p^{-1} -1/2)_+ \e)} (\log n)^{(\mn(\alpha
-(p^{-1} -1/2)_+ \e) +(1/2 -1/\theta)_+) (\cmn -1)}, \\
\text{где} B( S_{p,\theta}^\alpha B(I^d)) = B(L_p(I^d)) \cap
(\mathcal S_{p,\theta}^\alpha \mathcal B(I^d)).
\end{multline*}
Отсюда, пользуясь тем, что для $ f \in
L_q(I^d) $ при $ q \ge 2 $ соблюдается неравенство $ \|f
\|{L_2(I^d)} \le \|f\|_{L_q(I^d)}, $ в силу (3.2.1) при $ q \ge 2
$ и $ \alpha -(1/p -1/q)_+ \e >0 $ выводим
\begin{multline*}
d_n(\mathcal S_{p,\theta}^\alpha \mathcal B(I^d), L_q(I^d)) \ge
d_n(\mathcal S_{p,\theta}^\alpha \mathcal B(I^d), L_2(I^d)) \ge
d_n(B( S_{p,\theta}^\alpha B(I^d)), L_2(I^d))\\ \ge c_{40}
n^{-\mn(\alpha -(p^{-1} -1/2)_+ \e)} (\log n)^{(\mn(\alpha
-(p^{-1} -1/2)_+ \e) +(1/2 -1/\theta)_+) (\cmn -1)},
\end{multline*}
что совпадает с (3.3.19'). $ \square $

\end{document}